# Emergence of regularity for limit points of McKean–Vlasov particle systems

By **Robert Alexander Crowell** *of ETH Zürich*

**Abstract:** The empirical measure of an interacting $n$-particle system is a purely atomic random probability measure. In the limit as $n \to \infty$, we show for McKean–Vlasov systems with common noise that this measure becomes absolutely continuous with respect to Lebesgue measure for almost all times, almost surely. The density possesses good regularity properties, and we obtain a quantitative moment bound for its (random) fractional Sobolev norm. This result is obtained for dynamics with a bounded drift, bounded and Hölder-continuous diffusion coefficients and when the diffusion coefficient for the idiosyncratic noise is uniformly elliptic.

We directly study the sequence of particle systems via approximating their exchangeable dynamics by conditionally independent dynamics at the expense of an error. By using probabilistic means, the approximation and the error are controlled for each particle system, and the results are then passed to the large-system limit. The estimates thus obtained are combined via an analytic interpolation technique to derive a norm bound for the limiting random density. In this way, we obtain a regularity estimate for all cluster points, without requiring any knowledge of the dynamics of the limiting measure-valued flow.

**Robert Alexander Crowell**

rc@math.ethz.ch

Department of Mathematics
ETH Zürich
Rämistrasse 101
CH-8092 Zürich

**Acknowledgements:** This work is part of my PhD thesis written at ETH Zürich. I am deeply indebted to my advisor Martin Schweizer who gave me the freedom to pursue this line of research, and at the same time guided me with his immense experience throughout. I benefited greatly from his detailed feedback, and the arguments and the exposition in this manuscript improved markedly from his careful comments.

I would also like to thank participants at workshops held at ETH Zürich, Hammamet, CIRM Luminy, Oxford University and Imperial College London for their feedback.

Numerous discussions with my colleagues, in particular with my office-mate Florian Krach and with Andrew Allan, should not go unmentioned. They sharpened my reasoning.

**CONTENTS**

# INTRODUCTION

This paper is the first in a series considering the limiting behavior of particle systems of McKean–Vlasov type with common noise when the coefficients have only low regularity. We develop a regularity result for cluster points of the sequence $(\mu^n)_{n\in\mathbb{N}}$ of *empirical measure flows*, where $\mu^n = (\mu^n_t)_{t\in[0,T]}$ is associated with the $\mathbb{R}^d$-valued *n-particle system* given for $i \in \{1, \ldots, n\}$ and $t \in [0, T]$ by

$$\begin{cases} X^{i,n}_t = x_0 + \int_0^t b_s(X^{i,n}_s, \mu^n_s)\,\mathrm{d}s + \int_0^t \sigma_s(X^{i,n}_s, \mu^n_s)\,\mathrm{d}B^{i,n}_s + \int_0^t \bar{\sigma}_s(X^{i,n}_s, \mu^n_s)\,\mathrm{d}Z^n_s, \\ \mu^n_t = \frac{1}{n}\sum_{i=1}^n \delta_{X^{i,n}_t}. \end{cases} \quad (1)$$

The system (1) possesses a unique weak solution under very general assumptions on the coefficients. In addition, the process $\mu^n$ has continuous trajectories in the space $\mathcal{S}'(\mathbb{R}^d)$ of tempered distributions and is valued in the subset $\mathbf{M}_1^+(\mathbb{R}^d) \subseteq \mathcal{S}'(\mathbb{R}^d)$ of Borel probability measures on $\mathbb{R}^d$. Instead of the sequence of processes $(\mu^n)_{n\in\mathbb{N}}$, we study their laws $(\mathbf{P}^n_{x_0})_{n\in\mathbb{N}}$ together with a fixed $\mathcal{S}'(\mathbb{R}^d)$-valued process $\Lambda = (\Lambda_t)_{t\in[0,T]}$ satisfying $\mathbf{P}^n_{x_0} \circ \Lambda^{-1} = \mathrm{Law}(\mu^n)$. The sequence $(\mathbf{P}^n_{x_0})_{n\in\mathbb{N}}$ is tight and thus possesses cluster points.

Our main result in Theorem 1.5 establishes that *any cluster point* $\mathbf{P}^\infty_{x_0}$ of $(\mathbf{P}^n_{x_0})_{n\in\mathbb{N}}$ concentrates all its mass on flows of probability measures that are absolutely continuous with respect to Lebesgue measure, with a density that is in the Bessel potential or fractional Sobolev space $\mathsf{H}^w_r(\mathbb{R}^d)$ for some integrability $r > 1$ and regularity $w > 0$. We obtain this result under comparatively low regularity assumptions on the coefficients, namely, when $b$ is bounded, $\sigma, \bar{\sigma}$ are bounded and Hölder-continuous in an appropriate sense, and $\sigma$ is uniformly elliptic; see Assumptions 1.1 and 1.2 below. Qualitatively, Theorem 1.5 shows that

$$\mathbf{P}^\infty_{x_0}\left[\frac{\mathrm{d}\Lambda_t}{\mathrm{d}x} \text{ is in } \mathsf{H}^w_r(\mathbb{R}^d) \text{ for almost all } t \in (0, T]\right] = 1, \quad (2)$$

which is why we say that *regularity emerges* in the $n \to \infty$ limit. Indeed, as $\mu^n$ in (1) is purely atomic, we have that $\mathbf{P}^n_{x_0}[\mathrm{d}\Lambda_t/\mathrm{d}x$ exists for some $t \in (0, T]] = 0$ for all $n \in \mathbb{N}$. Hence each element of the sequence $(\mathbf{P}^n_{x_0})_{n\in\mathbb{N}}$ charges measure flows possessing no analytic regularity. However, (2) shows that the limit $\mathbf{P}^\infty_{x_0}$ of any convergent subsequence $(\mathbf{P}^{n_k}_{x_0})_{k\in\mathbb{N}}$ concentrates all its mass on measure flows possessing some analytic regularity. In fact, more can be said. If the system (1) includes the common noise $Z^n$, i.e. if $\bar{\sigma} \not\equiv 0$, then $\Lambda_t$ is a random probability measure under $\mathbf{P}^\infty_{x_0}$, and its density $p(t) = \mathrm{d}\Lambda_t/\mathrm{d}x$ is a random function in $\mathsf{H}^w_r(\mathbb{R}^d)$. Quantitatively, Theorem 1.5 gives a precise estimate for the regularity guaranteed by (2): We show that there exist a constant $c < \infty$ and appropriate parameters $r > 1$, $w > 0$, $q > 1$ and $\gamma < 1$ such that

$$\left\|\|p(t)\|_{\mathsf{H}^w_r(\mathbb{R}^d)}\right\|_{\mathbb{L}^q(\mathbf{P}^\infty_{x_0})} \leq c(1 \wedge t)^{-\gamma} \qquad \text{for all } t \in (0, T]. \quad (3)$$

The estimate (3) for cluster points of $(\mathbf{P}^n_{x_0})_{n\in\mathbb{N}}$ under such general assumptions is new. To understand why it is difficult to obtain and to appreciate why it is useful, it is worthwhile to explain how (3) can be proved in the case of very regular coefficients, e.g. if $b, \sigma, \bar{\sigma}$ are smooth. In that case, the classical *propagation of chaos* property establishes that as $n \to \infty$, the pair $(X^{1,n}, \mu^n)$ converges in an appropriate sense to a unique process $(X, U)$ solving the *McKean–Vlasov SDE* given for $t \in [0, T]$ by

$$\begin{cases} X_t = x_0 + \int_0^t b_s(X_s, U_s) \, \mathrm{d}s + \int_0^t \sigma_s(X_s, U_s) \, \mathrm{d}B_s + \int_0^t \bar{\sigma}_s(X_s, U_s) \, \mathrm{d}Z_s, \\ U_t = \mathrm{Law}(X_t \,|\, \mathcal{F}_t^Z), \end{cases} \quad (4)$$

so that (4) can be viewed as a large-system limiting model of (1). This provides an *indirect way* to obtain (3) by analyzing instead of the limit $U$ of the sequence $(\mu^n)_{n \in \mathbb{N}}$ from (1) the flow $(\mathrm{Law}(X_t \,|\, \mathcal{F}_t^Z))_{t \in [0,T]}$ in (4) of conditional time-marginal laws. This agenda is carried out in Crisan and McMurray [15], where Malliavin calculus is used to deduce that if $\sigma$ is in addition uniformly elliptic, then $\mathrm{Law}(X_t \,|\, \mathcal{F}_t^Z)$ possesses a smooth density relative to Lebesgue measure for all $t > 0$, almost surely. Thus with regular coefficients, the emergence of regularity can be derived from the propagation of chaos.

Under the weaker regularity assumptions we impose here, the propagation of chaos need not hold; so we must work directly with the sequence $(\mathbf{P}_{x_0}^n)_{n \in \mathbb{N}}$. As a consequence, the result in (3) is conceptually different from an estimate derived via (4). Our result is an a priori regularity estimate for $\Lambda$ under *any* of the potentially many cluster points of $(\mathbf{P}_{x_0}^n)_{n \in \mathbb{N}}$, and this can be used derive a propagation of chaos property from the emergence of regularity, i.e., to reverse the above reasoning. While we leave the details of this application the follow-up papers [18] and [17] based on the author's PhD thesis [19], we offer an appetizer to help appreciate the usefulness of (3). Cluster points of $(\mathbf{P}_{x_0}^n)_{n \in \mathbb{N}}$ lead to natural solution candidates for the stochastic Fokker–Planck equation, which is a certain SPDE associated to (4). Verifying that a candidate is in fact a solution is classical if the coefficients are regular (see e.g. Vaillancourt [53]); but if they are not, it is beyond established means. The estimate (3) is a crucial part to close this gap when the drift $b$ is irregular. Our strategy in that forthcoming work is to compensate the low regularity of $b$ with the emerging regularity of $\Lambda$. Once this is achieved, (4) can be solved and given meaning as a large-system limit of (1).

To obtain the key result (3), we devise a new *direct strategy* to work with the sequence $(\mathbf{P}_{x_0}^n)_{n \in \mathbb{N}}$ via an *approximation and interpolation* scheme. This has its foundations in the seminal contribution of Fournier and Printems [25], with subsequent works such as Debussche and Romito [21] or Bally and Caramellino [5]; see Section 7.2 for a fuller discussion and additional references. Adapting this scheme to the particle systems (1) first requires careful adjustments because the common noise $Z^n$ renders $\Lambda$ under $\mathbf{P}_{x_0}^\infty$ random, and second also new ideas. In broad terms, we start with an Euler-type approximation of the particle systems (1), and use this to define an approximated empirical measure. By jointly passing the empirical measure and its approximation to narrow limits, we can decompose $(\Lambda_t)_{t \in [0,T]}$ into a sum of two-parameter processes $(A_{\varepsilon,t} + E_{\varepsilon,t})_{0 \leq \varepsilon \leq t \leq T}$ of an approximation and an error, respectively. To derive regularity estimates for the terms $A$ and $E$, we need probabilistic arguments building on the particle systems and a limiting procedure. The nature of the Euler-type approximation of (1) lets us represent $A_{\varepsilon,t}$ as a Gaussian mixture law, and for elliptic $\sigma$ and $\varepsilon > 0$, this measure has a smooth density. However, as $\varepsilon \to 0$, it becomes singular, but the rate at which it does can be bounded. The error $E_{\varepsilon,t} = \Lambda_t - A_{\varepsilon,t}$, on the other hand, is a signed measure. To analyze its behavior, we use a Kolmogorov-type regularity argument at the level of the particle systems that we then pass to narrow limits, to find the rate at which $E_{\varepsilon,t}$ vanishes as $\varepsilon \to 0$. Here we need the Hölder-regularity of $\sigma$, $\bar{\sigma}$. To obtain an estimate for $\Lambda_t$, we then apply an interpolation argument from harmonic analysis to trade off the rate at which the smooth part $A_{\varepsilon,t}$ becomes singular against the rate at which $E_{\varepsilon,t}$ vanishes as $\varepsilon \to 0$. However, the randomness from the common noise makes the estimates we have random, available only for certain pairs $(\varepsilon, t)$, and possibly fail on nullsets. This restricts the way we can combine the estimates. We thus build in Proposition 5.1 a work-horse interpolation result to overcome this difficulty. While this is in the spirit of Debussche and Romito [21] or Bally and Caramellino [5], it is more subtle and must be tailored specifically to meet our needs.

The approach and tools we develop are not limited to McKean–Vlasov dynamics, although this much-studied case is arguably interesting and well-suited to demonstrate our ideas. We expect that our strategy to obtain (3) can be flexible enough to study other particle systems in a similar fashion; e.g. Kurtz and Xiong [35] or Crisan et al. [16]. A comparison to the literature and a further outlook are given at the end in Section 7.

**Organization**

This paper is structured as follows. The probabilistic setup, the main result in Theorem 1.5 and the outline of its proof are contained in Section 1. Section 2 constructs in (1.16) the decomposition of $\Lambda_t$ into $A_{\varepsilon,t}$ and $E_{\varepsilon,t}$. Sections 3 and 4 derive in Propositions 3.2 and 4.1, respectively, bounds for $A_{\varepsilon,t}$ and $E_{\varepsilon,t}$. In Section 5, we present the abstract interpolation result in Proposition 5.1. Section 6 combines all these ingredients into the proof of the main result in Theorem 1.5. Finally, Section 7 concludes with a discussion and an outlook. The Appendices A and B contain useful results from functional analysis and probability, as well as proofs omitted from the main body.

**Notation**

We collect the essential notation. Appendix A contains more background and references.

- Throughout $c(a_1, \ldots, a_k)$ denotes a generic finite constant which depends on the quantities $a_1, \ldots, a_k$, and which may change from line to line.

- The space of *Schwartz functions* is denoted by $\mathcal{S} = \mathcal{S}(\mathbb{R}^d)$ and endowed with its usual Fréchet topology. We write $\mathcal{S}' = \mathcal{S}'(\mathbb{R}^d)$ for the topological vector space of *(tempered) distributions*, and $\lambda[\phi] := \langle \lambda; \phi \rangle_{\mathcal{S}' \times \mathcal{S}}$ for the duality pairing of $\lambda \in \mathcal{S}'$ and $\phi \in \mathcal{S}$. The Fourier transform $\mathscr{F}$ is a linear automorphism in $\mathcal{S}'$ with inverse $\mathscr{F}^{-1}$.

- For a measure-space $(X, \mathcal{X}, \chi)$ and a Banach space $(B, \|\cdot\|_B)$, we denote the associated *Bochner space* by $\mathbb{L}^p(\chi; B) = \mathbb{L}^p((X, \mathcal{X}, \chi); B)$. If $B = \mathbb{R}$, we write $\mathbb{L}^p(\chi) = \mathbb{L}^p(X, \mathcal{X}, \chi)$ for the usual *Lebesgue space*.

- The *Bessel potential space* with integrability $r \in (1, \infty)$ and regularity $s \in \mathbb{R}$ is defined by $\mathsf{H}^s_r = \mathsf{H}^s_r(\mathbb{R}^d) := \{f \in \mathcal{S}'(\mathbb{R}^d) : \|f\|_{\mathsf{H}^s_r(\mathbb{R}^d)} < \infty\}$, where $\|f\|_{\mathsf{H}^s_r(\mathbb{R}^d)} := \|J^s f\|_{\mathbb{L}^r(\mathrm{d}x)}$ with $J^s f = \mathscr{F}^{-1}(h^s(\mathscr{F}f))$ for $f \in \mathcal{S}'$ and $h(\xi) := (1 + |\xi|^2)^{1/2}$.

- If $(X, \mathcal{X})$ is a measurable space then the set of *probability measures* is denoted by $\mathbf{M}_1^+(X) = \mathbf{M}_1^+(X, \mathcal{X})$. For $\mu, \nu \in \mathbf{M}_1^+(X)$ with $\int_X |x| \mathrm{d}\mu, \int_X |x| \mathrm{d}\nu < \infty$, we define the *Kantorovich–Rubinstein* distance by $\mathsf{d}_{\mathbb{W}_1}(\mu, \nu) := \sup\{\int f\,\mathrm{d}(\mu - \nu) : \|f\|_{\mathsf{Lip}} \leq 1\}$. We write $\mathcal{P}_{\mathrm{wk}^*} = \mathcal{P}_{\mathrm{wk}^*}(\mathbb{R}^d)$ for the set $\mathbf{M}_1^+ = \mathbf{M}_1^+(\mathbb{R}^d)$ with the *narrow topology*.

- Note that $\mathbf{M}_1^+ \subseteq \mathcal{S}'$ and that $\mathcal{P}_{\mathrm{wk}^*} \hookrightarrow \mathcal{S}'$, which is to say that the embedding is continuous. In particular, if $\mu \in \mathbf{M}_1^+$ and $\phi \in \mathcal{S}$, then $\mu[\phi] = \langle \mu; \phi \rangle_{\mathcal{S}' \times \mathcal{S}} = \int_{\mathbb{R}^d} \phi(x)\,\mu(\mathrm{d}x)$.

# 1 SETUP AND MAIN RESULT

## 1.1 Finite particle systems and assumptions

Let $T < \infty$ be the terminal time and $d \in \mathbb{N}$ the dimension. The basis for our study is the family of *n-particle systems* of McKean–Vlasov type formalized for each $n \in \mathbb{N}$ by

$$\begin{cases} X_t^{i,n} = x_0 + \int_0^t b_s(X_s^{i,n}, \mu_s^n)\,\mathrm{d}s + \int_0^t \sigma_s(X_s^{i,n}, \mu_s^n)\,\mathrm{d}B_s^{i,n} + \int_0^t \bar{\sigma}_s(X_s^{i,n}, \mu_s^n)\,\mathrm{d}Z_s^n\,, \\ \mu_t^n = \frac{1}{n}\sum_{i=1}^n \delta_{X_t^{i,n}}\,, \qquad\qquad\qquad\qquad\qquad \text{for all } i = 1, \ldots, n \text{ and } t \in [0, T]. \end{cases} \quad (1.1)$$

The $n$-particle system (1.1) consists of $n$ coupled stochastic differential equations (SDEs), each with values in $\mathbb{R}^d$. We refer to $X^{i,n}$, the $i$'th component of the system, as a *particle*, to $B^{1,n}, \ldots, B^{n,n}$ as the *idiosyncratic noise* and to $Z^n$ as the *common noise*. To state the assumptions on the coefficients, we need on $[0,T] \times \mathbb{R}^d \times \mathbf{M}_1^+(\mathbb{R}^d)$ the product $\sigma$-algebra $\mathcal{E}$ generated by the Lebesgue-measurable sets of $[0,T] \times \mathbb{R}^d$ and the Borel-measurable sets of $\mathcal{P}_{\mathrm{wk}^*}(\mathbb{R}^d)$.

**Assumption 1.1** | The drift coefficient $b : [0,T] \times \mathbb{R}^d \times \mathbf{M}_1^+(\mathbb{R}^d) \to \mathbb{R}^d$ is bounded and measurable relative to $\mathcal{E}$.

**Assumption 1.2** | The diffusion coefficients $\sigma : [0,T] \times \mathbb{R}^d \times \mathbf{M}_1^+(\mathbb{R}^d) \to \mathbb{R}^{d \times d}$ as well as $\bar{\sigma} : [0,T] \times \mathbb{R}^d \times \mathbf{M}_1^+(\mathbb{R}^d) \to \mathbb{R}^{d \times m}$ are bounded and measurable relative to $\mathcal{E}$. In addition,

(E) there is $\kappa > 0$ such that the uniform ellipticity condition

$$z^\top \big((\sigma_t \sigma_t^\top)(x, \mu)\big) z \geq \kappa\,|z|^2$$

holds for all $z \in \mathbb{R}^d$, $t \in [0,T]$ and $(x, \mu) \in \mathbb{R}^d \times \mathbf{M}_1^+(\mathbb{R}^d)$;

(H) there exist $\beta > 0$ and a constant $C < \infty$ such that the Hölder-regularity conditions

$$|\sigma_t(x, \mu) - \sigma_t(x', \mu')| \leq C\big(|x - x'|^\beta + \mathsf{d}_{\mathbb{W}_1}(\mu, \mu')^\beta\big)\,,$$
$$|\bar{\sigma}_t(x, \mu) - \bar{\sigma}_{t'}(x', \mu')| \leq C\big(|t - t'|^\beta + |x - x'|^\beta + \mathsf{d}_{\mathbb{W}_1}(\mu, \mu')^\beta\big)\,,$$

hold for all $t, t' \in [0, T]$ and $(x, \mu), (x', \mu') \in \mathbb{R}^d \times \mathbf{M}_1^+(\mathbb{R}^d)$.

## 1.2 Probabilistic setup for the particle systems

The $n$-particle system (1.1) is solved in the usual weak sense of SDEs; cf. Stroock and Varadhan [49, Ch. 6]. To construct such a solution, we write

$$\Omega^n_{\text{par}} := \mathbf{C}([0,T]; \mathbb{R}^{nd} \times \mathbb{R}^m)$$

with the coordinate process representing the particles and the common noise, i.e.,

$$(t,\omega) \mapsto \omega_t =: (X^n_t, Z^n_t)(\omega) = (X^{1,n}_t, \ldots, X^{n,n}_t, Z^n_t)(\omega) \qquad \text{for } (t,\omega) \in [0,T] \times \Omega^n_{\text{par}}.$$

We endow $\Omega^n_{\text{par}}$ with the filtration $\mathbb{G}^n := (\mathcal{G}^n_t)_{0 \leq t \leq T}$ given by the right-continuous version of the canonical raw filtration, i.e.,

$$\mathcal{G}^n_t := \bigcap_{0 < \varepsilon < T-t} \sigma\Big((X^n_s, Z^n_s) : 0 \leq s < t + \varepsilon\Big) \qquad \text{for } t \in [0,T),$$

$$\mathcal{G}^n := \mathcal{G}^n_T := \sigma\Big((X^n_s, Z^n_s) : 0 \leq s \leq T\Big).$$

**Lemma 1.3** | *Under Assumptions 1.1 and 1.2, there exist a unique probability measure $\mathbb{P}^n_{x_0}$ on $(\Omega^n_{\text{par}}, \mathcal{G}^n)$ and independent standard $d$-dimensional $(\mathbb{G}^n, \mathbb{P}^n_{x_0})$-Brownian motions $B^{1,n}, \ldots, B^{n,n}$ such that $Z^n$ is a standard $m$-dimensional $(\mathbb{G}^n, \mathbb{P}^n_{x_0})$-Brownian motion that is independent of $B^n := (B^{1,n}, \ldots, B^{n,n})$, and $(X^n, B^n, Z^n)$ solves (1.1). Moreover, the law of $X^n$ under $\mathbb{P}^n_{x_0}$ is exchangeable and we have*

$$\mathbf{E}^{\mathbb{P}^n_{x_0}}\left[\sup_{t \in [0,T]} |X^{i,n}_t|^q\right] =: c_q < \infty \qquad \text{for any } q > 1, \tag{1.2}$$

*where the constant $c_q$ depends only on $q$, $T$ and the bounds for $b$, $\sigma$, $\bar{\sigma}$, but not on $n$.*

By exchangeability, we mean that for any permutation $\pi$ of $[n] := \{1, \ldots, n\}$, we have

$$\text{Law}_{\mathbb{P}^n_{x_0}}(X^{1,n}, \ldots, X^{n,n}) = \text{Law}_{\mathbb{P}^n_{x_0}}(X^{\pi(1),n}, \ldots, X^{\pi(n),n}).$$

In the sequel, we denote the weak solution from Lemma 1.3 by

$$\boldsymbol{X}^n = \Big(\Omega^n_{\text{par}}, \mathcal{G}^n, \mathbb{G}^n, \mathbb{P}^n_{x_0}, (X^n_t)_{t \in [0,T]}, (B^n_t)_{t \in [0,T]}, (Z^n_t)_{t \in [0,T]}\Big). \tag{1.3}$$

Note that Assumption 1.2(E) also allows $\bar{\sigma} = 0$, in which case there is no common noise $Z^n$ in (1.1) and the above simplifies in the obvious way.

***Proof of Lemma 1.3*** Under Assumptions 1.1 and 1.2, existence and uniqueness in law follow from standard results on elliptic SDEs with continuous diffusion coefficients $\sigma$ and $\bar{\sigma}$, see Karatzas and Shreve [33, Sec. 5.4.D], paired with a Girsanov change of measure [33, Sec. 5.3.B], which is well defined since the drift coefficient $b$ is uniformly bounded and $\sigma$ is uniformly elliptic so $\sigma^{-1}$ is bounded; cf. (B.2). Exchangeability of the solution is shown in Vaillancourt [53, Thm. 2]. Standard estimates for SDEs with bounded coefficients show that for any $q > 1$, there exists $c_q < \infty$ satisfying $\mathbf{E}^{\mathbb{P}^n_{x_0}}[\sup_{t \in [0,T]} |X^{i,n}_t|^q] \leq c_q$ for all $n \in \mathbb{N}$ and $i \in [n]$. □

## 1.3 Probabilistic setup for the measure flow

From the $n$-particle system (1.1) we have the *empirical measure flow* $\mu^n := (\mu^n_t)_{t \in [0,T]}$. The sequence $(\mu^n)_{n \in \mathbb{N}}$ of these processes is a key object in this work. To focus on it, we transfer it to a canonical setup. For this, let

$$\Omega := \mathbf{C}\Big([0,T]; \mathcal{S}'(\mathbb{R}^d)\Big)$$

denote the *canonical trajectory space*. Its topological vector space structure is induced by convergence uniformly on $[0,T]$ for bounded sets in $\mathcal{S}(\mathbb{R}^d)$; see e.g. Kallianpur and Xiong [32, Ch. 2.4] or Appendix A. We denote the canonical process on $\Omega$ by

$$(t,\omega) \mapsto \omega_t =: \Lambda_t(\omega) = \Lambda_t \qquad \text{for } (t,\omega) \in [0,T] \times \Omega,$$

We endow $\Omega$ with the right-continuous version $\mathbb{F} := (\mathcal{F}_t)_{0 \leq t \leq T}$ of the the canonical filtration generated by $\Lambda$, i.e.,

$$\mathcal{F}_t := \bigcap_{0 < \varepsilon < T-t} \sigma(\Lambda_s : 0 \leq s < t + \varepsilon) \qquad \text{for } t \in [0,T),$$
$$\mathcal{F} := \mathcal{F}_T := \sigma(\Lambda_s : 0 \leq s \leq T).$$

The set $\mathbf{M}_1^+(\mathbb{R}^d)$ is not closed in the topology of $\mathcal{S}'(\mathbb{R}^d)$; see Remark A.2. We therefore need to establish part 1) of the following technical result, the proof of which is of subordinate importance and can be skipped at first reading.

**Lemma 1.4** | *Under Assumptions 1.1 and 1.2, let $\mathbf{X}^n$ be the weak solution from Lemma 1.3 of the $n$-particle system (1.1). Then there exists a version of the process $\mu^n = (\mu^n_t)_{t \in [0,T]}$ which has continuous trajectories in the strong topology of $\mathcal{S}'(\mathbb{R}^d)$. In particular,*

$$\mathbf{P}^n_{x_0} := \mathbb{P}^n_{x_0} \circ (\mu^n)^{-1} \tag{1.4}$$

*is a well-defined probability measure on $(\Omega, \mathcal{F})$. Moreover:*
  1) *The sequence $(\mathbf{P}^n_{x_0})_{n \in \mathbb{N}}$ is tight and any narrow cluster point $\mathbf{P}^\infty_{x_0}$ of $(\mathbf{P}^n_{x_0})_{n \in \mathbb{N}}$ satisfies*

$$\mathbf{P}^\infty_{x_0}\Big[\,\mathbf{C}\big([0,T]; \mathcal{S}'(\mathbb{R}^d) \cap \mathbf{M}_1^+(\mathbb{R}^d)\big)\,\Big] = 1\,, \tag{1.5}$$

*that is, $\mathbf{P}^\infty_{x_0}$ concentrates its mass on the set of probability-measure-valued processes.*
  2) *For any $\varepsilon, K > 0$ and $q > 1$, we have for all $n \in \mathbb{N} \cup \{\infty\}$ the concentration bound*

$$\mathbf{P}^n_{x_0}\Big[\omega \in \Omega : \Lambda_t\big[-K,K\big]^d < 1 - \varepsilon \text{ for some } t \in [0,T]\Big] \leq \frac{c_q}{\varepsilon K^q}\,, \tag{1.6}$$

*where $c_q$ is the constant from (1.2), which depends on $q$ but not on $K$, $\varepsilon$ nor $\mathbf{P}^n_{x_0}$.*

**Proof** The fact that $\mu^n$ has continuous trajectories in $\mathcal{S}'$ follows by combining the result of Mitoma [39, 40] with the observation that for each $\phi \in \mathcal{S}$, the trajectories of the process $\mu^n[\phi] := (\mu^n_t[\phi])_{t \in [0,T]}$ with $\mu^n_t[\phi] = \frac{1}{n} \sum_{i=1}^n \phi(X^{i,n}_t)$ are continuous in $\mathbb{R}$.

For part 1), the tightness of $(\mathbf{P}^n_{x_0})_{n \in \mathbb{N}}$ follows from the more general result in Lemma 2.6 and in particular the estimate (2.27) and the property (2.32). The concentration of mass property (1.5) is deduced from the concentration inequality (1.6) which we establish first.

For part 2), for $\delta > 0$ and each $K \in \mathbb{N}$, let $\phi_{\delta,K} \in \mathbf{C}_c^\infty(\mathbb{R}^d)$ be a function satisfying $\mathbb{1}_{[-(K-\delta), K-\delta]^d} \leq \phi_{\delta,K} \leq \mathbb{1}_{[-K,K]^d}$. For any $n \in \mathbb{N}$, $\Lambda$ is $\mathbf{P}_{x_0}^n$-a.s. $\mathbf{M}_1^+(\mathbb{R}^d)$-valued, and for any positive measure $\lambda(\omega)$, we have $\lambda(\omega)[\phi_{\delta,K}] \leq \lambda(\omega)[[-K,K]^d]$, hence

$$\left\{\omega : \lambda(\omega)\big[[-K,K]^d\big] < 1 - \varepsilon\right\} \subseteq \left\{\omega : \lambda(\omega)[\phi_{\delta,K}] < 1 - \varepsilon\right\}. \tag{1.7}$$

With $F_{\varepsilon,K}$ defined as above, we thus get

$$\mathbf{P}_{x_0}^n[F_{\varepsilon,K}] \leq \mathbf{P}_{x_0}^n\Big[\Lambda_t[\phi_{\delta,K}] < 1 - \varepsilon \text{ for some } t \in [0,T]\Big] \tag{1.8}$$
$$\leq \mathbf{P}_{x_0}^n\left[\inf_{t \in [0,T]} \Lambda_t[\phi_{\delta,K}] < 1 - \varepsilon\right].$$

We claim that for each $n \in \mathbb{N}$,

$$\left\{\omega \in \Omega_{\text{par}}^n : \inf_{t \in [0,T]} \mu_t^n(\omega)[\phi_{\delta,K}] < 1 - \varepsilon\right\} \subseteq \left\{\omega \in \Omega_{\text{par}}^n : \frac{1}{n}\sum_{i=1}^n \mathbb{1}_{\{\sup_{t \in [0,T]} |X_t^{i,n}| \geq K - \delta\}} \geq \varepsilon\right\}.$$

Indeed, if $\mu_t^n[\phi_{\delta,K}] < 1 - \varepsilon$ for some $t \in [0,T]$, then a fraction of at least $\varepsilon$ of particles must have exited $[-(K-\delta), K-\delta]^d$ by time $t$, and so the inclusion follows. Therefore, using (1.4) to pass the estimate from the joint law $\mathbb{P}_{x_0}^n$ of the whole particle system to the representative particle model $\mathbf{P}_{x_0}^n$ and using Markov's inequality yields

$$\mathbf{P}_{x_0}^n\left[\inf_{t \in [0,T]} \Lambda_t[\phi_{\delta,K}] < 1 - \varepsilon\right] = \mathbb{P}_{x_0}^n\left[\inf_{t \in [0,T]} \mu_t^n[\phi_{\delta,K}] < 1 - \varepsilon\right]$$
$$\leq \mathbb{P}_{x_0}^n\left[\frac{1}{n}\sum_{i=1}^n \mathbb{1}_{\{\sup_{t \in [0,T]} |X_t^{i,n}| \geq K - \delta\}} \geq \varepsilon\right]$$
$$\leq \frac{1}{\varepsilon}\mathbf{E}^{\mathbb{P}_{x_0}^n}\left[\frac{1}{n}\sum_{i=1}^n \mathbb{1}_{\{\sup_{t \in [0,T]} |X_t^{i,n}| \geq K - \delta\}}\right]$$
$$= \frac{1}{n\varepsilon}\sum_{i=1}^n \mathbb{P}_{x_0}^n\left[\sup_{t \in [0,T]} |X_t^{i,n}| \geq K - \delta\right].$$

In the last line, Markov's inequality and exchangeability now lead for any $q > 1$ to

$$\frac{1}{n\varepsilon}\sum_{i=1}^n \mathbb{P}_{x_0}^n\left[\sup_{t \in [0,T]} |X_t^{i,n}| \geq K - \delta\right] \leq \frac{1}{n\varepsilon(K-\delta)^q}\sum_{i=1}^n \mathbf{E}^{\mathbb{P}_{x_0}^n}\left[\sup_{t \in [0,T]} |X_t^{i,n}|^q\right]$$
$$= \frac{1}{\varepsilon(K-\delta)^q}\mathbf{E}^{\mathbb{P}_{x_0}^n}\left[\sup_{t \in [0,T]} |X_t^{1,n}|^q\right].$$

From (1.2), we have $c_q = \mathbf{E}^{\mathbb{P}_{x_0}^n}[\sup_{t \in [0,T]} |X_t^{i,n}|^q] < \infty$, uniformly in $n \in \mathbb{N}$, and so

$$\sup_{n \in \mathbb{N}} \mathbf{P}_{x_0}^n\left[\inf_{t \in [0,T]} \Lambda_t[\phi_{\delta,K}] < 1 - \varepsilon\right] \leq \frac{c_q}{\varepsilon(K-\delta)^q}. \tag{1.9}$$

To pass this inequality to narrow limit points, note that for each fixed $\phi \in \mathcal{S}$, the evaluation map $\mathcal{S}' \ni \lambda \mapsto \langle \lambda; \phi \rangle_{\mathcal{S}' \times \mathcal{S}} = \lambda[\phi] \in \mathbb{R}$ is continuous for the strong topology of $\mathcal{S}'$. For each $t \in [0,T]$, the set $\{\Lambda_t[\phi_{\delta,K}] < 1 - \varepsilon\}$ is thus open in $\mathbf{C}([0,T]; \mathcal{S}')$, and hence so is the set $\{\inf_{t\in[0,T]} \Lambda_t[\phi_{\delta,K}] < 1-\varepsilon\} = \bigcup_{t\in\mathbb{Q}\cap[0,T]}\{\Lambda_t[\phi_{\delta,K}] < 1-\varepsilon\}$ because $t \mapsto \Lambda_t[\phi_{\delta,K}]$ is continuous. By the Portmanteau theorem and (1.9), we then have

$$\mathbf{P}_{x_0}^\infty\left[\inf_{t\in[0,T]} \Lambda_t[\phi_{\delta,K}] < 1-\varepsilon\right] \leq \frac{c_q}{\varepsilon(K-\delta)^q}, \quad (1.10)$$

and letting $\delta \to 0$ yields the right-hand side in (1.6). It remains to show that $\Lambda$ is also $\mathbf{P}_{x_0}^\infty$-a.s. $\mathsf{M}_1^+(\mathbb{R}^d)$-valued so that we can use (1.7) and hence also (1.8) with $n = \infty$. Once this is done, (1.6) follows from (1.8)–(1.10). So take $\phi \in \mathcal{S}$ with $\phi \geq 0$ and note that the set $\{\inf_{t\in[0,T]}\Lambda_t[\phi] \geq 0\} = \bigcap_{t\in[0,T]\cap\mathbb{Q}}\{\Lambda_t[\phi] \geq 0\}$ is closed. So the Portmanteau theorem shows that $\mathbf{P}_{x_0}^\infty$-a.s., $\Lambda_t$ is a positive measure simultaneously for all $t \in [0,T]$.

We can now complete the proof of part 1). Choose the functions $\phi_{\delta,K}$ to satisfy in addition $\phi_{\delta,K+1} \geq \phi_{\delta,K}$ for all $K \in \mathbb{N}$. Then (1.10) and the monotonicity of measures give

$$\mathbf{P}_{x_0}^\infty\left[\bigcap_{K\in\mathbb{N}}\left\{\inf_{t\in[0,T]}\Lambda_t[\phi_{\delta,K}] < 1-\varepsilon\right\}\right] = \lim_{K\to\infty} \mathbf{P}_{x_0}^\infty\left[\inf_{t\in[0,T]}\Lambda_t[\phi_{\delta,K}] < 1-\varepsilon\right] = 0$$

so $\phi_{\delta,K} \uparrow \mathbb{1}$ gives $\mathbf{P}_{x_0}^\infty[\inf_{t\in[0,T]}\Lambda_t[\mathbb{R}^d] < 1-\varepsilon] = 0$. Letting $\varepsilon \to 0$, we obtain (1.5), and the proof of part 1) is complete. $\square$

## 1.4 Main result

By Lemma 1.4, the canonical process $\Lambda$ on $(\Omega, \mathcal{F}, \mathbf{P}_{x_0}^n)$, with $\mathbf{P}_{x_0}^n$ defined via (1.4), represents the dynamics of the empirical measure in the sense that for all $n \in \mathbb{N}$, we have $\mathrm{Law}_{\mathbb{P}_{x_0}^n}(\mu^n) = \mathrm{Law}_{\mathbf{P}_{x_0}^n}(\Lambda)$. Accordingly, we term $\mathbf{P}_{x_0}^n$ the *empirical measure law*. Note the distinction between $\mathbb{P}_{x_0}^n$ and $\mathbf{P}_{x_0}^n$: the former denotes the joint law of $n$ particles and the common Brownian motion, the latter models the flow of the empirical measure of the particles. In particular, $\Lambda_t$ is a purely atomic measure for all $t \in [0,T]$, $\mathbf{P}_{x_0}^n$-a.s. Lemma 1.4 also shows that under any narrow cluster point $\mathbf{P}_{x_0}^\infty$ of $(\mathbf{P}_{x_0}^n)_{n\in\mathbb{N}}$, of which there exists at least one due to tightness, the process $\Lambda$ remains $\mathsf{M}_1^+(\mathbb{R}^d)$-valued.

We next state our main result. It shows in particular that under any cluster point $\mathbf{P}_{x_0}^\infty$, we have for almost all $t \in (0,T)$ that $\Lambda_t$ is absolutely continuous with respect to Lebesgue measure, $\mathbf{P}_{x_0}^\infty$-a.s.

**Theorem 1.5** | *Impose Assumptions 1.1 and 1.2. Then there exist real numbers $w > 0$ and $r > 1$ such that for any cluster point $\mathbf{P}_{x_0}^\infty$ of the sequence $(\mathbf{P}_{x_0}^n)_{n\in\mathbb{N}}$, we have*

$$\mathbf{P}_{x_0}^\infty\left[\left\{\omega \in \Omega : \frac{\mathrm{d}\Lambda_t(\omega)}{\mathrm{d}x} = \bar{p}(t,\omega) \text{ is in } \mathsf{H}_r^w(\mathbb{R}^d) \text{ for almost all } t \in (0,T]\right\}\right] = 1, \quad (1.11)$$

*where $\bar{p} : (0,T] \times \Omega \to \mathsf{H}_r^w(\mathbb{R}^d)$ is a strongly measurable function.*

*More precisely, for each $t \in (0,T]$, we have $\mathbf{P}_{x_0}^\infty$-a.s. that $\Lambda_t \ll \mathrm{d}x$, and for any $q \geq 1$, $w$ and $r$ can be chosen in such a way that there exists $1 > \gamma > 0$ such that the function $t \mapsto (\omega \mapsto \mathrm{d}\Lambda_t(\omega)/\mathrm{d}x)$ defines a strongly measurable map*

$$p : (0,T] \to \mathbb{L}^q\Big((\Omega, \mathcal{F}, \mathbf{P}_{x_0}^\infty); \mathsf{H}_r^w(\mathbb{R}^d)\Big) \quad (1.12)$$

*which satisfies for some constant $c_{\text{Thm. 1.5}} < \infty$ the bound*

$$\left\| \|p(t)\|_{\mathsf{H}_r^w(\mathbb{R}^d)} \right\|_{\mathbb{L}^q(\mathbf{P}_{x_0}^\infty)} \leq c_{\text{Thm. 1.5}} (1 \wedge t)^{-\gamma} \tag{1.13}$$

*for all $t \in (0, T]$. In (1.13), $c_{\text{Thm. 1.5}}$ depends on $r$, $w$, $q$, $\gamma$, but not on $t$ and neither on $\mathbf{P}_{x_0}^\infty$.*

*Finally, $\bar{p}$ is unique up to $(\mathrm{d}t \otimes \mathrm{d}\mathbf{P}_{x_0}^\infty)$-a.e. equality and we have that $\bar{p}(t, \cdot) = p(t)$ in $\mathbb{L}^q((\Omega, \mathcal{F}, \mathbf{P}_{x_0}^\infty); \mathsf{H}_r^w(\mathbb{R}^d, \mathrm{d}x))$, for almost every $t \in (0, T]$. In particular, $\bar{p}(t, \omega)$ is the density of $\Lambda_t(\omega)$ $(\mathrm{d}t \otimes \mathrm{d}\mathbf{P}_{x_0}^\infty)$-a.e., and $\bar{p}(t, \cdot)$ satisfies the bound (1.13) for almost every $t \in (0, T]$.*

Sections 2–6 are concerned with proving Theorem 1.5. A discussion of related results, applications and extensions is contained in Section 7.

**Remark 1.6** | The constant $c_{\text{Thm. 1.5}}$ appearing in Theorem 1.5 depends next to $r$, $w$, $q$ and $\gamma$ also on the dimension $d$ and the time horizon $T$, as well as on the primitive quantities from Assumptions 1.1 and 1.2(H), namely the bounds for $\|b\|_\infty$, $\|\bar{\sigma}\|_\infty$, $\|\sigma\|_\infty$, the $\beta$-Hölder-seminorm bounds for $\sup_{t \in [0,T]}[\sigma_t]_\beta$ and $[\bar{\sigma}]_\beta$, and finally the ellipticity bound $\kappa$ from Assumption 1.2(E). For details, we refer to Remark 6.2 after the proof of Theorem 1.5.

## 1.5 Extended probabilistic setting for the proof and strategy

We now introduce two important ingredients for the proof of Theorem 1.5. The first is an extended probabilistic setup which provides the setting for the proof. The second is an approximation scheme which is the backbone of our argument. These ingredients also allow us to look ahead by outlining the main steps in the proof.

**Setting — An extended probabilistic setup** In Section 1.3, we introduced the basic probabilistic setting for Theorem 1.5. To prove the theorem, we extend this by an auxiliary construction, much like a scaffolding put in place to support our work. This scaffolding, although useful, obscures the view on the edifice we subsequently erect under it, and is thus removed in the final step of the proof.

We first introduce the *time simplex*

$$[0, T]^\Delta := \{(\varepsilon, t) \ : \ 0 \leq \varepsilon \leq t \leq T\} \subseteq \mathbb{R}^2$$

and define

$$\bar{\Omega} := \mathbf{C}\big([0, T]^\Delta; \mathcal{S}'(\mathbb{R}^d) \times \mathcal{S}'(\mathbb{R}^d) \times \mathbb{R}^m\big),$$

which is equipped with the topology of convergence uniformly on $[0, T]^\Delta$ for bounded sets in $\mathcal{S}'(\mathbb{R}^d) \times \mathcal{S}'(\mathbb{R}^d) \times \mathbb{R}^m$. From that, we get on $\bar{\Omega}$ the Borel-$\sigma$-algebra which we denote by $\bar{\mathcal{F}}$. On $\bar{\Omega}$, we define the process $(\bar{\Lambda}_t, A_{\varepsilon,t}, \bar{Z}_t)_{(\varepsilon,t) \in [0,T]^\Delta}$ by

$$\big((\varepsilon, t), \bar{\omega}\big) \mapsto (\bar{\Lambda}_t, A_{\varepsilon,t}, Z_t)(\bar{\omega}) = (\bar{\omega}_{0,t}^1, \bar{\omega}_{\varepsilon,t}^2, \bar{\omega}_{0,t}^3) \in \mathcal{S}'(\mathbb{R}^d) \times \mathcal{S}'(\mathbb{R}^d) \times \mathbb{R}^m, \tag{1.14}$$

where $\bar{\omega} = (\bar{\omega}^1, \bar{\omega}^2, \bar{\omega}^2) \in \bar{\Omega}$. In particular, we have

$$A_{\varepsilon,t}(\bar{\omega}) = \bar{\omega}_{\varepsilon,t}^2 \in \mathcal{S}'(\mathbb{R}^d). \tag{1.15}$$

Note that both $\bar{\Lambda}$ and $Z$ are constant in $\varepsilon$. Let us now define

$$E_{\varepsilon,t} := \bar{\Lambda}_t - A_{\varepsilon,t},$$

which gives us the *interpolation decomposition*

$$\bar{\Lambda}_t = A_{\varepsilon,t} + E_{\varepsilon,t} \qquad \text{in } \mathcal{S}'(\mathbb{R}^d) \text{ for } (\varepsilon, t) \in [0, T]^\Delta. \tag{1.16}$$

In this decomposition of $\bar{\Lambda} = (\bar{\Lambda}_t)_{t \in [0,T]}$, we call $A := (A_{\varepsilon,t})_{(\varepsilon,t) \in [0,T]^\Delta}$ the *approximation* and $E := (E_{\varepsilon,t})_{(\varepsilon,t) \in [0,T]^\Delta}$ the *error*, for reasons that will become evident momentarily.

**Step 1 — Approximation scheme**  Recall from Lemma 1.3 the weak solution

$$\boldsymbol{X}^n = \left(\Omega^n_{\text{par}}, \mathcal{G}^n, \mathbb{G}^n, \mathbb{P}^n_{x_0}, (X^n_t)_{t \in [0,T]}, (B^n_t)_{t \in [0,T]}, (Z^n_t)_{t \in [0,T]}\right)$$

of the $n$-particle system $X^n = (X^{1,n}, \ldots, X^{n,n})$ in (1.1). Fix $(\varepsilon, t) \in [0, T]^\Delta$. Then, as a process indexed by $s \in [t - \varepsilon, t]$, we have from (1.1) the dynamics

$$X^{i,n}_s = X^{i,n}_{t-\varepsilon} + \int_{t-\varepsilon}^s b_r(X^{i,n}_r, \mu^n_r)\,\mathrm{d}r + \int_{t-\varepsilon}^s \sigma_r(X^{i,n}_r, \mu^n_r)\,\mathrm{d}B^{i,n}_r + \int_{t-\varepsilon}^s \bar{\sigma}_r(X^{i,n}_r, \mu^n_r)\,\mathrm{d}Z^n_r. \tag{1.17}$$

Inspired by Fournier and Printems [25], we study for $s \in [t - \varepsilon, t]$ the auxiliary process

$$Y^{i,n;\varepsilon,t}_s := X^{i,n}_{t-\varepsilon} + \int_{t-\varepsilon}^s \sigma_r(X^{i,n}_{t-\varepsilon}, \mu^n_{t-\varepsilon})\,\mathrm{d}B^{i,n}_r + \bar{\sigma}_{t-\varepsilon}(X^{i,n}_{t-\varepsilon}, \mu^n_{t-\varepsilon})(Z^n_s - Z^n_{t-\varepsilon}). \tag{1.18}$$

Note that the arguments of $\sigma_r$ and $\bar{\sigma}_{t-\varepsilon}$ are now evaluated at the left endpoint $t-\varepsilon$. From our weak solution $\boldsymbol{X}^n$ and (1.18), this gives for each $n \in \mathbb{N}$ a family of $(\mathbb{R}^d)^n$-valued processes

$$Y^{n;\varepsilon,t} := (Y^{n;t,\varepsilon}_s)_{s \in [t-\varepsilon,t]} := (Y^{1,n;t,\varepsilon}_s, \ldots, Y^{n,n;t,\varepsilon}_s)_{s \in [t-\varepsilon,t]}$$

parametrized by $(\varepsilon, t) \in [0, T]^\Delta$, each defined on $(\Omega^n_{\text{par}}, \mathcal{G}^n, \mathbb{G}^n, \mathbb{P}^n_{x_0})$.

In the seminal contribution [25] and several subsequent works including those of Debussche and Fournier [20], Debussche and Romito [21], Bally [4], Bally and Caramellino [5], Romito [46], the *random variable* $Y^{i,n;\varepsilon,t}_t$ acts as a type of Gaussian approximation for the stochastic integral part in $X^{i,n}_t$, which is used to show that the law of $X^n_t$ possesses a density with respect to Lebesgue measure. In a similar spirit, the $n$ pairs of equations (1.17) and (1.18) suggest an approximation scheme for the *empirical measure* from (1.1),

$$\mu^n_t = \frac{1}{n} \sum_{i=1}^n \delta_{X^{i,n}_t}, \tag{1.19}$$

by the *approximated empirical measure* defined in analogy by

$$\nu^n_{\varepsilon,t} := \frac{1}{n} \sum_{i=1}^n \delta_{Y^{i,n;\varepsilon,t}_t}. \tag{1.20}$$

**Step 2 — The laws of the decomposition**  On the initial setup $(\Omega, \mathcal{F})$ with coordinate process $\Lambda$, we want to study the behavior of cluster points $\mathbf{P}^\infty_{x_0}$ of the sequence $(\mathbf{P}^n_{x_0})_{n \in \mathbb{N}}$. To this end, we *lift* $\mathbf{P}^\infty_{x_0}$ to some $\bar{\mathbf{P}}^\infty_{x_0}$ on the extended space $(\bar{\Omega}, \bar{\mathcal{F}})$ such that $\bar{\mathbf{P}}^\infty_{x_0}$ is a cluster point of the sequence $(\bar{\mathbf{P}}^n_{x_0})_{n \in \mathbb{N}}$ given by

$$\bar{\mathbf{P}}^n_{x_0} := \mathbb{P}^n_{x_0} \circ \left((\mu^n_t, \nu^n_{\varepsilon,t}, Z^n_t)_{(\varepsilon,t) \in [0,T]^\Delta}\right)^{-1}. \tag{1.21}$$

To ensure that (1.21) makes sense and that this sequence has cluster points, we show in Lemma 2.6 below first that each triple $(\mu^n_t, \nu^n_{\varepsilon,t}, Z^n_t)_{(\varepsilon,t) \in [0,T]^\Delta}$ has under $\mathbb{P}^n_{x_0}$ a modification with continuous trajectories in $\mathcal{S}'(\mathbb{R}^d) \times \mathcal{S}'(\mathbb{R}^d) \times \mathbb{R}^m$, and that the now well-defined sequence $(\bar{\mathbf{P}}^n_{x_0})_{n \in \mathbb{N}}$ is tight. By construction, we then also obtain

$$\mathrm{Law}_{\bar{\mathbf{P}}^\infty_{x_0}}(\bar{\Lambda}) = \mathrm{Law}_{\mathbf{P}^\infty_{x_0}}(\Lambda), \tag{1.22}$$

which shows that $\bar{\mathbf{P}}^\infty_{x_0}$ can be viewed as a lifting of $\mathbf{P}^\infty_{x_0}$, and for any $n \in \mathbb{N}$ also

$$\mathrm{Law}_{\bar{\mathbf{P}}^n_{x_0}}(A) = \mathrm{Law}_{\mathbb{P}^n_{x_0}}(\nu^n) \quad \text{and} \quad \mathrm{Law}_{\bar{\mathbf{P}}^n_{x_0}}(E) = \mathrm{Law}_{\mathbb{P}^n_{x_0}}(\mu^n - \nu^n). \tag{1.23}$$

All this is established in (2.25)–(2.37).

In the next two steps, we study in detail the properties of the parts $A$ and $E$ of the decomposition (1.16) of $\bar{\Lambda}$ under the measure $\bar{\mathbf{P}}^\infty_{x_0}$.

**Step 3 — Study of the regular part** As a cluster point, the measure $\bar{\mathbf{P}}^\infty_{x_0}$ is a limiting object, obtained along a subsequence $(\mathbf{P}^{n_k}_{x_0})_{k \in \mathbb{N}}$. To study the distribution of $A$ under $\bar{\mathbf{P}}^\infty_{x_0}$, we use (1.23) and the structure of $\nu^{n_k}$ to provide in Proposition 3.1 a probabilistic representation of $A_{\varepsilon,t}$ under $\bar{\mathbf{P}}^\infty_{x_0}$ as a limiting Gaussian mixture. With this representation, we then proceed to show in Proposition 3.2 that for all $0 < \varepsilon \le t \le T$, the approximation $A_{\varepsilon,t}$ is $\mathbf{P}^\infty_{x_0}$-a.s. in some $\mathsf{H}^s_r(\mathbb{R}^d)$ with $s > 0$, i.e. given by a function, and thus constitutes the regular part of the decomposition (1.16). More precisely, we have for each $\varepsilon \in (0, t]$ that $\mathbf{P}^\infty_{x_0}$-a.s.,

$$\|A_{\varepsilon,t}\|_{\mathsf{H}^s_r(\mathbb{R}^d)} \le a_A \varepsilon^{-p_A(s,r,d)} \quad \text{for a positive function } p_A \text{ and a constant } c_A,$$

with $p_A$ and $c_A$ both nonrandom. Instrumental for this is the uniform ellipticity of $\sigma$ from Assumption 1.2(E). The function $p_A$ is explicit and gives a deterministic quantitative bound on the rate at which $A_{\varepsilon,t}$ becomes singular as $\varepsilon \to 0$.

**Step 4 — Study of the distributional part** We next consider the error $E_{\varepsilon,t}$ under $\bar{\mathbf{P}}^\infty_{x_0}$. For this, we obtain an estimate in some $\mathsf{H}^{-u}_r$-norm, where now $u > 0$ so that $E_{\varepsilon,t}$ is given by a distribution. More precisely, we show in Proposition 4.1 that for all $\xi < \beta$ and *certain but not all* $\varepsilon \in (0, t]$, we have $\mathbf{P}^\infty_{x_0}$-a.s. that

$$\|E_{\varepsilon,t}\|_{\mathsf{H}^{-u}_r(\mathbb{R}^d)} \le C_E \varepsilon^{p_E(\xi)} \quad \text{for a positive function } p_E \text{ and a } random \text{ constant } C_E.$$

This estimate relies on the $\beta$-Hölder-regularity of $\sigma$ and $\bar{\sigma}$ from Assumption 1.2(H) and employs a Kolmogorov-type regularity argument together with properties of the space $\mathsf{H}^{-u}_r(\mathbb{R}^d)$. Here again the nonrandom function $p_E$ is explicit so that the estimate gives us a quantitative bound on the rate at which $E_{\varepsilon,t}$ vanishes as $\varepsilon \to 0$. The constant in this bound is random, but we can show that $C_E \in \mathbb{L}^q(\bar{\mathbf{P}}^\infty_{x_0})$ for all $q \ge 1$, which implies that $C_E < \infty$ almost surely. Finally, $p_E$ is increasing so that $\xi$ can be used to improve the rate at which $E_{\varepsilon,t}$ vanishes as $\varepsilon \to 0$, as long as $\xi < \beta$.

Having obtained in Step 1 the decomposition (1.16) of $\bar{\Lambda}$ into a sum of $A$ and $E$, and in Steps 3 and 4 bounds on the parts $A$ and $E$, we next present an analytic tool that allows us to combine our estimates to obtain an improved bound for their sum $\bar{\Lambda}$.

**Step 5 — An abstract interpolation result** Our starting point is a distribution $\lambda \in \mathcal{S}'$ and a family of linear decompositions $\lambda = a_\varepsilon + e_\varepsilon$ indexed by appropriate $\varepsilon > 0$ into a 'regular part' $a_\varepsilon$ and a 'distributional part' $e_\varepsilon$. For the respective parts, we assume that we know the rates at which $\|a_\varepsilon\|_{\mathsf{H}^s_r(\mathbb{R}^d)} \to \infty$ and $\|e_\varepsilon\|_{\mathsf{H}^{-u}_r(\mathbb{R}^d)} \to 0$ as $\varepsilon \to 0$. For $a_\varepsilon$, the rate depends on the primitive parameters $r$, $s$, $d$, while for $e_\varepsilon$, it may depend on an additional parameter $\xi$.

In this setting, interpolation theory provides an efficient tool to trade off the rate at which the regular part becomes singular against the rate at which the distributional part vanishes, allowing us to obtain for $w \in (-u, s)$ a bound for the $\mathsf{H}^w_r(\mathbb{R}^d)$-norm of $\lambda$, i.e., to show that $\lambda$ lies in the *intermediate space* $\mathsf{H}^w_r(\mathbb{R}^d)$; in symbols,

$$\mathcal{S}(\mathbb{R}^d) \hookrightarrow \mathsf{H}^s_r(\mathbb{R}^d) \hookrightarrow \mathsf{H}^w_r(\mathbb{R}^d) \hookrightarrow \mathsf{H}^{-u}_r(\mathbb{R}^d) \hookrightarrow \mathcal{S}'(\mathbb{R}^d).$$

Moreover, for certain ranges of the parameters $r$, $s$, $u$, $\xi$, we can show that $w > 0$, i.e. that $\lambda$ is a function, not merely a distribution.

The interpolation result from Step 5 is the tool we need to complete the proof of Theorem 1.5. It allows us to combine the decomposition and the estimates from Steps 2–4 in a way that respects the restrictions we face: The estimates from Steps 3 and 4 may fail on certain nullsets that depend on $(\varepsilon, t) \in [0, T]^\Delta$, and the estimates in Step 4 are random and available only for certain pairs $(\varepsilon, t)$.

**Step 6 — Proof of Theorem 1.5** For a set $\bar{\Omega}\backslash N$ of full $\bar{\mathbf{P}}^\infty_{x_0}$-measure, we apply the interpolation result obtained in Step 5 to the distribution $\Lambda_t(\omega) \in \mathcal{S}'$, pointwise for each $\omega \in \bar{\Omega}\backslash N$. The first ingredient for the interpolation is the linear decomposition (1.16) of $\bar{\Lambda}_t(\omega)$ from Step 2. The second ingredient is provided by Steps 3 and 4 from which $\|A_{\varepsilon,t}\|_{\mathsf{H}^s_r(\mathbb{R}^d)} \to \infty$ and $\|E_{\varepsilon,t}\|_{\mathsf{H}^{-u}_r(\mathbb{R}^d)} \to 0$ at rates which are independent of $\omega$, but with constants that may depend on $\omega$.

The $\beta$-Hölder-regularity of $\sigma$ and $\bar{\sigma}$ from Assumption 1.2(H) creates a gap between these rates; see (5.24) below. This allows us to make appropriate choices of $\xi$, which controls the rate at which $E_{\varepsilon,t}$ vanishes, and of $r$ and $s$, which determine the rate at which $A_{\varepsilon,t}$ explodes, to obtain a *random bound* on $\|\bar{\Lambda}_t\|_{\mathsf{H}^w_r(\mathbb{R}^d)}$ for all intermediate values $0 < w < w_0(d, \xi, s, r)$. Moreover, the estimates from Steps 3 and 4 and the interpolation result from Step 5 allow us to obtain in addition a moment estimate for the random norm $\|\bar{\Lambda}_t\|_{\mathsf{H}^w_r(\mathbb{R}^d)}$, namely that $\|\bar{\Lambda}_t\|_{\mathsf{H}^w_r(\mathbb{R}^d)} \in \mathbb{L}^q(\bar{\mathbf{P}}^\infty_{x_0})$ for all $q \geq 1$.

All this is obtained on the extended probabilistic setup, i.e. for $\bar{\Lambda}$ on $(\bar{\Omega}, \bar{\mathcal{F}}, \bar{\mathbf{P}}^\infty_{x_0})$. To complete the proof, we remove this temporary extension by projecting the estimate onto our initial setup $(\Omega, \mathcal{F}, \mathbf{P}^\infty_{x_0})$. In this way, we get a norm bound for $\Lambda$, and this ultimately leads to the estimate (1.13) from Theorem 1.5.

## 2 THE LAWS OF THE DECOMPOSITION

In Section 1.5, as Step 1 of the proof of Theorem 1.5, we already constructed in (1.20) for each $n \in \mathbb{N}$ and $(\varepsilon, t) \in [0, T]^\Delta$ the approximation $\nu^n_{\varepsilon,t}$ of the empirical measure $\mu^n_t$ from (1.19). In this section, we carry out Step 2. This needs an estimate for the regularity of

$$(\varepsilon, t) \mapsto \nu^n_{\varepsilon,t} = \frac{1}{n} \sum_{i=1}^n \delta_{Y^{i,n;\varepsilon,t}_t}$$

as an $\mathcal{S}'$-valued flow of the two parameters $(\varepsilon, t) \in [0, T]^\Delta$. With such an estimate, we then show that the sequence $(\mu^n, \nu^n, Z^n)$ is tight, and this allows us to control the joint law of the interpolation decomposition from (1.16). The main building-block for all this is Proposition 2.2 below.

## 2.1 A key regularity estimate

The following simple *bounded coefficients estimate* is used below.

**Lemma 2.1** | *Suppose the $\mathbb{R}^k$-valued process $S = (S_t)_{t \in [0,T]}$ has the form*

$$S_t = S_0 + \int_0^t a_r \, \mathrm{d}r + \int_0^t m_r \, \mathrm{d}W_r \qquad \text{with } S_0 \in \mathbb{R}^k,$$

*where $W = (W_t)_{t \in [0,T]}$ is a Brownian motion and $a = (a_t)_{t \in [0,T]}$, $m = (m_t)_{t \in [0,T]}$ are predictable processes which are bounded uniformly in $t$ and $\omega$. For any $0 \leq t \leq t' \leq T$ and any $q \in [1, \infty)$, we then have the estimate*

$$\Big\| \sup_{t \leq s \leq t'} |S_s - S_t| \Big\|_{\mathbb{L}^q} \leq \|a\|_\infty |t' - t| + c\|m\|_\infty |t' - t|^{\frac{1}{2}}$$

*for a constant $c = c(q) < \infty$ which depends only on $q$.*

**Proof** We use the triangle inequality. For $A := \int a_r \, \mathrm{d}r$, we clearly have

$$\sup_{t \leq s \leq t'} |A_s| \leq \int_t^{t'} |a_r| \, \mathrm{d}r \leq \|a\|_\infty |t' - t| \qquad \text{a.s.}$$

For $M := \int m_r \, \mathrm{d}W_r$ we use the Burkholder–Davis–Gundy inequality to obtain

$$\Big\| \sup_{t \leq s \leq t'} |M_s| \Big\|_{\mathbb{L}^q} \leq c(q) \Big\| \Big( \int_t^{t'} |m_r|^2 \, \mathrm{d}r \Big)^{\frac{1}{2}} \Big\|_{\mathbb{L}^q} \leq c(q) \|m\|_\infty |t' - t|^{1/2}.$$

$\square$

We also rely on the following key regularity estimate. Recall that $\beta$ is the Hölder-exponent of $\sigma, \bar{\sigma}$ from Assumption 1.2(H).

**Proposition 2.2** | *Under Assumptions 1.1 and 1.2 and for all $q \geq 1$, there exists a constant $c = c(q, b, \sigma, \bar{\sigma}, T, \beta) < \infty$ such that*

$$\|Y_{t'}^{i,n;\varepsilon',t'} - Y_t^{i,n;\varepsilon,t}\|_{\mathbb{L}^q(\mathbb{P}_{x_0}^n)} \leq c|(t - \varepsilon, t) - (t' - \varepsilon', t')|^{\beta/2} \qquad (2.1)$$

*for all $(\varepsilon, t), (\varepsilon', t') \in [0, T]^\Delta$, $i \in [n]$ and $n \in \mathbb{N}$.*

**Proof** Without loss of generality, assume $t - \varepsilon \leq t' - \varepsilon'$. For each $n \in \mathbb{N}$ and $i \in [n]$, the dynamics of $(Y_s^{i,n;\varepsilon,t})_{s \in [t-\varepsilon,t],(\varepsilon,t) \in [0,T]^\Delta}$ in (1.18) give

$$Y_{t'}^{i,n;\varepsilon',t'} - Y_t^{i,n;\varepsilon,t} = F_0^{i,n} + F_1^{i,n} - F_2^{i,n} \qquad (2.2)$$

# EMERGENCE OF REGULARITY

with

$$F_0^{i,n} := X_{t'-\varepsilon'}^{i,n} - X_{t-\varepsilon}^{i,n}, \tag{2.3}$$

$$F_1^{i,n} := \int_{t'-\varepsilon'}^{t'} \sigma_r(X_{t'-\varepsilon'}^{i,n}, \mu_{t'-\varepsilon'}^n) \, dB_r^{i,n} + \int_{t'-\varepsilon'}^{t'} \bar{\sigma}_{t'-\varepsilon'}(X_{t'-\varepsilon'}^{i,n}, \mu_{t'-\varepsilon'}^n) \, dZ_r^n, \tag{2.4}$$

$$F_2^{i,n} := \int_{t-\varepsilon}^{t} \sigma_r(X_{t-\varepsilon}^{i,n}, \mu_{t-\varepsilon}^n) \, dB_r^{i,n} + \int_{t-\varepsilon}^{t} \bar{\sigma}_{t-\varepsilon}(X_{t-\varepsilon}^{i,n}, \mu_{t-\varepsilon}^n) \, dZ_r^n. \tag{2.5}$$

To establish (2.1), we derive bounds for $F_0^{i,n}$ and $F_1^{i,n} - F_2^{i,n}$ in $\mathbb{L}^q(\mathbb{P}_{x_0}^n)$ which are uniform in $i \in [n]$ and $n \in \mathbb{N}$.

**Step 1** Starting from (2.3) and using the dynamics of $X^n$ in (1.17), we get with Assumptions 1.1 and 1.2 from the bounded coefficients estimate in Lemma 2.1 that

$$\begin{aligned}
\|F_0^{i,n}\|_{\mathbb{L}^q(\mathbb{P}_{x_0}^n)} &\leq \|b\|_\infty |(t'-\varepsilon') - (t-\varepsilon)| + c(q)(\|\sigma\|_\infty + \|\bar{\sigma}\|_\infty)|(t'-\varepsilon') - (t-\varepsilon)|^{\frac{1}{2}} \\
&\leq \left(\|b\|_\infty T^{\frac{1}{2}} + c(q)(\|\sigma\|_\infty + \|\bar{\sigma}\|_\infty)\right)|(t'-\varepsilon') - (t-\varepsilon)|^{\frac{1}{2}} \\
&\leq T^{\frac{1-\beta}{2}}|(t'-\varepsilon') - (t-\varepsilon)|^{\frac{\beta}{2}}\left(\|b\|_\infty T^{\frac{1}{2}} + c(q)(\|\sigma\|_\infty + \|\bar{\sigma}\|_\infty)\right).
\end{aligned} \tag{2.6}$$

**Step 2** Clearly, we always have $t - \varepsilon \leq t$ and $t' - \varepsilon' \leq t'$. Recall also our assumption that $t - \varepsilon \leq t' - \varepsilon'$, so that $t - \varepsilon \leq t' - \varepsilon' \leq t'$. To control the term $F_1^{i,n} - F_2^{i,n}$, we distinguish three cases, depending on the position of $t$ relative to $t - \varepsilon \leq t' - \varepsilon' \leq t'$.

**Step 2.1** The case $t - \varepsilon \leq t \leq t' - \varepsilon' \leq t'$ where $t$ is to the left of $t' - \varepsilon'$ is clearest. The boundedness of $\sigma, \bar{\sigma}$ from Assumption 1.2 lets us use the bounded coefficients estimate in Lemma 2.1 to find with the assumed ordering $t \leq t' - \varepsilon' \leq t'$ that

$$\begin{aligned}
\|F_1^{i,n}\|_{\mathbb{L}^q(\mathbb{P}_{x_0}^n)} &\leq c(q)(\|\sigma\|_\infty + \|\bar{\sigma}\|_\infty)|t' - (t'-\varepsilon')|^{\frac{1}{2}} \\
&\leq c(q)(\|\sigma\|_\infty + \|\bar{\sigma}\|_\infty)|t' - t|^{\frac{1}{2}} \\
&\leq T^{\frac{1-\beta}{2}}|t' - t|^{\frac{\beta}{2}} c(q)(\|\sigma\|_\infty + \|\bar{\sigma}\|_\infty) \\
&\leq c(T, q, \beta, \sigma, \bar{\sigma})|t' - t|^{\frac{\beta}{2}}.
\end{aligned} \tag{2.7}$$

An analogous argument using the ordering $t - \varepsilon \leq t \leq t' - \varepsilon'$ yields

$$\begin{aligned}
\|F_2^{i,n}\|_{\mathbb{L}^q(\mathbb{P}_{x_0}^n)} &\leq c(q)(\|\sigma\|_\infty + \|\bar{\sigma}\|_\infty)|t - (t-\varepsilon)|^{\frac{1}{2}} \\
&\leq c(T, q, \beta, \sigma, \bar{\sigma})|(t'-\varepsilon') - (t-\varepsilon)|^{\frac{\beta}{2}}.
\end{aligned} \tag{2.8}$$

The decomposition (2.2) with the bounds (2.6)–(2.8) now gives

$$\|Y_{t'}^{i,n;\varepsilon',t'} - Y_t^{i,n;\varepsilon,t}\|_{\mathbb{L}^q(\mathbb{P}_{x_0}^n)} \leq c(T, q, \beta, b, \sigma, \bar{\sigma}) \max\{|t' - t|^{\frac{\beta}{2}}, |(t'-\varepsilon') - (t-\varepsilon)|^{\frac{\beta}{2}}\}. \tag{2.9}$$

Note that the bound in (2.9) comes from the lengths of the stochastic integrals, namely $(t'-\varepsilon') - (t-\varepsilon)$ for $F_0^{i,n}$ and $t' - (t'-\varepsilon') \leq t' - t$ in the present case for $F_1^{i,n}$ and similarly $t - (t'-\varepsilon') \leq (t'-\varepsilon') - (t-\varepsilon)$ for $F_2^{i,n}$.

**Step 2.2** The case $t - \varepsilon \leq t' - \varepsilon' \leq t \leq t'$ where $t$ lies between $t' - \varepsilon'$ and $t'$ requires more care because the intervals $[t-\varepsilon, t]$ and $[t'-\varepsilon', t']$ overlap. In a first step, we split $[t-\varepsilon, t]$ into $[t-\varepsilon, t'-\varepsilon']$ and $[t'-\varepsilon', t]$, and similarly $[t'-\varepsilon', t']$ into $[t'-\varepsilon', t]$ and $[t, t']$. Splitting the integrals in $F_1^{i,n}$ and $F_2^{i,n}$ from (2.4) and (2.5) accordingly gives

$$F_1^{i,n} - F_2^{i,n} = (G_1^{i,n} + G_2^{i,n}) - (G_3^{i,n} + G_4^{i,n}) \tag{2.10}$$

with

$$G_1^{i,n} := \int_{t'-\varepsilon'}^{t} \sigma_r(X_{t'-\varepsilon'}^{i,n}, \mu_{t'-\varepsilon'}^n) \, dB_r^{i,n} + \int_{t'-\varepsilon'}^{t} \bar{\sigma}_{t'-\varepsilon'}(X_{t'-\varepsilon'}^{i,n}, \mu_{t'-\varepsilon'}^n) \, dZ_r^n,$$

$$G_2^{i,n} := \int_{t}^{t'} \cdots dB_r^{i,n} + \int_{t}^{t'} \cdots dZ_r^n,$$

$$G_3^{i,n} := \int_{t-\varepsilon}^{t'-\varepsilon'} \sigma_r(X_{t-\varepsilon}^{i,n}, \mu_{t-\varepsilon}^n) \, dB_r^{i,n} + \int_{t-\varepsilon}^{t'-\varepsilon'} \bar{\sigma}_{t-\varepsilon}(X_{t-\varepsilon}^{i,n}, \mu_{t-\varepsilon}^n) \, dZ_r^n,$$

$$G_4^{i,n} := \int_{t'-\varepsilon'}^{t} \cdots dB_r^{i,n} + \int_{t'-\varepsilon'}^{t} \cdots dZ_r^n.$$

For the non-overlapping terms $G_2^{i,n}$ and $G_3^{i,n}$, we proceed as in (2.7) and (2.8) to find

$$\|G_2^{i,n}\|_{\mathbb{L}^q(\mathbb{P}_{x_0}^n)} \le c(T, q, \beta, \sigma, \bar{\sigma})|t' - t|^{\frac{\beta}{2}}, \tag{2.11}$$

$$\|G_3^{i,n}\|_{\mathbb{L}^q(\mathbb{P}_{x_0}^n)} \le c(T, q, \beta, \sigma, \bar{\sigma})|(t' - \varepsilon') - (t - \varepsilon)|^{\frac{\beta}{2}}. \tag{2.12}$$

For the overlapping terms $G_1^{i,n}$ and $G_4^{i,n}$ on $[t'-\varepsilon', t]$, the triangle inequality gives

$$\|G_1^{i,n} - G_4^{i,n}\|_{\mathbb{L}^q(\mathbb{P}_{x_0}^n)} \le \|G_5^{i,n}\|_{\mathbb{L}^q(\mathbb{P}_{x_0}^n)} + \|G_6^{i,n}\|_{\mathbb{L}^q(\mathbb{P}_{x_0}^n)} \tag{2.13}$$

with individual terms given by

$$G_5^{i,n} := \int_{t'-\varepsilon'}^{t} \Big( \sigma_r(X_{t'-\varepsilon'}^{i,n}, \mu_{t'-\varepsilon'}^n) - \sigma_r(X_{t-\varepsilon}^{i,n}, \mu_{t-\varepsilon}^n) \Big) \, dB_r^{i,n}, \tag{2.14}$$

$$G_6^{i,n} := \int_{t'-\varepsilon'}^{t} \Big( \bar{\sigma}_{t'-\varepsilon'}(X_{t'-\varepsilon'}^{i,n}, \mu_{t'-\varepsilon'}^n) - \bar{\sigma}_{t-\varepsilon}(X_{t-\varepsilon}^{i,n}, \mu_{t-\varepsilon}^n) \Big) \, dZ_r^n. \tag{2.15}$$

For $G_5^{i,n}$, we note that since $(x, \mu) \mapsto \sigma_r(x, \mu)$ is by Assumption 1.2(H) $\beta$-Hölder-continuous, uniformly for $r \in [0, T]$, we get the estimate

$$|\sigma_r(X_{t'-\varepsilon'}^{i,n}, \mu_{t'-\varepsilon'}^n) - \sigma_r(X_{t-\varepsilon}^{i,n}, \mu_{t-\varepsilon}^n)| \le \sup_{r \in [0,T]} [\sigma_r]_\beta \Big( |X_{t'-\varepsilon'}^{i,n} - X_{t-\varepsilon}^{i,n}|^\beta + d_{\mathbb{W}_1}(\mu_{t'-\varepsilon'}^n - \mu_{t-\varepsilon}^n)^\beta \Big),$$

where $[\sigma_r]_\beta$ denotes the usual $\beta$-Hölder-seminorm. Note that this bound is independent of $r$. Hence, using the Burkholder–Davis–Gundy inequality shows that

$$\|G_5^{i,n}\|_{\mathbb{L}^q(\mathbb{P}_{x_0}^n)} \le c(q) \, \mathbf{E}^{\mathbb{P}_{x_0}^n} \left[ \left( \int_{t'-\varepsilon'}^{t} |\sigma_r(X_{t'-\varepsilon'}^{i,n}, \mu_{t'-\varepsilon'}^n) - \sigma_r(X_{t-\varepsilon}^{i,n}, \mu_{t-\varepsilon}^n)|^2 \, dr \right)^{\frac{q}{2}} \right]^{\frac{1}{q}}$$

$$\le c(q) \sup_{r \in [0,T]} [\sigma_r]_\beta |t - (t' - \varepsilon')|^{\frac{1}{2}} \Big( \big\| |X_{t'-\varepsilon'}^{i,n} - X_{t-\varepsilon}^{i,n}|^\beta \big\|_{\mathbb{L}^q(\mathbb{P}_{x_0}^n)}$$

$$+ \| d_{\mathbb{W}_1}(\mu_{t'-\varepsilon'}^n, \mu_{t-\varepsilon}^n)^\beta \|_{\mathbb{L}^q(\mathbb{P}_{x_0}^n)} \Big).$$

For the second term, we proceed similarly. Since $(r, x, \mu) \mapsto \bar{\sigma}_r(x, \mu)$ is by Assumption 1.2(H) $\beta$-Hölder-continuous, we get

$$|\bar{\sigma}_{t'-\varepsilon'}(X^{i,n}_{t'-\varepsilon'}, \mu^n_{t'-\varepsilon'}) - \bar{\sigma}_{t-\varepsilon}(X^{i,n}_{t-\varepsilon}, \mu^n_{t-\varepsilon})| \leq [\bar{\sigma}]_\beta \Big(|(t'-\varepsilon') - (t-\varepsilon)|^\beta + |X^{i,n}_{t'-\varepsilon'} - X^{i,n}_{t-\varepsilon}|^\beta$$
$$+ \mathsf{d}_{\mathbb{W}_1}(\mu^n_{t'-\varepsilon'} - \mu^n_{t-\varepsilon})^\beta \Big).$$

This estimate again is independent of $r$; so the Burkholder–Davis–Gundy inequality gives

$$\|G^{i,n}_6\|_{\mathbb{L}^q(\mathbb{P}^n_{x_0})} \leq c(q) \, \mathbf{E}^{\mathbb{P}^n_{x_0}} \left[ \left( \int_{t'-\varepsilon'}^t |\bar{\sigma}_{t'-\varepsilon'}(X^{i,n}_{t'-\varepsilon'}, \mu^n_{t'-\varepsilon'}) - \bar{\sigma}_{t-\varepsilon}(X^{i,n}_{t-\varepsilon}, \mu^n_{t-\varepsilon})|^2 \, dr \right)^{\frac{q}{2}} \right]^{\frac{1}{q}}$$
$$\leq c(q)[\bar{\sigma}]_\beta |t - (t'-\varepsilon')|^{\frac{1}{2}} \Big( |(t'-\varepsilon') - (t-\varepsilon)|^\beta + \left\| |X^{i,n}_{t'-\varepsilon'} - X^{i,n}_{t-\varepsilon}|^\beta \right\|_{\mathbb{L}^q(\mathbb{P}^n_{x_0})}$$
$$+ \| \mathsf{d}_{\mathbb{W}_1}(\mu^n_{t'-\varepsilon'}, \mu^n_{t-\varepsilon})^\beta \|_{\mathbb{L}^q(\mathbb{P}^n_{x_0})} \Big).$$

Combining these estimates via (2.13) and then using that for $0 < \beta < 1$, the concavity of $x \mapsto x^\beta$ implies that $\|U^\beta\|_{\mathbb{L}^q(\mathbb{P}^n_{x_0})} \leq \|U\|^\beta_{\mathbb{L}^q(\mathbb{P}^n_{x_0})}$ now shows

$$\|G^{i,n}_1 - G^{i,n}_4\|_{\mathbb{L}^q(\mathbb{P}^n_{x_0})} \leq c(q) \Big( \sup_{r \in [0,T]} [\sigma_r]_\beta + [\bar{\sigma}]_\beta \Big) |t - (t'-\varepsilon')|^{\frac{1}{2}}$$
$$\times \Big( |(t'-\varepsilon') - (t-\varepsilon)|^\beta + \left\| |X^{i,n}_{t'-\varepsilon'} - X^{i,n}_{t-\varepsilon}|^\beta \right\|_{\mathbb{L}^q(\mathbb{P}^n_{x_0})}$$
$$+ \| \mathsf{d}_{\mathbb{W}_1}(\mu^n_{t'-\varepsilon'}, \mu^n_{t-\varepsilon})^\beta \|_{\mathbb{L}^q(\mathbb{P}^n_{x_0})} \Big)$$
$$\leq c(q) \Big( \sup_{r \in [0,T]} [\sigma_r]_\beta + [\bar{\sigma}]_\beta \Big) |t - (t'-\varepsilon')|^{\frac{1}{2}}$$
$$\times \Big( |(t'-\varepsilon') - (t-\varepsilon)|^\beta + \|X^{i,n}_{t'-\varepsilon'} - X^{i,n}_{t-\varepsilon}\|^\beta_{\mathbb{L}^q(\mathbb{P}^n_{x_0})}$$
$$+ \| \mathsf{d}_{\mathbb{W}_1}(\mu^n_{t'-\varepsilon'}, \mu^n_{t-\varepsilon})\|^\beta_{\mathbb{L}^q(\mathbb{P}^n_{x_0})} \Big). \quad (2.16)$$

The second term in the brackets on the right-hand side, we recognize as $F^{i,n}_0$ from (2.3) and hence can use (2.6) to bound it by $c(q, \sigma, \bar{\sigma}, b, T, \beta)|(t'-\varepsilon') - (t-\varepsilon)|^{\beta/2}$. For the first term in the brackets, we use $|(t'-\varepsilon') - (t-\varepsilon)|^\beta \leq T^{\beta/2}|(t'-\varepsilon') - (t-\varepsilon)|^{\beta/2}$ which has the same form. Finally, for the third term in the brackets, we use the definition (A.11) of the Kantorovich–Rubinstein distance and the definition (1.19) of $\mu^n$ as the empirical measure to get

$$\mathsf{d}_{\mathbb{W}_1}(\mu^n_{t'-\varepsilon'}, \mu^n_{t-\varepsilon}) = \sup_{\|f\|_{\mathsf{Lip}} \leq 1} \left| \int_{\mathbb{R}^d} f(x) \, (\mu^n_{t'-\varepsilon'} - \mu^n_{t-\varepsilon})(dx) \right|$$
$$= \sup_{\|f\|_{\mathsf{Lip}} \leq 1} \left| \frac{1}{n} \sum_{j=1}^n \Big( f(X^{j,n}_{t'-\varepsilon'}) - f(X^{j,n}_{t-\varepsilon}) \Big) \right|$$
$$\leq \frac{1}{n} \sum_{j=1}^n |X^{j,n}_{t'-\varepsilon'} - X^{j,n}_{t-\varepsilon}|. \quad (2.17)$$

Using the triangle inequality and the fact that the $X^{j,n}$, $j = 1, \ldots, n$, are identically distributed because they are exchangeable by Lemma 1.3 therefore yields

$$\| \mathsf{d}_{\mathbb{W}_1}(\mu^n_{t'-\varepsilon'}, \mu^n_{t-\varepsilon}) \|_{\mathbb{L}^q(\mathbb{P}^n_{x_0})} \leq \frac{1}{n} \sum_{j=1}^n \|X^{j,n}_{t'-\varepsilon'} - X^{j,n}_{t-\varepsilon}\|_{\mathbb{L}^q(\mathbb{P}^n_{x_0})} = \|X^{i,n}_{t'-\varepsilon'} - X^{i,n}_{t-\varepsilon}\|_{\mathbb{L}^q(\mathbb{P}^n_{x_0})}.$$

In consequence, the third term in the brackets on the right-hand side of (2.16) can be bounded exactly like the first term, and so we end up with

$$\|G_1^{i,n} - G_4^{i,n}\|_{\mathbb{L}^q(\mathbb{P}_{x_0}^n)} \leq c(q,\sigma,\bar{\sigma},b,T,\beta)|t-(t'-\varepsilon')|^{\frac{1}{2}}|(t'-\varepsilon')-(t-\varepsilon)|^{\frac{\beta}{2}} \qquad (2.18)$$

$$\leq c(q,\sigma,\bar{\sigma},b,T,\beta)|(t'-\varepsilon')-(t-\varepsilon)|^{\frac{\beta}{2}}. \qquad (2.19)$$

The decomposition of $F_1^{i,n} - F_2^{i,n}$ in (2.10) and the bounds (2.11) and (2.12) thus imply

$$\|F_1^{i,n} - F_2^{i,n}\|_{\mathbb{L}^q(\mathbb{P}_{x_0}^n)} \leq c(q,\sigma,\bar{\sigma},b,T,\beta)\max\{|t'-t|^{\frac{\beta}{2}}, |(t'-\varepsilon')-(t-\varepsilon)|^{\frac{\beta}{2}}\}.$$

To conclude, we combine this with the bound for $F_0^{i,n}$ in (2.6) and use the decomposition of $Y_{t'}^{i,n;\varepsilon',t'} - Y_t^{i,n;\varepsilon,t}$ in (2.2) to find

$$\|Y_{t'}^{i,n;\varepsilon',t'} - Y_t^{i,n;\varepsilon,t}\|_{\mathbb{L}^q(\mathbb{P}_{x_0}^n)} \leq c(q,\sigma,\bar{\sigma},b,T,\beta)\max\{|t'-t|^{\frac{\beta}{2}}, |(t'-\varepsilon')-(t-\varepsilon)|^{\frac{\beta}{2}}\}, \qquad (2.20)$$

which has the same form as (2.9).

**Step 2.3** The case $t - \varepsilon \leq t' - \varepsilon' \leq t' \leq t$ where $t$ lies to the right of $t'$ can be treated in entire analogy to the preceding case. We split the interval $[t-\varepsilon, t]$ into $[t-\varepsilon, t'-\varepsilon']$, $[t'-\varepsilon', t']$, and $[t', t]$. Splitting the integrals in $F_2^{i,n}$ from (2.5) accordingly gives

$$F_1^{i,n} - F_2^{i,n} = F_1^{i,n} - (H_1^{i,n} + H_2^{i,n} + H_3^{i,n}) \qquad (2.21)$$

with $F_1^{i,n}$ as in (2.4) and

$$H_1^{i,n} := \int_{t-\varepsilon}^{t'-\varepsilon'} \sigma_r(X_{t-\varepsilon}^{i,n}, \mu_{t-\varepsilon}^n)\,dB_r^{i,n} + \int_{t-\varepsilon}^{t'-\varepsilon'} \bar{\sigma}_{t-\varepsilon}(X_{t-\varepsilon}^{i,n}, \mu_{t-\varepsilon}^n)\,dZ_r^n,$$

$$H_2^{i,n} := \int_{t'-\varepsilon'}^{t'} \cdots dB_r^{i,n} + \int_{t'-\varepsilon'}^{t'} \cdots dZ_r^n,$$

$$H_3^{i,n} := \int_{t'}^{t} \cdots dB_r^{i,n} + \int_{t'}^{t} \cdots dZ_r^n.$$

The non-overlapping terms are $H_1^{i,n} = G_3^{i,n}$ and $H_3^{i,n}$, which is almost $G_2^{i,n}$ except that $t$ and $t'$ are interchanged. So (2.12) and (2.11) can again be used, and it remains to look at the difference of $F_1^{i,n}$ and $H_2^{i,n}$ which are both integrals over $[t'-\varepsilon', t']$. But this is almost the same situation as with $G_1^{i,n}$ and $G_4^{i,n}$, except that $t$ in the integrals is replaced by $t'$. Following the same steps as for (2.13) up to (2.19) therefore yields

$$\|F_1^{i,n} - H_2^{i,n}\|_{\mathbb{L}^q(\mathbb{P}_{x_0}^n)} \leq c(q,\sigma,\bar{\sigma},b,T,\beta)|(t'-\varepsilon')-(t-\varepsilon)|^{\frac{\beta}{2}}, \qquad (2.22)$$

and combining this with the estimates (2.12) and (2.11) used for $H_1^{i,n}$ and $H_3^{i,n}$ leads via (2.21) and the decomposition (2.2) for $Y_{t'}^{i,n;\varepsilon',t'} - Y_t^{i,n;\varepsilon,t}$ to

$$\|Y_{t'}^{i,n;\varepsilon',t'} - Y_t^{i,n;\varepsilon,t}\|_{\mathbb{L}^q(\mathbb{P}_{x_0}^n)} \leq c(q,\sigma,\bar{\sigma},b,T,\beta)\max\{|t'-t|^{\frac{\beta}{2}}, |(t'-\varepsilon')-(t-\varepsilon)|^{\frac{\beta}{2}}\}, \qquad (2.23)$$

which again has the same form as (2.9).

**Step 3** To conclude the proof, we observe that in all cases, the resulting estimates (2.9), (2.20), (2.23) have the same form. Since all norms on $\mathbb{R}^2$ are equivalent, we can slightly change the constant to transform the max-norm into the Euclidean norm so that

$$\|Y_{t'}^{i,n;\varepsilon',t'} - Y_t^{i,n;\varepsilon,t}\|_{\mathbb{L}^q(\mathbb{P}_{x_0}^n)} \leq c(q,\sigma,\bar{\sigma},b,T,\beta,d)|(t-\varepsilon,t) - (t'-\varepsilon',t')|^{\frac{\beta}{2}}.$$

This gives (2.1) and thus establishes the result. □

**Remark 2.3** | Note that in the case $\varepsilon = 0$ and $t = t'$, we get from (2.18) the more potent estimate $\|G_1^{i,n} - G_4^{i,n}\|_{\mathbb{L}^q(\mathbb{P}_{x_0}^n)} \leq c(q,\sigma,\bar{\sigma},b,T,\beta)|\varepsilon'|^{(1+\beta)/2}$. This is used in Lemma 4.2 below.

## 2.2 Tightness

We next recall a version of the Sobolev embedding theorem for Bessel potential spaces. For this, let $0 < \rho < 1$ and define for $f : \mathbb{R}^d \to \mathbb{R}$ the seminorm

$$[f]_\rho := \sup_{x \neq y} \frac{|f(x) - f(y)|}{|x-y|^\rho},$$

where $\mathbb{R}^d$ is equipped with the usual Euclidean norm $|\cdot|$. This gives the $\rho$-Hölder space

$$\mathsf{H\ddot{o}l}_\rho(\mathbb{R}^d) := \{f : \mathbb{R}^d \to \mathbb{R} \ : \ \|f\|_{\mathsf{H\ddot{o}l}_\rho} < \infty\}$$

with its usual norm $\|\cdot\|_{\mathsf{H\ddot{o}l}_\rho} := \|\cdot\|_\infty + [\cdot]_\rho$.

**Lemma 2.4** | *Let $r \in (1,\infty)$ with conjugate Hölder exponent $r'$, i.e., $1/r + 1/r' = 1$, and set $u := 1 + d/r'$. Then the embedding $\mathsf{H}_{r'}^u(\mathbb{R}^d) \hookrightarrow \mathsf{H\ddot{o}l}_\rho(\mathbb{R}^d)$ is continuous for all $0 < \rho < 1$, i.e., there exists a constant $c = c(\rho,r,u) < \infty$ with $\|f\|_{\mathsf{H\ddot{o}l}_\rho} \leq c\|f\|_{\mathsf{H}_{r'}^u}$ for all $f \in \mathsf{H}_{r'}^u(\mathbb{R}^d)$.*

**Proof** Details are found in Appendix B. □

**Lemma 2.5** | *Let $s \in \mathbb{R}$ and $r \in (1,\infty)$ with conjugate Hölder exponent $r'$. Then $\mathsf{H}_r^s(\mathbb{R}^d)$ is a Banach space. Moreover, its topological dual $(\mathsf{H}_r^s(\mathbb{R}^d))'$ is isomorphic to $\mathsf{H}_{r'}^{-s}(\mathbb{R}^d)$, i.e., there exists $c = c(r,s) \in (0,\infty)$ with $c^{-1}\|f\|_{(\mathsf{H}_r^s)'} \leq \|f\|_{\mathsf{H}_{r'}^{-s}} \leq c\|f\|_{(\mathsf{H}_r^s)'}$ for all $f \in \mathsf{H}_{r'}^{-s}(\mathbb{R}^d)$.*

**Proof** See Bergh and Löfström [6, Cor. 6.2.8]. □

From the regularity estimate in Proposition 2.2, we obtain a *joint* tightness property.

**Lemma 2.6** | *Let Assumptions 1.1 and 1.2 be in force. For each $n \in \mathbb{N}$, the family $(\nu_{\varepsilon,t}^n)_{(\varepsilon,t) \in [0,T]^\Delta}$ is a set of $\mathcal{S}'$-valued random variables. Moreover,*

$$[0,T]^\Delta \ni (\varepsilon,t) \mapsto (\mu_t^n, \nu_{\varepsilon,t}^n, Z_t^n) \in \mathcal{S}'(\mathbb{R}^d) \times \mathcal{S}'(\mathbb{R}^d) \times \mathbb{R}^m \tag{2.24}$$

*on $(\Omega_{\mathrm{par}}^n, \mathbb{G}^n, \mathbb{P}_{x_0}^n)$ possesses continuous trajectories, up to modification. Finally,*

$$\bar{\mathbf{P}}_{x_0}^n := \mathrm{Law}_{\mathbb{P}_{x_0}^n}\left((\mu_t^n, \nu_{\varepsilon,t}^n, Z_t^n)_{(\varepsilon,t) \in [0,T]^\Delta}\right) = \mathbb{P}_{x_0}^n \circ \left((\mu_t^n, \nu_{\varepsilon,t}^n, Z_t^n)_{(\varepsilon,t) \in [0,T]^\Delta}\right)^{-1} \tag{2.25}$$

*is a well-defined measure on $(\bar{\Omega}, \bar{\mathcal{F}})$ and the sequence $(\bar{\mathbf{P}}_{x_0}^n)_{n \in \mathbb{N}}$ is tight.*

The proof of Lemma 2.6 is an application of the Kolmogorov continuity and tightness criteria. The first allows us to deduce almost sure path-regularity properties of a stochastic process taking values in a Polish space from an appropriate $\mathbb{L}^q$-bound on its increments; see for instance Kallenberg [31, Thm. 4.23]. The second lets us deduce tightness of a sequence of such processes if the $\mathbb{L}^q$-bound is uniform in the index $n$ of the sequence; see [31, Thm. 23.7]. The choice of $q$ is connected to the dimension of the index set of the processes under consideration. In the present setting, the index set is the time simplex $[0,T]^\Delta$, which is two-dimensional.

***Proof of Lemma 2.6*** Fix $r > 1$ and let $u = 1 + d/r'$ with $r'$ conjugate to $r$. Using first Lemma 2.5, then the definition of the approximated empirical measure $\nu^n$ from (1.20), the triangle inequality and finally Lemma 2.4 gives

$$\|\nu^n_{\varepsilon,t} - \nu^n_{\varepsilon',t'}\|_{\mathsf{H}_r^{-u}(\mathbb{R}^d)} \leq c(r) \sup_{\|\phi\|_{\mathsf{H}_{r'}^u} \leq 1} |(\nu^n_{\varepsilon,t} - \nu^n_{\varepsilon',t'})[\phi]|$$

$$= c(r) \sup_{\|\phi\|_{\mathsf{H}_{r'}^u} \leq 1} \left| \frac{1}{n} \sum_{i=1}^n \left( \phi(Y^{i,n;\varepsilon,t}_t) - \phi(Y^{i,n;\varepsilon',t'}_{t'}) \right) \right|$$

$$\leq c(r,\rho) \frac{1}{n} \sum_{i=1}^n |Y^{i,n;\varepsilon,t}_t - Y^{i,n}_{\varepsilon',t'}|^\rho.$$

We now take the $\mathbb{L}^q(\mathbb{P}^n_{x_0})$-norm for $q \geq 1$ chosen below, apply the triangle inequality and $\|U^\rho\|_{\mathbb{L}^q(\mathbb{P}^n_{x_0})} \leq \|U\|^\rho_{\mathbb{P}^n_{x_0}}$ because $0 < \rho < 1$, and finally use the uniform $\mathbb{L}^q(\mathbb{P}^n_{x_0})$-bound (2.1) from Proposition 2.2 to obtain for all $n \in \mathbb{N}$ and whenever $(\varepsilon,t),(\varepsilon',t') \in [0,T]^\Delta$ that

$$\left\| \|\nu^n_{\varepsilon,t} - \nu^n_{\varepsilon',t'}\|_{\mathsf{H}_r^{-u}} \right\|_{\mathbb{L}^q(\mathbb{P}^n_{x_0})} \leq c(r,\rho) \frac{1}{n} \sum_{i=1}^n \|Y^{i,n;\varepsilon,t}_t - Y^{i,n;\varepsilon',t'}_{t'}\|^\rho_{\mathbb{L}^q(\mathbb{P}^n_{x_0})}$$

$$\leq c(q,b,\sigma,\bar{\sigma},T,\beta,r,\rho)|(t-\varepsilon,t) - (t'-\varepsilon',t')|^{\frac{\rho\beta}{2}}. \quad (2.26)$$

From the definition (1.18), we have $Y^{i,n;0,t}_t = X^{i,n}_t$ and therefore by (1.20) that $\nu^n_{0,t} = \mu^n_t$. Using (2.26) for $\varepsilon = \varepsilon' = 0$ thus gives for all $n \in \mathbb{N}$ that

$$\left\| \|\mu^n_t - \mu^n_{t'}\|_{\mathsf{H}_r^{-u}(\mathbb{R}^d)} \right\|_{\mathbb{L}^q(\mathbb{P}^n_{x_0})} \leq c(q,b,\sigma,\bar{\sigma},T,\beta,r,\rho)|t - t'|^{\frac{\rho\beta}{2}}. \quad (2.27)$$

In addition, since $Z^n$ is a $(\mathbb{G}^n, \mathbb{P}^n_{x_0})$-Brownian motion, we also have

$$\|Z^n_t - Z^n_{t'}\|_{\mathbb{L}^q(\mathbb{P}^n_{x_0})} \leq c(q)|t-t'|^{\frac{1}{2}} \leq c(q)T^{\frac{1-\rho\beta}{2}}|t-t'|^{\frac{\rho\beta}{2}}. \quad (2.28)$$

We now consider the Polish space $(B, \|\cdot\|_B)$ given by

$$B_{r,u} := \mathsf{H}_r^{-u}(\mathbb{R}^d) \times \mathsf{H}_r^{-u}(\mathbb{R}^d) \times \mathbb{R}^m \quad (2.29)$$

with the product norm $\|(x_1,x_2,z)\|_B := (\|x_1\|^r_{\mathsf{H}_r^{-u}} + \|x_2\|^r_{\mathsf{H}_r^{-u}} + |z|^r)^{1/r}$. By combining the estimates (2.26)–(2.28), we then find a constant $c = c(q,b,\sigma,\bar{\sigma},T,\beta,r,\rho) < \infty$ such that for all $n \in \mathbb{N}$ and all $(\varepsilon,t),(\varepsilon',t') \in [0,T]^\Delta$, we have

$$\left\| \|(\mu^n_t,\nu^n_{\varepsilon,t},Z^n_t) - (\mu^n_{t'},\nu^n_{\varepsilon',t'},Z^n_{t'})\|_B \right\|_{\mathbb{L}^q(\mathbb{P}^n_{x_0})} \leq c|(t-\varepsilon,t) - (t'-\varepsilon',t')|^{\frac{\rho\beta}{2}}. \quad (2.30)$$

Let us now choose and fix $q > 2/(\rho\beta/2) = 4/(\rho\beta)$. For each $n \in \mathbb{N}$, we can then apply the Kolmogorov criterion in [31, Thm. 4.23] with the $\mathbb{L}^q(\mathbb{P}^n_{x_0})$-bound (2.30) to find that the process

$$(\mu^n_t, \nu^n_{t-\varepsilon,t}, Z^n_t)_{(t-\varepsilon,t)\in\{(t-\varepsilon,t):(\varepsilon,t)\in[0,T]^\Delta\}} \tag{2.31}$$

has a modification with continuous trajectories in $B$. Since the bound (2.30) is uniform over $n \in \mathbb{N}$, upon choosing a continuous modification of (2.31), the Kolmogorov tightness criterion in [31, Thm. 23.7] implies that the sequence

$$\left(\mathrm{Law}_{\mathbb{P}^n_{x_0}}\left((\mu^n_t, \nu^n_{t-\varepsilon,t}, Z^n_t)_{(t-\varepsilon,t)\in\{(t-\varepsilon,t):(\varepsilon,t)\in[0,T]^\Delta\}}\right)\right)_{n\in\mathbb{N}}$$

is tight in $\mathbf{M}_1^+(\mathbf{C}(\{(t-\varepsilon,t):(\varepsilon,t)\in[0,T]^\Delta\}; B))$.

Next, observe that the map $(t-\varepsilon,t) \mapsto (t-(t-\varepsilon),t) = (\varepsilon,t)$ is a smooth bijection between $\{(t-\varepsilon,t):(\varepsilon,t)\in[0,T]^\Delta\}$ and $[0,T]^\Delta$. Therefore, also the process

$$(\mu^n_t, \nu^n_{\varepsilon,t}, Z^n_t)_{(\varepsilon,t)\in[0,T]^\Delta}$$

has a modification with continuous trajectories in $B$, and so

$$\bar{\mathbf{P}}^n_{x_0} := \mathrm{Law}_{\mathbb{P}^n_{x_0}}\left((\mu^n_t, \nu^n_{\varepsilon,t}, Z^n_t)_{(\varepsilon,t)\in[0,T]^\Delta}\right) = \mathbb{P}^n_{x_0} \circ \left((\mu^n_t, \nu^n_{\varepsilon,t}, Z^n_t)_{(\varepsilon,t)\in[0,T]^\Delta}\right)^{-1} \tag{2.32}$$

is a well-defined probability measure on $\bar{\Omega} = \mathbf{C}([0,T]^\Delta; B)$ and $(\bar{\mathbf{P}}^n_{x_0})_{n\in\mathbb{N}}$ is tight in $\mathbf{M}_1^+(\mathbf{C}([0,T]^\Delta; B_{r,u}))$. Finally, the inclusion $B \hookrightarrow \mathcal{S}'(\mathbb{R}^d) \times \mathcal{S}'(\mathbb{R}^d) \times \mathbb{R}^m$ is continuous so that the sequence in (2.31) is tight in $\mathbf{M}_1^+(\mathbf{C}([0,T]^\Delta; \mathcal{S}'(\mathbb{R}^d) \times \mathcal{S}'(\mathbb{R}^d) \times \mathbb{R}^m))$ as well. $\square$

## 2.3 Construction of the law of the decomposition

We defined in Sections 1.3 and 1.5 the canonical space $\Omega = \mathbf{C}([0,T]; \mathcal{S}')$ and its extension

$$\bar{\Omega} = \mathbf{C}\left([0,T]^\Delta; \mathcal{S}'(\mathbb{R}^d) \times \mathcal{S}'(\mathbb{R}^d) \times \mathbb{R}^m\right),$$

which we equipped with its Borel-$\sigma$-algebra $\bar{\mathcal{F}}$. From the definition (2.25) of $\bar{\mathbf{P}}^n_{x_0}$ and the definition (1.14) of the canonical process $\bar{\Omega}$, we have

$$\mathrm{Law}_{\bar{\mathbf{P}}^n_{x_0}}(\bar{\Lambda}) = \mathrm{Law}_{\mathbf{P}^n_{x_0}}(\Lambda) = \mathrm{Law}_{\mathbb{P}^n_{x_0}}(\mu^n), \quad \mathrm{Law}_{\bar{\mathbf{P}}^n_{x_0}}(Z) = \mathrm{Law}_{\mathbb{P}^n_{x_0}}(Z^n) \tag{2.33}$$

and in particular also

$$\mathrm{Law}_{\bar{\mathbf{P}}^n_{x_0}}(A) = \mathrm{Law}_{\mathbb{P}^n_{x_0}}(\nu^n). \tag{2.34}$$

In Theorem 1.5, we consider the cluster points of the sequence $(\mathbf{P}^n_{x_0})_{n\in\mathbb{N}}$. However, we now have the additional sequence $(\bar{\mathbf{P}}^n_{x_0})_{n\in\mathbb{N}}$ on the extended space, and this sequence if tight by Lemma 2.6. To connect these two objects, we use the following two-stage construction. First,

*fix a cluster point $\mathbf{P}^\infty_{x_0}$ of the sequence $(\mathbf{P}^n_{x_0})_{n\in\mathbb{N}}$ and* (2.35)
*a subsequence $(n_k)_{k\in\mathbb{N}}$ such that $\mathbf{P}^{n_k}_{x_0} \to \mathbf{P}^\infty_{x_0}$ narrowly as $k \to \infty$.*

By Lemma 2.6, the corresponding sequence $(\bar{\mathbf{P}}_{x_0}^{n_k})_{k \in \mathbb{N}}$ is again tight so that in a second step, we subsequently

*fix a cluster point $\bar{\mathbf{P}}_{x_0}^{\infty}$ of the subsequence $(\bar{\mathbf{P}}_{x_0}^{n_k})_{k \in \mathbb{N}}$ and* (2.36)
*a subsubsequence $(n_{k'})_{k' \in \mathbb{N}}$ of $(n_k)_{k \in \mathbb{N}}$ such that $\bar{\mathbf{P}}_{x_0}^{n_{k'}} \to \bar{\mathbf{P}}_{x_0}^{\infty}$ narrowly as $k' \to \infty$.*

Note that by construction, $\bar{\Lambda}$ under $\bar{\mathbf{P}}_{x_0}^{\infty}$ has the same law as $\Lambda$ under the cluster point $\mathbf{P}_{x_0}^{\infty}$ that we started out with in (2.35), i.e.,

$$\mathrm{Law}_{\bar{\mathbf{P}}_{x_0}^{\infty}}(\bar{\Lambda}) = \mathrm{Law}_{\mathbf{P}_{x_0}^{\infty}}(\Lambda).$$ (2.37)

The measure $\bar{\mathbf{P}}_{x_0}^{\infty}$ on $(\bar{\Omega}, \bar{\mathcal{F}})$ is the lifting we sought. It lets us study $\Lambda$ via $\bar{\Lambda}$ using its decomposition into $A_{\varepsilon,t} + E_{\varepsilon,t}$ from (1.16). This completes Step 2 of our proof.

## 3 THE REGULAR PART

In Sections 1 and 2, we constructed on the extended space $(\bar{\Omega}, \bar{\mathcal{F}}, \bar{\mathbf{P}}_{x_0}^{\infty})$ a decomposition of $\bar{\Lambda}$ into a regular part $A$ and a distributional part $E$. The goal of this section is to obtain quantitative norm estimates of $A$. For this, we first derive a probabilistic representation for $A$ under $\bar{\mathbf{P}}_{x_0}^{\infty}$ which is a narrow cluster point of the sequence $(\bar{\mathbf{P}}_{x_0}^{n})_{n \in \mathbb{N}}$. This is achieved using the identity $\mathrm{Law}_{\bar{\mathbf{P}}_{x_0}^{n}}(A) = \mathrm{Law}_{\mathbb{P}_{x_0}^{n}}(\nu^n)$ for $n \in \mathbb{N}$ from (2.34), and the approximations (1.18) and (1.20) which are sufficiently simple to allow a semi-explicit formula for $A$ as a Gaussian mixture in terms of $\Lambda$ and $Z$; see Proposition 3.1. We then use analytic means to obtain from this representation of $A$ a quantitative norm bound.

### 3.1 Probabilistic representation

To state the representation result, let us define for each $(\varepsilon, t) \in [0, T]^{\Delta}$ the maps

$$M_{\varepsilon,t} : \mathbb{R}^d \times \mathbf{M}_1^+(\mathbb{R}^d) \times \mathbf{C}([0,T]; \mathbb{R}^d) \to \mathbb{R}^d,$$

$$(x, \lambda, z) \mapsto M_{\varepsilon,t}(x, \lambda, z) := x + \bar{\sigma}_{t-\varepsilon}(x, \lambda)(z_t - z_{t-\varepsilon}),$$ (3.1)

and

$$\Sigma_{\varepsilon,t} : \mathbb{R}^d \times \mathbf{M}_1^+(\mathbb{R}^d) \to \mathbb{R}^{d \times d},$$

$$(x, \lambda) \mapsto \Sigma_{\varepsilon,t}(x, \lambda) := \int_{t-\varepsilon}^{t} (\sigma_s \sigma_s^{\top})(x, \lambda) \, \mathrm{d}s.$$ (3.2)

Once we equip $\mathbf{M}_1^+ \subseteq \mathcal{S}'$ with the topology inherited from $\mathcal{S}'$, and the respective domains with their product topologies, these maps are immediately seen to be Borel-measurable since $\sigma, \bar{\sigma}$ are measurable in view of Assumption 1.2.

Denote by $\mathcal{N}(m, v)$ the $d$-dimensional normal law with mean $m \in \mathbb{R}^d$ and covariance matrix $v \in \mathbb{R}^{d \times d}$, and by $g(\,\cdot\,; m, v)$ its density on $\mathbb{R}^d$ relative to Lebesgue measure.

**Proposition 3.1** | *Under Assumptions 1.1 and 1.2, let $\bar{\mathbf{P}}_{x_0}^\infty$ be a cluster point of $(\bar{\mathbf{P}}_{x_0}^n)_{n \in \mathbb{N}}$. For each fixed pair $(\varepsilon, t) \in [0, T]^\Delta$, we have $\bar{\mathbf{P}}_{x_0}^\infty$-a.s. the representation*

$$A_{\varepsilon,t}[\phi] = \int_{\mathbb{R}^d} \mathbb{E}\Big[\phi\big(U_{\Lambda_{t-\varepsilon}, Z_{\cdot \wedge t}}(x)\big)\Big] \Lambda_{t-\varepsilon}(\mathrm{d}x) \quad \textit{for all } \phi \in \mathcal{S}\,. \tag{3.3}$$

*In* (3.3), $U_{\Lambda_{t-\varepsilon}, Z_{\cdot \wedge t}}(x)$ *is a normal random variable in $\mathbb{R}^d$ which has mean $M_{\varepsilon,t}(x, \Lambda_{t-\varepsilon}, Z_{\cdot \wedge t})$ as in* (3.1) *and variance $\Sigma_{\varepsilon,t}(x, \Lambda_{t-\varepsilon})$ as in* (3.2), *more compactly*

$$U_{\Lambda_{t-\varepsilon}, Z_{\cdot \wedge t}} \stackrel{(\mathrm{d})}{=} \mathcal{N}\big(M_{\varepsilon,t}(x, \Lambda_{t-\varepsilon}, Z_{\cdot \wedge t}), \Sigma_{\varepsilon,t}(x, \Lambda_{t-\varepsilon})\big)\,,$$

*and $\mathbb{E}$ denotes the expectation with respect to this normal law.*

Before giving the proof of Proposition 3.1, let us develop a heuristic understanding of why we should expect the mixture representation (3.3). Suppose that $\Lambda_{t-\varepsilon}$ charges a small set of size $\mathrm{d}x$ with positive $\bar{\mathbf{P}}_{x_0}^\infty$-probability. Since $\mu_{t-\varepsilon}^n$ converges in law to $\Lambda_{t-\varepsilon}$ as $n \to \infty$, we expect that $\mu_{t-\varepsilon}^n[\mathrm{d}x] = \frac{1}{n}\#\{i : X_{t-\varepsilon}^{i,n} \in (x, x+\mathrm{d}x]\}$ should tend to $\Lambda_{t-\varepsilon}[(x, x+\mathrm{d}x]] > 0$ as $n \to \infty$. In particular, $\#\{i : X_{t-\varepsilon}^{i,n} \in (x, x+\mathrm{d}x]\}$, the number of particles in $(x, x+\mathrm{d}x]$, must tend to infinity with positive $\bar{\mathbf{P}}_{x_0}^\infty$-probability. Consider now the approximating particles $Y_t^{i,n;\varepsilon,t}$ from (1.18) starting from within $(x, x+\mathrm{d}x]$ at time $t-\varepsilon$. By time $t$, they will all have been perturbed by the common component $\bar{\sigma}_{t-\varepsilon}(x, \Lambda_{t-\varepsilon})(Z_t^n - Z_{t-\varepsilon}^n)$. Moreover, conditionally on $Z^n$ and $\Lambda_{t-\varepsilon}$, the approximating particles are independent and perturbed only by their idiosyncratic Brownian increments $\int_{t-\varepsilon}^t \sigma_r(x, \Lambda_{t-\varepsilon})\,\mathrm{d}B_r^{i,n}$. In the weak limit, we then expect a conditional law of large numbers to be valid for the particles that started in $(x, x+\mathrm{d}x]$, giving rise to the Gaussian law with mean and covariance $M_{\varepsilon,t}(x, \Lambda_{t-\varepsilon}, Z_{\cdot \wedge t}) = x + \bar{\sigma}_{t-\varepsilon}(x, \Lambda_{t-\varepsilon})(Z_t^n - Z_{t-\varepsilon}^n)$ and $\Sigma_{\varepsilon,t}(x, \Lambda_{t-\varepsilon}) = \int_{t-\varepsilon}^t (\sigma_r \sigma_r^\top)(x, \Lambda_{t-\varepsilon})\,\mathrm{d}r$, respectively. These are precisely the functions in (3.1) and (3.2). The mixture is then obtained by averaging the local $\mathrm{d}x$-snapshot with respect to the probability measure $\Lambda_{t-\varepsilon}$, as in (3.3).

*Proof of Proposition 3.1* The proof consists of several steps. We first derive an appropriate representation for the approximating finite particle system (1.18), then verify that it carries over to the empirical measure (1.20), and finally show that it remains valid in the limit.

For the pair $(0, t) \in [0, T]^\Delta$, (3.1) and (3.2) give $M_{0,t}(x, \lambda, z) = x$ and $\Sigma_{0,t}(x, \lambda) = 0$ so that $U_{\Lambda_t, Z_{\cdot \wedge t}} \sim \mathcal{N}(x, 0) = \delta_x$. Hence (3.3) becomes the assertion that

$$A_{0,t}[\phi] = \int_{\mathbb{R}^d} \phi(x) \Lambda_t(\mathrm{d}x) \qquad \bar{\mathbf{P}}_{x_0}^\infty\text{-a.s.} \tag{3.4}$$

But for each $n \in \mathbb{N}$, (1.18) and (1.20) give $\mu_t^n = \nu_{0,t}^n$, and so the identity (2.34) implies $\mathrm{Law}_{\bar{\mathbf{P}}_{x_0}^n}(A_{0,t}) = \mathrm{Law}_{\mathbf{P}_{x_0}^n}(\Lambda_t) = \mathrm{Law}_{\mathbb{P}_{x_0}^n}(\mu_t^n)$. Now (3.4) follows by passing to the limit $\bar{\mathbf{P}}_{x_0}^\infty$ of an appropriate subsequence $(\bar{\mathbf{P}}_{x_0}^{n_k})_{k \in \mathbb{N}}$.

For the remainder of the proof, we fix $(\varepsilon, t)$ with $0 < \varepsilon \leq t \leq T$. Consider a fixed $n \in \mathbb{N}$. We take the limit only in Step 6 below.

**Step 1** We first derive a representation for the full $n$-particle system (1.18) under $\mathbb{P}_{x_0}^n$. For each $\phi \in \mathcal{S}$, the definition (1.20) of $\nu_{\varepsilon,t}^n$ gives

$$\nu_{\varepsilon,t}^n[\phi] = \frac{1}{n}\sum_{i=1}^n \phi(Y_t^{i,n;\varepsilon,t})\,, \tag{3.5}$$

which implies

$$\mathbf{E}^{\mathbb{P}^n_{x_0}}\left[\nu^n_{\varepsilon,t}[\phi]\,\Big|\,\mu^n_{t-\varepsilon}, Z^n_{\cdot\wedge t}\right] = \mathbf{E}^{\mathbb{P}^n_{x_0}}\left[\frac{1}{n}\sum_{i=1}^n \phi(Y^{i,n;\varepsilon,t}_t)\,\Big|\,\mu^n_{t-\varepsilon}, Z^n_{\cdot\wedge t}\right] \qquad \mathbb{P}^n_{x_0}\text{-a.s.} \qquad (3.6)$$

With the notation $X^n_{t-\varepsilon} = (X^{1,n}_{t-\varepsilon}, \ldots, X^{n,n}_{t-\varepsilon})$ from Section 1.2 and noting that the $\sigma$-field generated by $\mu^n_{t-\varepsilon}$ is contained in that generated by $X^n_{t-\varepsilon}$, we can rewrite the above as

$$\mathbf{E}^{\mathbb{P}^n_{x_0}}\left[\frac{1}{n}\sum_{i=1}^n \mathbf{E}^{\mathbb{P}^n_{x_0}}[\phi(Y^{i,n;\varepsilon,t}_t)\,|\,X^n_{t-\varepsilon}, Z^n_{\cdot\wedge t}]\,\Big|\,\mu^n_{t-\varepsilon}, Z^n_{\cdot\wedge t}\right] \qquad \mathbb{P}^n_{x_0}\text{-a.s.}$$

Using the definition of $Y^{i,n;\varepsilon,t}_t$ in (1.18), we note that conditionally on $X^n_{t-\varepsilon}, Z^n_{\cdot\wedge t}$, the vector $(Y^{1,n;\varepsilon,t}_t, \ldots, Y^{n,n;\varepsilon,t}_t)$ is under $\mathbb{P}^n_{x_0}$ jointly normal with conditionally independent components. Moreover, for each $i \in [n]$, the random variable $Y^{i,n;\varepsilon,t}_t$ conditionally on $X^{i,n}_{t-\varepsilon}, \mu^n_{t-\varepsilon}$ and $Z^n_{\cdot\wedge t}$ is under $\mathbb{P}^n_{x_0}$ independent of $X^{j,n}_{t-\varepsilon}$ for all $j \neq i$ because $B^{1,n}, \ldots, B^{n,n}$ are independent and independent of $Z^n$. We thus have $\mathbb{P}^n_{x_0}$-a.s. for each $i \in [n]$ that

$$\text{Law}_{\mathbb{P}^n_{x_0}}(Y^{i,n;\varepsilon,t}_t \mid X^n_{t-\varepsilon}, Z^n_{\cdot\wedge t}) = \text{Law}_{\mathbb{P}^n_{x_0}}(Y^{i,n;\varepsilon,t}_t \mid X^{i,n}_{t-\varepsilon}, \mu^n_{t-\varepsilon}, Z^n_{\cdot\wedge t})$$
$$= \mathcal{N}\Big(M_{\varepsilon,t}(X^{i,n}_{t-\varepsilon}, \mu^n_{t-\varepsilon}, Z_{\cdot\wedge t}), \Sigma_{\varepsilon,t}(X^{i,n}_{t-\varepsilon}, \mu^n_{t-\varepsilon})\Big). \qquad (3.7)$$

Define for $\phi \in \mathcal{S}$ and $(\lambda, z) \in \mathbf{M}^+_1(\mathbb{R}^d) \times \mathbf{C}([0,T]; \mathbb{R}^m)$ the map $F^\phi_{\varepsilon,t}(\,\cdot\,; \lambda, z) : \mathbb{R}^d \to \mathbb{R}$ by

$$F^\phi_{\varepsilon,t}(x; \lambda, z) := \int_{\mathbb{R}^d} \phi(y)\, g\Big(y; M_{\varepsilon,t}(x, \lambda, z), \Sigma_{\varepsilon,t}(x, \lambda)\Big)\,\mathrm{d}y, \qquad (3.8)$$

where we recall that $g(\,\cdot\,; m, v)$ is the density function of $\mathcal{N}(m, v)$. With this notation and (3.7), we see that $\mathbb{P}^n_{x_0}$-almost surely,

$$\mathbf{E}^{\mathbb{P}^n_{x_0}}\left[\frac{1}{n}\sum_{i=1}^n \phi(Y^{i,n;\varepsilon,t}_t)\,\Big|\,X^n_{t-\varepsilon}, Z^n_{\cdot\wedge t}\right] = \frac{1}{n}\sum_{i=1}^n F^\phi_{\varepsilon,t}(X^{i,n}_{t-\varepsilon}; \mu^n_{t-\varepsilon}, Z^n_{\cdot\wedge t}).$$

On the other hand, using the definition of the empirical measure $\mu^n_{t-\varepsilon}$ from (1.19), the last expression can be written as

$$\mathbf{E}^{\mathbb{P}^n_{x_0}}\left[\frac{1}{n}\sum_{i=1}^n \phi(Y^{i,n;\varepsilon,t}_t)\,\Big|\,X^n_{t-\varepsilon}, Z^n_{\cdot\wedge t}\right] = \int_{\mathbb{R}^d} F^\phi_{\varepsilon,t}(x; \mu^n_{t-\varepsilon}, Z^n_{\cdot\wedge t})\,\mu^n_{t-\varepsilon}(\mathrm{d}x). \qquad (3.9)$$

Note that this last term is measurable relative to the $\sigma$-algebra $\sigma(\mu^n_{t-\varepsilon}, Z^n_{\cdot\wedge t}) \subseteq \sigma(X^n_{t-\varepsilon}, Z^n_{\cdot\wedge t})$. Combining (3.6) with (3.9) then gives that $\mathbb{P}^n_{x_0}$-almost surely,

$$\mathbf{E}^{\mathbb{P}^n_{x_0}}\left[\nu^n_{\varepsilon,t}[\phi]\,\Big|\,\mu^n_{t-\varepsilon}, Z^n_{\cdot\wedge t}\right] = \int_{\mathbb{R}^d} F^\phi_{\varepsilon,t}(x; \mu^n_{t-\varepsilon}, Z^n_{\cdot\wedge t})\,\mu^n_{t-\varepsilon}(\mathrm{d}x). \qquad (3.10)$$

**Step 2** We next derive a concentration bound via Hoeffding's inequality. Recall that the family $\{Y^{1,n;\varepsilon,t}_t, \ldots, Y^{n,n;\varepsilon,t}_t\}$ is under $\mathbb{P}^n_{x_0}$ independent conditionally on $X^n_{t-\varepsilon}, Z^n_{\cdot\wedge t}$. Now define

$$\tilde{D}^n_{\delta,\phi} := \left\{\left|\frac{1}{n}\sum_{i=1}^n \Big(\phi(Y^{i,n;\varepsilon,t}_t) - \mathbf{E}^{\mathbb{P}^n_{x_0}}[\phi(Y^{i,n;\varepsilon,t}_t)\,|\,X^n_{t-\varepsilon}, Z^n_{\cdot\wedge t}]\Big)\right| > \delta\right\}$$

for $\delta > 0$. Note that $\phi(Y_t^{i,n;\varepsilon,t})$ for $i \in [n]$ are random variables that are bounded by $\|\phi\|_\infty$. Then Hoeffding's bound gives the concentration inequality

$$\mathbb{P}_{x_0}^n[\tilde{D}_{\delta,\phi}^n \mid X_{t-\varepsilon}^n, Z_{\cdot \wedge t}^n] \leq \exp\left(-\frac{n\delta^2}{\|\phi\|_\infty}\right) \qquad \mathbb{P}_{x_0}^n\text{-a.s.} \qquad (3.11)$$

We remark that this bound is independent of $\varepsilon$ since also $\|\phi\|_\infty$ is independent of $\varepsilon$, and this is why we chose to make the notation $\tilde{D}_{\delta,\phi}^n$ not depend on $\varepsilon$ either. Letting

$$D_{\delta,\phi}^n := \left\{ \left| \nu_{\varepsilon,t}^n[\phi] - \int_{\mathbb{R}^d} F_{\varepsilon,t}^\phi(x; \mu_{t-\varepsilon}^n, Z_{\cdot \wedge t}^n) \, \mu_{t-\varepsilon}^n(\mathrm{d}x) \right| > \delta \right\} \qquad (3.12)$$

and in view of (3.5) and (3.10), we have $\mathbb{P}_{x_0}^n$-a.s. that $D_{\delta,\phi}^n = \tilde{D}_{\delta,\phi}^n$. Together with the tower property of conditional expectations and (3.11), this gives

$$\mathbb{P}_{x_0}^n[D_{\delta,\phi}^n \mid \mu_{t-\varepsilon}^n, Z_{\cdot \wedge t}^n] = \mathbf{E}^{\mathbb{P}_{x_0}^n}\left[\mathbb{P}_{x_0}^n[\tilde{D}_{\delta,\phi}^n \mid X_{t-\varepsilon}^n, Z_{\cdot \wedge t}^n] \,\Big|\, \mu_{t-\varepsilon}^n, Z_{\cdot \wedge t}^n\right]$$

$$\leq \exp\left(-\frac{n\delta^2}{\|\phi\|_\infty}\right) \qquad \mathbb{P}_{x_0}^n\text{-a.s.} \qquad (3.13)$$

This is our concentration inequality for the approximated empirical measure.

**Step 3** To transfer the bound obtained in Step 2 from $\mathbb{P}_{x_0}^n$ to $\bar{\mathbf{P}}_{x_0}^n$, we note that all the quantities appearing in (3.12) are expressible in terms of the canonical coordinates on the space $\bar{\Omega}$. Indeed, in view of (3.12) and (2.34), define

$$D_{\delta,\phi} := \left\{ \left| A_{\varepsilon,t}[\phi] - \int_{\mathbb{R}^d} F_{\varepsilon,t}^\phi(x; \Lambda_{t-\varepsilon}, Z_{\cdot \wedge t}) \, \Lambda_{t-\varepsilon}(\mathrm{d}x) \right| > \delta \right\}. \qquad (3.14)$$

Then (3.13) yields in view of (2.33) that

$$\bar{\mathbf{P}}_{x_0}^n[D_{\delta,\phi} \mid \Lambda_{t-\varepsilon}, Z_{\cdot \wedge t}] = \mathbb{P}_{x_0}^n[D_{\delta,\phi}^n \mid \Lambda_{t-\varepsilon}, Z_{\cdot \wedge t}] \leq \exp\left(-\frac{n\delta^2}{\|\phi\|_\infty}\right) \qquad \bar{\mathbf{P}}_{x_0}^n\text{-a.s.} \qquad (3.15)$$

Since $n$ in (3.13) is arbitrary, (3.15) is valid for all $n \in \mathbb{N}$.

**Step 4** Recall $F_{\varepsilon,t}^\phi$ from (3.8). We claim that the map

$$(\lambda, z) \mapsto \int_{\mathbb{R}^d} F_{\varepsilon,t}^\phi(x; \lambda, z_{\cdot \wedge t}) \, \lambda(\mathrm{d}x) \qquad (3.16)$$

is continuous on $\mathcal{P}_{\mathrm{wk}^*}(\mathbb{R}^d) \times \mathbf{C}([0,T]; \mathbb{R}^m)$, where we recall that $\mathcal{P}_{\mathrm{wk}^*}(\mathbb{R}^d)$ denotes the space of probability measures on $\mathcal{B}(\mathbb{R}^d)$ equipped with narrow convergence. To argue this, let $(\lambda^k, z^k) \to (\lambda, z)$ in $\mathcal{P}_{\mathrm{wk}^*}(\mathbb{R}^d) \times \mathbf{C}([0,T]; \mathbb{R}^m)$. Then

$$\limsup_{k \to \infty} \left| \int_{\mathbb{R}^d} F_{\varepsilon,t}^\phi(x; \lambda^k, z^k) \, \lambda^k(\mathrm{d}x) - \int_{\mathbb{R}^d} F_{\varepsilon,t}^\phi(x; \lambda, z) \, \lambda(\mathrm{d}x) \right|$$

$$\leq \limsup_{k \to \infty} \left| \int_{\mathbb{R}^d} F_{\varepsilon,t}^\phi(x; \lambda^k, z^k) \left(\lambda^k(\mathrm{d}x) - \lambda(\mathrm{d}x)\right) \right|$$

$$+ \limsup_{k \to \infty} \left| \int_{\mathbb{R}^d} \left( F_{\varepsilon,t}^\phi(x; \lambda^k, z^k) - F_{\varepsilon,t}^\phi(x; \lambda, z) \right) \lambda(\mathrm{d}x) \right|. \qquad (3.17)$$

For the second summand, note that by Assumption 1.2(H), the maps $(\lambda, x) \mapsto M_{\varepsilon,t}(x, \lambda, z)$ and $\lambda \mapsto \Sigma_\varepsilon(x, \lambda)$ from (3.1) and (3.2) are continuous for all $(\varepsilon, t) \in [0, T]^\Delta$ and $x \in \mathbb{R}^d$. Hence so is the map $(\lambda, z) \mapsto F^\phi_{\varepsilon,t}(x; \lambda, z)$ from (3.8). We thus have for all $x \in \mathbb{R}^d$ that $\lim_{k \to \infty} F^\phi_{\varepsilon,t}(x; \lambda^k, z^k) = F^\phi_{\varepsilon,t}(x; \lambda, z)$, and so dominated convergence gives

$$\lim_{k \to \infty} \left| \int_{\mathbb{R}^d} \left( F^\phi_{\varepsilon,t}(x; \lambda^k, z^k) - F^\phi_{\varepsilon,t}(x; \lambda, z) \right) \lambda(\mathrm{d}x) \right| = 0.$$

We now turn to the first summand in (3.17). Since $(x, \lambda) \mapsto \sigma_t(x, \lambda)$ and $(x, \lambda) \mapsto \bar{\sigma}_t(x, \lambda)$ are Hölder-continuous functions by Assumption 1.2(H), the definitions (3.8) of $F^\phi_{\varepsilon,t}$, and (3.1) and (3.2) of $M_{\varepsilon,t}$ and $\Sigma_{\varepsilon,t}$ imply that $\mathcal{C} := \{x \mapsto F^\phi_{\varepsilon,t}(x; \lambda^k, z^k) : k \in \mathbb{N}\} \subseteq \mathbf{C}_b(\mathbb{R}^d)$ is a family of functions that is equicontinuous at each point $x \in \mathbb{R}^d$. But because $\lambda^k \to \lambda$ in $\mathcal{P}_{\mathrm{wk}^*}(\mathbb{R}^d)$ as $k \to \infty$, the convergence result of Stroock and Varadhan [49, Cor. 1.1.5] shows that

$$\limsup_{k \to \infty} \left| \int_{\mathbb{R}^d} F^\phi_{\varepsilon,t}(x; \lambda^k, z^k) \left( \lambda^k(\mathrm{d}x) - \lambda(\mathrm{d}x) \right) \right|$$
$$\leq \limsup_{k \to \infty} \sup_{F \in \mathcal{C}} \left| \int_{\mathbb{R}^d} F(x) \left( \lambda^k(\mathrm{d}x) - \lambda(\mathrm{d}x) \right) \right| = 0.$$

In conclusion, both summands on the right-hand side of (3.17) vanish, showing the continuity of (3.16) as claimed.

**Step 5** Let $(n_k)_{k \in \mathbb{N}}$ be a subsequence such that $(\bar{\mathbf{P}}^{n_k}_{x_0})_{k \in \mathbb{N}}$ converges to $\bar{\mathbf{P}}^\infty_{x_0}$ narrowly in $\mathbf{M}^+_1(\bar{\Omega})$ as $k \to \infty$. The space $\bar{\Omega}$ is endowed with a complicated topology and not separable. For technical reasons, we therefore confine ourselves to a more convenient subspace. For $r \in (1, \infty)$ and $u \in \mathbb{R}$, define $\bar{\Omega}_{r,u} := \mathbf{C}([0, T]^\Delta; B_{r,u})$ with $B_{r,u}$ is defined in (2.29). Clearly $\bar{\Omega}_{r,u} \hookrightarrow \bar{\Omega}$. It follows from the proof of Lemma 2.6 and in particular from (2.32) that there exist $r$ and $u$ such that $\mathrm{Law}_{\mathbb{P}^n_{x_0}}((\mu^n_t, \nu^n_{\varepsilon,t}, Z^n_t)_{(\varepsilon,t) \in [0,T]^\Delta})$ is in $\mathbf{M}^+_1(\Omega_{r,u})$ for all $n \in \mathbb{N}$, where $B$ is defined in (2.29). We can therefore assume that each $\bar{\mathbf{P}}^{n_k}_{x_0}$ as well as $\bar{\mathbf{P}}^\infty_{x_0}$ puts all mass on $\bar{\Omega}_{r,u}$, i.e., that $\bar{\mathbf{P}}^{n_k}_{x_0}[\bar{\Omega}_{r,u}] = 1$ for all $k \in \mathbb{N} \cup \{\infty\}$, where we set $n_\infty = \infty$. Since $\bar{\Omega}_{r,u}$ is Polish, we can apply the Skorohod representation theorem; see e.g. Pollard [44, Thm. IV.13]. We therefore find on an appropriate probability space $(\tilde{\Omega}, \tilde{\mathcal{F}}, \tilde{\mathbf{P}})$ a family of $\bar{\Omega}_{r,u}$-valued random variables $(\tilde{\Lambda}^k, \tilde{A}^k, \tilde{Z}^k)$ for $k \in \mathbb{N} \cup \{\infty\}$ which satisfy

$$\mathrm{Law}_{\tilde{\mathbf{P}}}(\tilde{\Lambda}^k, \tilde{A}^k, \tilde{Z}^k) = \mathrm{Law}_{\bar{\mathbf{P}}^{n_k}_{x_0}}(\bar{\Lambda}, A, Z) \qquad \text{for all } k \in \mathbb{N} \cup \{\infty\} \qquad (3.18)$$

and

$$\lim_{k \to \infty} (\tilde{\Lambda}^k, \tilde{A}^k, \tilde{Z}^k) = (\tilde{\Lambda}^\infty, \tilde{A}^\infty, \tilde{Z}^\infty) \qquad \text{in } \bar{\Omega}_{r,u}, \, \tilde{\mathbf{P}}\text{-a.s.} \qquad (3.19)$$

Since the space $\mathbf{C}^\infty_c(\mathbb{R}^d)$ of smooth and compactly supported functions is separable and dense in $\mathbf{H}^u_r(\mathbb{R}^d)$, we can find a $\tilde{\mathbf{P}}$-nullset $N \in \tilde{\mathcal{F}}$ such that for all $\tilde{\omega} \in \tilde{\Omega} \setminus N$, $\tilde{\Lambda}^k(\tilde{\omega})$ is $\mathbf{M}^+_1(\mathbb{R}^d)$-valued for all $k \in \mathbb{N} \cup \{\infty\}$, we have

$$\lim_{\delta \to 0} \sup_{\substack{t, t' \in [0,T] \\ |t - t'| < \delta}} |\tilde{\Lambda}^k_t(\tilde{\omega})[\phi] - \tilde{\Lambda}^k_{t'}(\tilde{\omega})[\phi]| = 0 \qquad \text{for all } \phi \in \mathbf{C}^\infty_c(\mathbb{R}^d) \text{ and } k \in \mathbb{N} \cup \{\infty\}, \qquad (3.20)$$

and

$$\lim_{k \to \infty} \sup_{t \in [0,T]} |\tilde{\Lambda}^k_t(\tilde{\omega})[\phi] - \tilde{\Lambda}^\infty_t(\tilde{\omega})[\phi]| = 0 \qquad \text{for all } \phi \in \mathbf{C}^\infty_c(\mathbb{R}^d). \qquad (3.21)$$

In addition, from the concentration bounds (1.9) and (1.10), we deduce that we can find a $\tilde{\mathbf{P}}$-nullset $\tilde{N} \supseteq N$ such that for each $\varepsilon > 0$, $\tilde{\omega} \in \tilde{\Omega} \setminus \tilde{N}$ and $k \in \mathbb{N} \cup \{\infty\}$, there exists a function $\phi_K \in \mathbf{C}_c^\infty(\mathbb{R}^d)$ with $0 \leq \phi_K \leq 1$ and $\phi_K = 1$ on $[-K, K]^d$ satisfying

$$\sup_{t \in [0,T]} \tilde{\Lambda}_t^k(\tilde{\omega})[1 - \phi_K] < \varepsilon. \tag{3.22}$$

We first claim that for each $k \in \mathbb{N} \cup \{\infty\}$ and $\tilde{\omega} \in \tilde{\Omega} \setminus \tilde{N}$, $\tilde{\Lambda}^k(\tilde{\omega})$ is in $\mathbf{C}([0,T]; \mathcal{P}_{\mathrm{wk}^*}(\mathbb{R}^d))$. To argue this, let $f \in \mathbf{C}_b(\mathbb{R}^d)$. Then for $\phi_K$ as in (3.22) and $t, t' \in [0, T]$, we get

$$\left|\left(\tilde{\Lambda}_t^k(\tilde{\omega}) - \tilde{\Lambda}_{t'}^k(\tilde{\omega})\right)[f]\right| \leq \left|\left(\tilde{\Lambda}_t^k(\tilde{\omega}) - \tilde{\Lambda}_{t'}^k(\tilde{\omega})\right)[\phi_K f]\right| + \left|\left(\tilde{\Lambda}_t^k(\tilde{\omega}) - \tilde{\Lambda}_{t'}^k(\tilde{\omega})\right)[(1 - \phi_K)f]\right|.$$

Since $\phi_K f$ is in $\mathbf{C}_b(\mathbb{R}^d)$ and compactly supported, we can find a function $\psi \in \mathbf{C}_c^\infty(\mathbb{R}^d)$ with $\|\phi_K f - \psi\|_\infty < \varepsilon$. Using this and then appealing once more to (3.22), we obtain

$$\left|\left(\tilde{\Lambda}_t^k(\tilde{\omega}) - \tilde{\Lambda}_{t'}^k(\tilde{\omega})\right)[f]\right| \leq 2\varepsilon + \left|\left(\tilde{\Lambda}_t^k(\tilde{\omega}) - \tilde{\Lambda}_{t'}^k(\tilde{\omega})\right)[\psi]\right| + \left|\left(\tilde{\Lambda}_t^k(\tilde{\omega}) - \tilde{\Lambda}_{t'}^k(\tilde{\omega})\right)[(1 - \phi_K)f]\right|$$
$$\leq 2\varepsilon + \left|\left(\tilde{\Lambda}_t^k(\tilde{\omega}) - \tilde{\Lambda}_{t'}^k(\tilde{\omega})\right)[\psi]\right| + 2\varepsilon\|f\|_\infty.$$

In view of (3.20) and since $\varepsilon > 0$ is arbitrary, we see that for each $f \in \mathbf{C}_b(\mathbb{R}^d)$, the function $t \mapsto \tilde{\Lambda}_t^k(\tilde{\omega})[f]$ is uniformly continuous on $[0, T]$. Since $f \in \mathbf{C}_b(\mathbb{R}^d)$ is arbitrary we deduce that $\tilde{\Lambda}^k(\tilde{\omega}) \in \mathbf{C}([0,T]; \mathcal{P}_{\mathrm{wk}^*}(\mathbb{R}^d))$. This gives the claim.

We next claim that for each $\tilde{\omega} \in \tilde{\Omega} \setminus \tilde{N}$, $\tilde{\Lambda}^k(\tilde{\omega}) \to \tilde{\Lambda}^\infty(\tilde{\omega})$ in $\mathbf{C}([0,T]; \mathcal{P}_{\mathrm{wk}^*}(\mathbb{R}^d))$ as $k \to \infty$. The argument is similar to before. Proceeding as above choosing first $\phi_K$ as in (3.22) for $k = \infty$ and then $\psi \in \mathbf{C}_c^\infty(\mathbb{R}^d)$ with $\|\phi_K f - \psi\|_\infty < \varepsilon$ gives

$$\left|\left(\tilde{\Lambda}_t^k(\tilde{\omega}) - \tilde{\Lambda}_t^\infty(\tilde{\omega})\right)[f]\right| \leq 2\varepsilon + \left|\left(\tilde{\Lambda}_t^k(\tilde{\omega}) - \tilde{\Lambda}_t^\infty(\tilde{\omega})\right)[\psi]\right| + \left|\left(\tilde{\Lambda}_t^k(\tilde{\omega}) - \tilde{\Lambda}_t^\infty(\tilde{\omega})\right)[(1 - \phi_K)f]\right|.$$

Since $\tilde{\Lambda}^k(\tilde{\omega})$ is $\mathbf{M}_1^+(\mathbb{R}^d)$-valued for all $k \in \mathbb{N} \cup \{\infty\}$, we get $\Lambda_t^k(\tilde{\omega})[1-\phi_K] = 1 - \Lambda_t^k(\tilde{\omega})[\phi_K]$. Therefore $|(\Lambda_t^k(\tilde{\omega}) - \Lambda_t^\infty(\tilde{\omega}))[1-\phi_K]| = |(\Lambda_t^k(\tilde{\omega}) - \Lambda_t^\infty(\tilde{\omega}))[\phi_K]|$. Hence (3.20) and (3.21) imply that there exists $k_0 \in \mathbb{N}$, such that $\sup_{t \in [0,T]} |\tilde{\Lambda}_t^k(\tilde{\omega})[1 - \phi_K]| < \varepsilon$ for all $k \geq k_0(\tilde{\omega})$ so that

$$\left|\left(\tilde{\Lambda}_t^k(\tilde{\omega}) - \tilde{\Lambda}_t^\infty(\tilde{\omega})\right)[f]\right| \leq 2\varepsilon + \left|\left(\tilde{\Lambda}_t^k(\tilde{\omega}) - \tilde{\Lambda}_t^\infty(\tilde{\omega})\right)[\psi]\right| + 2\varepsilon\|f\|_\infty \qquad \text{whenever } k \geq k_0(\tilde{\omega}).$$

Since $\varepsilon > 0$ is arbitrary, (3.20) gives $\lim_{k \to \infty} \sup_{t \in [0,T]} |(\tilde{\Lambda}_t^k(\tilde{\omega}) - \tilde{\Lambda}_t^\infty(\tilde{\omega}))[f]| = 0$, and since $f \in \mathbf{C}_b(\mathbb{R}^d)$ is arbitrary, also $\tilde{\Lambda}^k(\tilde{\omega}) \to \tilde{\Lambda}^\infty(\tilde{\omega})$ in $\mathbf{C}([0,T]; \mathcal{P}_{\mathrm{wk}^*}(\mathbb{R}^d))$ as $k \to \infty$. This yields the second claim.

Setting
$$\bar{\Omega}_{r,u}^* := \mathbf{C}\Big([0,T]^\Delta; \mathcal{P}_{\mathrm{wk}^*}(\mathbb{R}^d) \times \mathsf{H}_r^u(\mathbb{R}^d) \times \mathbb{R}^m\Big),$$

we get from the first claim and (3.18) that the measures $\bar{\mathbf{P}}_{x_0}^{n_k}$ for $k \in \mathbb{N} \cup \{\infty\}$ concentrate all their mass on $\bar{\Omega}_{r,u}^*$, and from the second claim and (3.19) that $\bar{\mathbf{P}}_{x_0}^{n_k} \to \bar{\mathbf{P}}_{x_0}^\infty$ narrowly in $\mathbf{M}_1^+(\bar{\Omega}_{r,u}^*)$ as $k \to \infty$.

**Step 6** We now pass to the limit in (3.15). As in Step 5, let $(n_k)_{k \in \mathbb{N}}$ be a subsequence such that $\bar{\mathbf{P}}_{x_0}^{n_k} \to \bar{\mathbf{P}}_{x_0}^\infty$ narrowly in $\mathbf{M}_1^+(\bar{\Omega})$ as $k \to \infty$. By Step 5, we can assume that the measures $\bar{\mathbf{P}}_{x_0}^{n_k}$ for $k \in \mathbb{N} \cup \{\infty\}$ concentrate all their mass on the Polish space $\bar{\Omega}_{r,u}^*$. Note that the set $D_{\delta,\phi}$ from Step 3 is open in $\bar{\Omega}^*_{r,u}$. Indeed, for each $\phi \in \mathcal{S}$ and $(\varepsilon, t) \in [0, T]^\Delta$, the restriction of the map $(\bar{\Lambda}, A, \bar{Z}) \mapsto A_{\varepsilon,t}[\phi]$ to $\bar{\Omega}^*_{r,u}$ is continuous, and also (3.16) is continuous on $\bar{\Omega}^*_{r,u}$ by Step 4. Note also that the tower property gives

$$\bar{\mathbf{P}}^{n_k}_{x_0}[D_{\delta,\phi}] = \mathbf{E}^{\bar{\mathbf{P}}^{n_k}_{x_0}}\left[\bar{\mathbf{P}}^{n_k}_{x_0}[D_{\delta,\phi}|\Lambda_{t-\varepsilon}, Z_{\cdot\wedge t}]\right] \qquad \text{for all } k \in \mathbb{N}.$$

Therefore the Portmanteau theorem and the concentration inequality (3.15) show that

$$\begin{aligned}\bar{\mathbf{P}}^\infty_{x_0}[D_{\delta,\phi}] &= \bar{\mathbf{P}}^\infty_{x_0}[D_{\delta,\phi} \cap \bar{\Omega}^*_{r,u}] \\ &\leq \liminf_{k\to\infty} \bar{\mathbf{P}}^{n_k}_{x_0}[D_{\delta,\phi} \cap \bar{\Omega}^*_{r,u}] \\ &= \liminf_{k\to\infty} \bar{\mathbf{P}}^{n_k}_{x_0}[D_{\delta,\phi}] \\ &\leq \liminf_{k\to\infty} \exp(-n_k\delta^2/\|\phi\|_\infty) = 0.\end{aligned}$$

As a consequence, we have $\bar{\mathbf{P}}^\infty_{x_0}[D_{\delta,\phi}] = 0$ for any $\delta > 0$. In other words, by (3.14),

$$A_{\varepsilon,t}[\phi] = \int_{\mathbb{R}^d} F^\phi_{\varepsilon,t}(x; \Lambda_{t-\varepsilon}, Z_{\cdot\wedge t})\, \Lambda_{t-\varepsilon}(\mathrm{d}x) \qquad \bar{\mathbf{P}}^\infty_{x_0}\text{-a.s.} \tag{3.23}$$

for each fixed $\phi \in \mathcal{S}$. Because $\mathcal{S}$ is separable, (3.23) extends to all of $\mathcal{S}$ up to a global $\bar{\mathbf{P}}^\infty_{x_0}$-nullset which may, however, depend on our initially chosen $\varepsilon$ and $t$.

***Step 7*** It remains to derive the probabilistic representation given in Proposition 3.1. Consider the function $\mathbb{R}^d \times \mathcal{S}' \times \mathbf{C}([0, T]; \mathbb{R}^m) \to \mathbf{M}^+_1(\mathbb{R}^d)$ given by mapping $(x, \lambda, z)$ to the normal law on $\mathbb{R}^d$ with mean vector $M_{\varepsilon,t}(x, \lambda, z)$ and covariance matrix $\Sigma_{\varepsilon,t}(x, \lambda)$ defined in (3.1) and (3.2), respectively. With the auxiliary probability space $([0,1], \mathcal{B}([0,1]), \mathrm{d}y)$, we can find a random variable $U : [0, 1] \times \mathbb{R}^d \times \mathcal{S}' \times \mathbf{C}([0, T]; \mathbb{R}^m) \to \mathbb{R}^d$ such that $\bar{\mathbf{P}}^\infty_{x_0}$-a.s.,

$$\mathrm{Law}_{([0,1], \mathcal{B}([0,1]), \mathrm{d}y)}\left(U(\,\cdot\,, x, \Lambda_{t-\varepsilon}, Z_{\cdot\wedge t})\right) = \mathcal{N}\big(M_{\varepsilon,t}(x, \Lambda_{t-\varepsilon}, Z_{\cdot\wedge t}), \Sigma_{\varepsilon,t}(x, \Lambda_{t-\varepsilon})\big),$$

simultaneously for all $x \in \mathbb{R}^d$; see for instance Kallenberg [31, Lem. 4.22]. Suppressing the $y$-argument and introducing the notation $U_{\lambda,z}(x) := U(x, \lambda, z)$, we then obtain from (3.23) and the definition (3.8) of $F^\phi_{\varepsilon,t}$ that $\bar{\mathbf{P}}^\infty_{x_0}$-a.s.,

$$A_{\varepsilon,t}[\phi] = \int_{\mathbb{R}^d} \mathbb{E}\big[\phi\big(U_{\Lambda_{t-\varepsilon}, Z_{\cdot\wedge t}}(x)\big)\big]\, \Lambda_{t-\varepsilon}(\mathrm{d}x) \quad \text{for all } \phi \in \mathcal{S}.$$

This is the claimed representation (3.3). $\square$

## 3.2 Estimates for the regular part

Via the probabilistic representation in Proposition 3.1, we derive a bound for $A$.

**Proposition 3.2** | *Impose Assumptions 1.1 and 1.2 and fix real numbers $1 < r < \infty$ and $s \geq 0$. Let $r'$ be conjugate to $r$. There exists a constant $c_{r,s} < \infty$ such that for any fixed pair $(\varepsilon, t) \in [0, T]^\Delta$ satisfying $0 < \varepsilon \leq 1 \wedge t$ and any narrow cluster point $\bar{\mathbf{P}}^\infty_{x_0}$ of $(\bar{\mathbf{P}}^n_{x_0})_{n \in \mathbb{N}}$, we have*

$$\|A_{\varepsilon,t}\|_{\mathsf{H}^s_r(\mathbb{R}^d)} \leq c_{r,s}\varepsilon^{-(d/r'+s)/2} \quad \bar{\mathbf{P}}^\infty_{x_0}\text{-a.s.} \tag{3.24}$$

*The constant $c_{r,s}$ depends on the initially fixed choices of $r$, $s$ and $d$.*

The proof of Proposition 3.2 is given in Section 3.2.2 after collecting some auxiliary results in Section 3.2.1. But first, let us make a comment.

Note that since Proposition 3.2 is valid for all $s \geq 0$, it follows that $\bar{\mathbf{P}}_{x_0}^\infty$-a.s., we have for each $(\varepsilon, t)$ with $0 < \varepsilon \leq t \leq T$ that $A_{\varepsilon,t} \in \bigcap_{s \in \mathbb{N}} \mathsf{H}_r^s(\mathbb{R}^d)$. This implies that the distribution $A_{\varepsilon,t} \in \mathcal{S}'(\mathbb{R}^d)$ can be identified with a smooth function; see e.g. Bergh and Löfström [6, Thm. 6.3.2]. It is this observation that leads us to call $A$ the regular part of the decomposition.

### 3.2.1 Preliminaries for the proof of Proposition 3.2

We need two technical lemmata for reference. As to not distract from our main endeavor, we defer their proofs to Appendix B.

**Lemma 3.3** | *Let $1 < r < \infty$, $s \geq 0$ and $f \in \mathsf{H}_r^s(\mathbb{R}^d)$.*

*1) Define by $(\tau_m f)(x) := f(x + m)$ for $m \in \mathbb{R}^d$ the translation of $f$. Then we have for all $m \in \mathbb{R}^d$ that $\|\tau_m f\|_{\mathsf{H}_r^s(\mathbb{R}^d)} = \|f\|_{\mathsf{H}_r^s(\mathbb{R}^d)}$, i.e., the $\mathsf{H}_r^s$-norm is translation-invariant.*

*2) Define by $f_\delta(x) := f(\delta x)$ for $\delta \geq 1$ the dilation of $f$. Then we have for all $\delta \geq 1$ that $\|f_\delta\|_{\mathsf{H}_r^s(\mathbb{R}^d)} \leq \delta^{s-d/r} \|f\|_{\mathsf{H}_r^s(\mathbb{R}^d)}$.*

For the following result, we lack a good reference and therefore include it.

**Lemma 3.4** | *Let $\mathcal{V} \subseteq \mathbb{R}^{d \times d}$ be a set of symmetric matrices that is bounded and satisfies the uniform ellipticity condition that*

$$z^\top V z \geq c |z|^2 \quad \text{for all } z \in \mathbb{R} \text{ and } V \in \mathcal{V}, \text{ for some } c > 0.$$

*Then for any $r > 1$ and $s \geq 0$, we have for the normal density $g(\,\cdot\,; m, v)$ that*

$$\sup_{V \in \mathcal{V}} \|g(\,\cdot\,; 0, V)\|_{\mathsf{H}_r^s(\mathbb{R}^d)} \leq c_{r,s}, \tag{3.25}$$

*where $c_{r,s} < \infty$ is a constant depending on $c$, and on $r$, $s$ and $d$.*

### 3.2.2 Proof of Proposition 3.2

To prove Proposition 3.2 and make use of Lemma 3.4 we rewrite the probabilistic representation (3.3) of $A_{\varepsilon,t}[\phi]$ in analytic terms. Introduce the short-hand

$$G_{\varepsilon,t}(y; x, \Lambda_{t-\varepsilon}, Z_{\cdot \wedge t}) := g\big(y; M_{\varepsilon,t}(x, \Lambda_{t-\varepsilon}, Z_{\cdot \wedge t}), \Sigma_{\varepsilon,t}(x, \Lambda_{t-\varepsilon})\big) \tag{3.26}$$

for the density of $\mathcal{N}(M_{\varepsilon,t}(x, \Lambda_{t-\varepsilon}, Z_{\cdot \wedge t}), \Sigma_{\varepsilon,t}(x, \Lambda_{t-\varepsilon}))$. We then obtain for the expectation in (3.3) the integral expression

$$A_{\varepsilon,t}[\phi] = \int_{\mathbb{R}^d} \left( \int_{\mathbb{R}^d} \phi(y) G_{\varepsilon,t}(y; x, \Lambda_{t-\varepsilon}, Z_{\cdot \wedge t}) \, \mathrm{d}y \right) \Lambda_{t-\varepsilon}(\mathrm{d}x). \tag{3.27}$$

The representation (3.27) shows that $A_{\varepsilon,t}$ can be viewed as a measurable function of $\Lambda_{t-\varepsilon}$ and $Z_{\cdot \wedge t}$. We can now finally give the proof of Proposition 3.2.

***Proof of Proposition 3.2*** Let $r > 1$ and $s \geq 0$ and $(\varepsilon, t) \in [0, T]^\Delta$ be fixed. Consider $\phi \in \mathcal{S}$. Then (3.27) and the duality between $\mathsf{H}_{r'}^{-s}(\mathbb{R}^d)$ and $\mathsf{H}_r^s(\mathbb{R}^d)$ from Lemma 2.5 give

$$\begin{aligned} |A_{\varepsilon,t}[\phi]| &= \left| \int_{\mathbb{R}^d} \left( \int_{\mathbb{R}^d} \phi(y) G_{\varepsilon,t}(y; x, \Lambda_{t-\varepsilon}, Z_{\cdot \wedge t}) \, \mathrm{d}y \right) \Lambda_{t-\varepsilon}(\mathrm{d}x) \right| \\ &\leq \int_{\mathbb{R}^d} \|\phi\|_{\mathsf{H}_{r'}^{-s}(\mathbb{R}^d)} \|G_{\varepsilon,t}(\,\cdot\,; x, \Lambda_{t-\varepsilon}, Z_{\cdot \wedge t})\|_{\mathsf{H}_r^s(\mathbb{R}^d)} \Lambda_{t-\varepsilon}(\mathrm{d}x). \end{aligned} \tag{3.28}$$

To bound the $\mathsf{H}_r^s$-norm of $G_{\varepsilon,t}(\,\cdot\,;x,\Lambda_{t-\varepsilon},Z_{\cdot\wedge t})$, fix a $\bar{\mathbf{P}}_{x_0}^\infty$-nullset $N \in \bar{\mathcal{F}}$ such that $\Lambda_t(\bar{\omega})$ is in $\mathbf{M}_1^+(\mathbb{R}^d)$ for all $\bar{\omega} \in \bar{\Omega}\setminus N$. Using the translation-invariance of the $\mathsf{H}_r^s$-norm in part 1) of Lemma 3.3 and recalling the short-hand (3.26), we find that on $\bar{\Omega}\setminus N$

$$\|G_{\varepsilon,t}(\,\cdot\,;x,\Lambda_{t-\varepsilon},Z_{\cdot\wedge t})\|_{\mathsf{H}_r^s(\mathbb{R}^d)} = \left\|g\big(\,\cdot\,;M_{\varepsilon,t}(x,\Lambda_{t-\varepsilon},Z_{\cdot\wedge t}),\Sigma_{\varepsilon,t}(x,\Lambda_{t-\varepsilon})\big)\right\|_{\mathsf{H}_r^s(\mathbb{R}^d)}$$
$$= \left\|g\big(\,\cdot\,;0,\Sigma_{\varepsilon,t}(x,\Lambda_{t-\varepsilon})\big)\right\|_{\mathsf{H}_r^s(\mathbb{R}^d)}. \qquad (3.29)$$

This greatly simplifies matters because we can now exploit the scaling properties of the $\mathsf{H}_r^s$-norm in part 2) of Lemma 3.3. Let

$$\tilde{\Sigma}_{\varepsilon,t}(x,\lambda) := \frac{1}{\varepsilon}\Sigma_{\varepsilon,t}(x,\lambda) = \frac{1}{\varepsilon}\int_{t-\varepsilon}^t (\sigma_r\sigma_r^\top)(x,\lambda)\,\mathrm{d}r\,.$$

Due to the uniform ellipticity and boundedness of $\sigma$ in Assumption 1.2(E), the computation

$$z^\top \tilde{\Sigma}_{\varepsilon,t}(x,\lambda)z = \frac{1}{\varepsilon}\int_{t-\varepsilon}^t z^\top(\sigma_r\sigma_r^\top)(x,\lambda)z\,\mathrm{d}r \geq \frac{1}{\varepsilon}\int_{t-\varepsilon}^t \kappa(z^\top z)\,\mathrm{d}s = \kappa|z|^2$$

shows that the family $\mathcal{V} := \{\tilde{\Sigma}_{\varepsilon,t}(x,\lambda) : 0 < \varepsilon \leq t, x \in \mathbb{R}^d, \lambda \in \mathbf{M}_1^+(\mathbb{R}^d)\}$ is uniformly elliptic and bounded. Now the explicit form of the Gaussian density gives

$$g\big(y;0,\Sigma_{\varepsilon,t}(x,\Lambda_{t-\varepsilon})\big) = g\big(y;0,\varepsilon\tilde{\Sigma}_{\varepsilon,t}(x,\Lambda_{t-\varepsilon})\big) = \varepsilon^{-d/2}g\big(\varepsilon^{-1/2}y;0,\tilde{\Sigma}_{\varepsilon,t}(x,\Lambda_{t-\varepsilon})\big)\,,$$

and combining this with the scaling properties of the $\mathsf{H}_r^s$-norm in part 2) of Lemma 3.3 and the auxiliary bound from Lemma 3.4 yields on $\bar{\Omega}\setminus N$, the bound

$$\left\|g\big(\,\cdot\,;0,\Sigma_{\varepsilon,t}(x,\Lambda_{t-\varepsilon})\big)\right\|_{\mathsf{H}_r^s(\mathbb{R}^d)} = \varepsilon^{-d/2}\left\|g\big(\varepsilon^{-1/2}\,\cdot\,;0,\tilde{\Sigma}_{\varepsilon,t}(x,\Lambda_{t-\varepsilon})\big)\right\|_{\mathsf{H}_r^s(\mathbb{R}^d)}$$
$$= \varepsilon^{-(s-d/r)/2-d/2}\left\|g\big(\,\cdot\,;0,\tilde{\Sigma}_{\varepsilon,t}(x,\Lambda_{t-\varepsilon})\big)\right\|_{\mathsf{H}_r^s(\mathbb{R}^d)}$$
$$\leq c(r,s)\varepsilon^{-(d/r'+s)/2}\,. \qquad (3.30)$$

Returning to (3.28) and noting that $\Lambda_t(\bar{\omega}) \in \mathbf{M}_1^+(\mathbb{R}^d)$ on $\bar{\Omega}\setminus N$, we see via (3.29) and (3.30) that

$$|A_{\varepsilon,t}[\phi]| \leq c(r,s)\|\phi\|_{\mathsf{H}_{r'}^{-s}(\mathbb{R}^d)}\varepsilon^{-(d/r'+s)/2} \qquad \bar{\mathbf{P}}_{x_0}^\infty\text{-a.s.} \qquad (3.31)$$

Now the dual representation of the $\mathsf{H}_r^s$-norm from Lemma 2.5 shows that $\bar{\mathbf{P}}_{x_0}^\infty$-a.s.,

$$\|A_{\varepsilon,t}\|_{\mathsf{H}_r^s(\mathbb{R}^d)} \leq c(r,s)\sup\{|A_{\varepsilon,t}[\phi]| : \phi \in \mathcal{S}, \|\phi\|_{\mathsf{H}_{r'}^{-s}(\mathbb{R}^d)} = 1\} \leq c(r,s)\varepsilon^{-(d/r'+s)/2}\,,$$

as claimed. □

## 4 THE DISTRIBUTIONAL PART

In this section, we carry out Step 4 of the proof of Theorem 1.5 as discussed in Section 1.5. We want to obtain a quantitative norm estimate for the distributional part $E$ in the decomposition (1.16) of $\bar{\Lambda}$. This is achieved via a Kolmogorov-type argument at the level of the $n$-particle system (1.1). The bound obtained is uniform in $n$ and stable under narrow convergence so that we can pass it to the limit $n \to \infty$. In more detail, we get the following result.

**Proposition 4.1** | *Impose Assumptions 1.1 and 1.2. Let $\bar{\mathbf{P}}_{x_0}^\infty$ be a cluster point of $(\bar{\mathbf{P}}_{x_0}^n)_{n \in \mathbb{N}}$, fix reals $r \in (1, \infty)$, $q \in [1, \infty)$ and $\xi \in (0, \beta)$ with $\beta$ as in Assumption 1.2(H), take $r'$ conjugate to $r$ and set $u := 1 + d/r'$. Then for all $\alpha > 0$ and $t \in (0, T]$, there exists a random variable*

$$\bar{C}_{t,r,\alpha,\xi} \in \mathbb{L}^q(\bar{\Omega}, \bar{\mathcal{F}}, \bar{\mathbf{P}}_{x_0}^\infty) \tag{4.1}$$

*satisfying $\bar{\mathbf{P}}_{x_0}^\infty$-a.s. simultaneously for all $\varepsilon \in \{2^{-\alpha n}(1 \wedge t) : n \in \mathbb{N}_0\}$ the bound*

$$\|E_{\varepsilon,t}\|_{\mathsf{H}_r^{-u}(\mathbb{R}^d)} \leq \bar{C}_{t,r,\alpha,\xi}\, \varepsilon^{(1+\xi)/2}\,. \tag{4.2}$$

*In addition, for fixed $\alpha$, there is for each $t \in (0,T]$ a choice of $\bar{C}_{t,r,\alpha,\xi}$ such that*

$$\sup_{t \in [0,T]} \|\bar{C}_{t,r,\alpha,\xi}\|_{\mathbb{L}^q(\bar{\Omega},\bar{\mathcal{F}},\bar{\mathbf{P}}_{x_0}^\infty)} < \infty\,.$$

We begin with an auxiliary lemma establishing a uniform bound at the level of the particle systems, from which we subsequently deduce Proposition 4.1.

**Lemma 4.2** | *Impose Assumptions 1.1 and 1.2. Then for any $q \in [1,\infty)$ and with $\beta$ as in Assumption 1.2(H), there exists a constant $c = c(q, \sigma, \bar{\sigma}, b, T, \beta) < \infty$ such that for all $t \in (0,T]$, we have*

$$\|X_t^{i,n} - Y_t^{i,n;\varepsilon,t}\|_{\mathbb{L}^q(\mathbb{P}_{x_0}^n)} \leq c\varepsilon^{(1+\beta)/2} \tag{4.3}$$

*for all $0 \leq \varepsilon \leq t$, $n \in \mathbb{N}$ and $i \in [n]$.*

***Proof*** This argument is very similar to the proof of Proposition 2.2. We first use (1.17) and (1.18) for each $i \in [n]$ to write

$$X_t^{i,n} - Y_t^{i,n;\varepsilon,t} = J_0^{i,n} + J_1^{i,n} + J_2^{i,n} \tag{4.4}$$

with

$$J_0^{i,n} := \int_{t-\varepsilon}^{t} b_r(X_r^{i,n}, \mu_r^n)\,\mathrm{d}r\,,$$
$$J_1^{i,n} := \int_{t-\varepsilon}^{t} \Big(\sigma_r(X_r^{i,n}, \mu_r^n) - \sigma_r(X_{t-\varepsilon}^{i,n}, \mu_{t-\varepsilon}^n)\Big)\mathrm{d}B_r^{i,n}\,,$$
$$J_2^{i,n} := \int_{t-\varepsilon}^{t} \Big(\bar{\sigma}_r(X_r^{i,n}, \mu_r^n) - \bar{\sigma}_{t-\varepsilon}(X_{t-\varepsilon}^{i,n}, \mu_{t-\varepsilon}^n)\Big)\mathrm{d}Z_r^n\,.$$

Since $b$ is bounded by Assumption 1.1, we control $J_0^{i,n}$ by

$$\|J_0^{i,n}\|_{\mathbb{L}^q(\mathbb{P}_{x_0}^n)} \leq \|b\|_\infty \varepsilon \leq \|b\|_\infty T^{(1-\beta)/2} \varepsilon^{(1+\beta)/2}\,, \tag{4.5}$$

uniformly for $i \in [n]$ and $n \in \mathbb{N}$. The terms $J_1^{i,n}$ and $J_2^{i,n}$ are almost like $G_5^{i,n}$ and $G_6^{i,n}$ in (2.14) and (2.15), and the sum $J_1^{i,n} + J_2^{i,n}$ is of the same form as $G_1^{i,n} - G_4^{i,n}$ in Step 2.2 of the proof of Proposition 2.2. In the case here, the integrals are taken over $(t - \varepsilon, t]$ instead of $(t' - \varepsilon', t]$ there, with the arguments of $\sigma_r$ and $\bar{\sigma}_r$ in the first integrals being $(X_r^{i,n}, \mu_r^n)$ here instead of $(X_{t'-\varepsilon'}^{i,n}, \mu_{t'-\varepsilon'}^n)$ there. As pointed out in Remark 2.3, with $\varepsilon'$ there replaced by $\varepsilon$ here, we get a good estimate in this case. Using directly the bound (2.18) with the corresponding replacements and noting that we also have $|r - (t - \varepsilon)| \leq \varepsilon$ when $r \in [t - \varepsilon, t]$, we obtain

$$\|J_1^{i,n} + J_2^{i,n}\|_{\mathbb{L}^q(\mathbb{P}_{x_0}^n)} \leq c(q, \sigma, \bar{\sigma}, b, T, \beta) \varepsilon^{(1+\beta)/2} .$$

Combining this with (4.5) and (4.4) gives (4.3). □

We now come to the proof of Proposition 4.1. While the uniform bound of Lemma 4.2 is essentially the same as that in Proposition 2.2, which we used as an ingredient for the Kolmogorov continuity criterion, the estimate (4.3) is not strong enough to deduce (4.2) in Proposition 4.1 by the same route. Indeed, for such a strategy to be successful, we should need for $\varepsilon \neq \varepsilon'$ a bound on $\|(X_t^{i,n} - Y_{i,n}^{i,n;\varepsilon,t}) - (X_t^{i,n} - Y_{i,n}^{i,n;\varepsilon',t})\|_{\mathbb{L}^q(\mathbb{P}_{x_0}^n)}$ of the order $|\varepsilon - \varepsilon'|^{(1+\beta)/2}$. Such a bound is, however, not true, as can be seen by considering the case of nonconstant $\sigma, \bar{\sigma}$. In fact, from this case, we see that we can only hope for a bound of order $|\varepsilon - \varepsilon'|^{1/2}$.

On the other hand, for $0 = \varepsilon' < \varepsilon \leq t$, we get the strong bound (4.3) in Lemma 4.2 which is enough to mimic the ideas used in classical proofs of the Kolmogorov criterion; see Cohen and Elliott [14, Thm. A.4.1]. In this way, we can show $\bar{\mathbf{P}}_{x_0}^\infty$-almost sure *local Hölder-continuity at zero of order* $(1 + \xi)/2$ *on a countable set* for any $\xi < \beta$, which is exactly the bound in (4.2) in Proposition 4.1 that we now prove.

***Proof of Proposition 4.1*** We start with the weak solution $\boldsymbol{X}^n$ of the $n$-particle system from Lemma 1.3 and the auxiliary process $(Y_s^{n;\varepsilon,t})_{s \in [t-\varepsilon,t],(\varepsilon,t) \in [0,T]^\Delta}$ defined in (1.18). Observe that if $r'$ is conjugate to $r$ and $u = 1 + d/r'$, then by Lemma 2.4, any $\phi \in \mathsf{H}_{r'}^u(\mathbb{R}^d)$ is $\rho$-Hölder-continuous for any $0 < \rho < 1$ and we have the bound $[\phi]_\rho \leq c(r, \rho) \|\phi\|_{\mathsf{H}_{r'}^u(\mathbb{R}^d)}$. Therefore, using the definitions of $\mu^n$ and $\nu^n$ from (1.19) and (1.20), we see that

$$|(\mu_t^n - \nu_{\varepsilon,t}^n)[\phi]| = \left| \frac{1}{n} \sum_{i=1}^n \left( \phi(X_t^{i,n}) - \phi(Y_t^{i,n;\varepsilon,t}) \right) \right| \leq [\phi]_\rho \frac{1}{n} \sum_{i=1}^n |X_t^{i,n} - Y_t^{i,n;\varepsilon,t}|^\rho . \quad (4.6)$$

Together with the dual characterization of the $\mathsf{H}_r^{-u}$-norm from Lemma 2.5, we get

$$\|\mu_t^n - \nu_{\varepsilon,t}^n\|_{\mathsf{H}_r^{-u}(\mathbb{R}^d)} \leq c(r) \sup\{|(\mu_t^n - \nu_{\varepsilon,t}^n)[\phi]| : \|\phi\|_{\mathsf{H}_{r'}^u(\mathbb{R}^d)} = 1\}$$

$$\leq c(r,\rho) \frac{1}{n} \sum_{i=1}^n |X_t^{i,n} - Y_t^{i,n;\varepsilon,t}|^\rho$$

for $0 < \rho < 1$. The fact that $\rho < 1$ implies $\|U^\rho\|_{\mathbb{L}^q(\mathbb{P}_{x_0}^n)} \leq \|U\|_{\mathbb{L}^q(\mathbb{P}_{x_0}^n)}^\rho$, and the triangle inequality now give

$$\left\| \|\mu_t^n - \nu_{\varepsilon,t}^n\|_{\mathsf{H}_r^{-u}(\mathbb{R}^d)} \right\|_{\mathbb{L}^q(\mathbb{P}_{x_0}^n)} \leq c(r,\rho) \frac{1}{n} \sum_{i=1}^n \|X_t^{i,n} - Y_t^{i,n;\varepsilon,t}\|_{\mathbb{L}^q(\mathbb{P}_{x_0}^n)}^\rho . \quad (4.7)$$

On the right-hand side of this inequality, we use the uniform estimate (4.3) from Lemma 4.2. Taking into account the constant $c(r, \rho)$ from (4.7) yields for any $q \geq 1$ a new constant $c = c(q, \sigma, \bar{\sigma}, b, T, \beta, r, \rho) < \infty$ such that

$$\sup_{n \in \mathbb{N}} \left\| \|\mu_t^n - \nu_{\varepsilon,t}^n\|_{\mathsf{H}_r^{-u}(\mathbb{R}^d)} \right\|_{\mathbb{L}^q(\mathbb{P}_{x_0}^n)} \leq c\varepsilon^{\rho(1+\beta)/2} \tag{4.8}$$

for all $(\varepsilon, t) \in [0, T]^\Delta$. Note that the constant $c$ in (4.8) does not depend on $t \in [0, T]$. Since (2.33) implies that $\mathrm{Law}_{\bar{\mathbb{P}}_{x_0}^n}(\Lambda_t, A_{\varepsilon,t}) = \mathrm{Law}_{\mathbb{P}_{x_0}^n}(\mu_t^n, \nu_{\varepsilon,t}^n)$ for any $n \in \mathbb{N}$ and since $E_{\varepsilon,t} = \Lambda_t - A_{\varepsilon,t}$ by the definition in (1.16), the bound (4.8) is equivalent to

$$\sup_{n \in \mathbb{N}} \left\| \|E_{\varepsilon,t}\|_{\mathsf{H}_r^{-u}(\mathbb{R}^d)} \right\|_{\mathbb{L}^q(\bar{\mathbb{P}}_{x_0}^n)} \leq c\varepsilon^{\rho(1+\beta)/2}$$

for all $(\varepsilon, t) \in [0, T]^\Delta$. Thus if $\bar{\mathbb{P}}_{x_0}^{n_k} \to \bar{\mathbb{P}}_{x_0}^\infty$ narrowly as $k \to \infty$, then

$$\left\| \|E_{\varepsilon,t}\|_{\mathsf{H}_r^{-u}(\mathbb{R}^d)} \right\|_{\mathbb{L}^q(\bar{\mathbb{P}}_{x_0}^\infty)} \leq \liminf_{k \to \infty} \left\| \|E_{\varepsilon,t}\|_{\mathsf{H}_r^{-u}(\mathbb{R}^d)} \right\|_{\mathbb{L}^q(\bar{\mathbb{P}}_{x_0}^{n_k})} \leq c\varepsilon^{\rho(1+\beta)/2}. \tag{4.9}$$

Let us now fix $0 < \xi < \beta$ and $\alpha > 0$. With these choices, we can choose $0 < \rho < 1$ sufficiently large such that $\rho(1+\beta) > 1 + \xi$. We next introduce the shorthand notation

$$(0, 1 \wedge t]_\alpha := \{2^{-\alpha n}(1 \wedge t) : n \in \mathbb{N}_0\}$$

and define the real-valued random variable

$$\bar{C}_{t,r,\alpha,\xi} := \sup \left\{ \frac{\|E_{\varepsilon,t}\|_{\mathsf{H}_r^{-u}(\mathbb{R}^d)}}{\varepsilon^{(1+\xi)/2}} : \varepsilon \in (0, 1 \wedge t]_\alpha \right\}. \tag{4.10}$$

Because $\sup_{n \in \mathbb{N}} |U_n| \leq \sum_{n \in \mathbb{N}} |U_n|$, (4.9) and $\rho(1+\beta) > 1 + \xi$ imply

$$\|\bar{C}_{t,r,\alpha,\xi_0}\|_{\mathbb{L}^q(\bar{\mathbb{P}}_{x_0}^\infty)} \leq \left\| \sum_{\varepsilon \in (0,1 \wedge t]_\alpha} \frac{\|E_{\varepsilon,t}\|_{\mathsf{H}_r^{-u}(\mathbb{R}^d)}}{\varepsilon^{(1+\xi)/2}} \right\|_{\mathbb{L}^q(\bar{\mathbb{P}}_{x_0}^\infty)}$$

$$\leq \sum_{\varepsilon \in (0,1 \wedge t]_\alpha} \left\| \frac{\|E_{\varepsilon,t}\|_{\mathsf{H}_r^{-u}(\mathbb{R}^d)}}{\varepsilon^{(1+\xi)/2}} \right\|_{\mathbb{L}^q(\bar{\mathbb{P}}_{x_0}^\infty)}$$

$$\leq \sum_{\varepsilon \in (0,1 \wedge t]_\alpha} c\varepsilon^{\rho(1+\beta)/2 - (1+\xi)/2}$$

$$= \sum_{n \in \mathbb{N}_0} c\left(2^{-\alpha n}(1 \wedge t)\right)^{\rho(1+\beta)/2 - (1+\xi)/2}$$

$$=: c_{q,\rho,\beta,\xi,\alpha,q,r}(1 \wedge t)^{\rho(1+\beta)/2 - (1+\xi)/2} < \infty, \tag{4.11}$$

where the constant $c_{q,\rho,\beta,\xi,\alpha,q,r}$ depends only on $\alpha$, $q$, $\rho$, $r$ and the gap $\rho(1+\beta) - (1+\xi)$. The definition (4.10) and the bound (4.11) show that $\bar{\mathbb{P}}_{x_0}^\infty$-a.s., the (random) map

$$(0, 1 \wedge t]_\alpha \ni \varepsilon \mapsto \|E_{\varepsilon,t}\|_{\mathsf{H}_r^{-u}(\mathbb{R}^d)} \in \mathbb{R}$$

is locally Hölder-continuous at zero of order $(1 + \xi)/2$ with (random) Hölder-norm $\bar{C}_{t,r,\alpha,\xi} \in \mathbb{L}^q(\bar{\Omega}, \bar{\mathcal{F}}, \bar{\mathbb{P}}_{x_0}^\infty)$. This proves (4.1) and (4.2). In addition, (4.11) shows that

$$\sup_{t \in (0,T]} \|\bar{C}_{t,r,\alpha,\xi}\|_{\mathbb{L}^q(\bar{\mathbb{P}}_{x_0}^\infty)} \leq c_{q,\rho,\beta,\xi,\alpha,q,r}(1 \wedge T)^{\rho(1+\beta)/2 - (1+\xi)/2} < \infty,$$

and the proof is complete. $\square$

# 5 AN ABSTRACT INTERPOLATION RESULT

Recall from Appendix A the Bessel potential spaces $\mathsf{H}_r^s(\mathbb{R}^d)$ for $r \in (1, \infty)$ and $s \in \mathbb{R}$. For $u, s > 0$ and $w \in (-u, s)$, we have the chain of continuous inclusions

$$\mathcal{S}(\mathbb{R}^d) \hookrightarrow \mathsf{H}_r^s(\mathbb{R}^d) \hookrightarrow \mathsf{H}_r^w(\mathbb{R}^d) \hookrightarrow \mathsf{H}_r^{-u}(\mathbb{R}^d) \hookrightarrow \mathcal{S}'(\mathbb{R}^d);$$

see (5.7) below. In this setting, $\mathsf{H}_r^w(\mathbb{R}^d)$ is called an *intermediate space* between $\mathsf{H}_r^s(\mathbb{R}^d)$ and $\mathsf{H}_r^{-u}(\mathbb{R}^d)$. Recall also that $\mathsf{H}_r^0(\mathbb{R}^d) = \mathsf{L}^r(\mathbb{R}^d, \mathrm{d}x)$; so for $w \geq 0$, elements of $\mathsf{H}_r^w(\mathbb{R}^d)$ are functions, not merely distributions.

The following abstract interpolation result, on which we heavily rely, lets us obtain bounds for an intermediate $\mathsf{H}_r^w(\mathbb{R}^d)$-norm if we have an appropriate linear decomposition with quantitative information on the norms of the summands in the respective end-point spaces $\mathsf{H}_r^{-u}(\mathbb{R}^d)$ and $\mathsf{H}_r^s(\mathbb{R}^d)$. This is achieved by using classical interpolation theory.

**Proposition 5.1** | *Consider reals $s > 0$, $1 < r < 2$ and $\xi > 0$, take $r'$ conjugate to $r$ and set $u := 1 + d/r'$. Let $\lambda \in \mathcal{S}'(\mathbb{R}^d)$ be a distribution and assume that there exists $\varepsilon_0 > 0$ such that for each $\varepsilon \in (0, \varepsilon_0]$, we have a decomposition*

$$\lambda = a_\varepsilon + e_\varepsilon \tag{5.1}$$

*for some $a_\varepsilon$ and $e_\varepsilon$ in $\mathcal{S}'(\mathbb{R}^d)$ satisfying*

$$\|a_\varepsilon\|_{\mathsf{H}_r^s(\mathbb{R}^d)} \leq c_a \varepsilon^{-(d/r'+s)/2} \quad \text{and} \quad \|e_\varepsilon\|_{\mathsf{H}_r^{-u}(\mathbb{R}^d)} \leq c_e \varepsilon^{(1+\xi)/2} \tag{5.2}$$

*with constants $c_a$, $c_e < \infty$ that may depend on $\varepsilon_0$, but are independent of $\varepsilon \neq \varepsilon_0$. Then the following hold:*

*1) There exists $w_0 = w_0(d, \xi, s, r) \in (-u, s)$ such that $\lambda \in \mathsf{H}_r^w(\mathbb{R}^d)$ for all $w < w_0$. More specifically, we have*

$$\|\lambda\|_{\mathsf{H}_r^w(\mathbb{R}^d)} \leq c(c_a \varepsilon_0^{-(d/r'+s)/2} + c_e \varepsilon_0^{(1+\xi)/2}) \tag{5.3}$$

*for a constant $c = c(w, d, \xi, s, r) < \infty$ independent of $\varepsilon_0$.*

*2) For all $\xi > 0$, there exist combinations of parameters $s > 0$ and $1 < r < 2$ satisfying*

$$\gamma := (d/r' + s)/2 < 1$$

*and such that if the decomposition (5.1) and the bounds in (5.2) hold, then we have $w_0(d, \xi, s, r) \in (0, s)$ in part 1). In particular, $\lambda$ is then a function in some $\mathsf{H}_r^w(\mathbb{R}^d)$ with $w > 0$.*

*3) For all $\xi > 0$, $s > 0$ and $1 < r < 2$, there exists a number $\alpha_0 = \alpha_0(\xi, r, s) > 0$ such that if the decomposition in (5.1) and the norm bounds in (5.2) are valid only for*

$$\varepsilon \in \{2^{-\alpha_0 n} \varepsilon_0 : n \in \mathbb{N}_0\} =: (0, \varepsilon_0]_{\alpha_0} \subseteq (0, \varepsilon_0],$$

*then the conclusions from parts 1) and 2) already hold.*

The proof of Proposition 5.1 is given in Section 5.2, after recalling some background and notation from interpolation theory in Section 5.1 .

Proposition 5.1 falls into a larger family of regularity results for probability laws which grew out of the seminal contribution of Fournier and Printems [25]. Results by Debussche and Romito [21] and Bally and Caramellino [5] are closely related to ours; see Section 7 for a fuller discussion, a comparison and further references. Proposition 5.1, however, is more tailored to our needs. In Section 6 below, we apply it to the (random) decomposition $\Lambda_t = A_{\varepsilon,t} + E_{\varepsilon,t}$. We see from Proposition 4.1, where we obtained an $\mathsf{H}_r^{-u}(\mathbb{R}^d)$-norm bound for the (random) term $E_{\varepsilon,t}$, that part 3) of Proposition 5.1 is vital for this application to succeed. In fact, the norm bound we obtained in Proposition 4.1 holds only on sets of the form $(0, 1 \wedge t]_\alpha$ for $\alpha > 0$, but not for $(0, 1 \wedge t]$.

## 5.1 Background from interpolation theory

We start by recalling some background from classical real interpolation theory. Our exposition is specialized to the present case. We refer to Bergh and Löfström [6, Ch.3, Ch. 6] for more details on the topic. Other good accounts include Adams and Fournier [1, Ch. 7], Lunardi [38, Ch. 1] or the classic books Triebel [52, Ch. 2.4] and [51, Ch. 1].

Let $(\mathrm{X}, \|\cdot\|_\mathrm{X}) \hookrightarrow (\mathrm{Z}, \|\cdot\|_\mathrm{Z})$ be continuously embedded Banach spaces. A Banach space $(\mathrm{Y}, \|\cdot\|_\mathrm{Y})$ is called an *intermediate space between* X *and* Z if $\mathrm{X} \hookrightarrow \mathrm{Y} \hookrightarrow \mathrm{Z}$.

For $z \in \mathrm{Z}$ and $t > 0$, we define the *K-functional*

$$K(t, z) = \inf\{\|z - x\|_\mathrm{Z} + t\|x\|_\mathrm{X} : x \in \mathrm{X}\}. \tag{5.4}$$

Moreover, for $1 \leq p < \infty$ and $\theta \in (0, 1)$, we set

$$\|z\|_{\theta,p,(\mathrm{Z},\mathrm{X});K} := \left(\int_0^\infty \left(t^{-\theta} K(t, z)\right)^p \frac{dt}{t}\right)^{1/p}$$

and define the *K-interpolation space*

$$(\mathrm{Z}, \mathrm{X})_{\theta,p;K} := \left\{z \in \mathrm{Z} : \|z\|_{\theta,p,(\mathrm{Z},\mathrm{X});K} < \infty\right\}.$$

Then $((\mathrm{Z}, \mathrm{X})_{\theta,p;K}, \|\cdot\|_{\theta,p,(\mathrm{Z},\mathrm{X});K})$ is an intermediate space between X and Z; [6, Thm. 3.1.2].

There is a norm equivalent to $\|\cdot\|_{\theta,p,(\mathrm{Z},\mathrm{X});K}$ which can be described in discrete terms and turns out to be more convenient for our purposes. Indeed, defining

$$\|z\|_{\theta,p,(\mathrm{Z},\mathrm{X});k} := \left(\sum_{n=-\infty}^\infty 2^{-pn\theta} K(2^n, z)^p\right)^{1/p}, \tag{5.5}$$

one can prove the bounds

$$\underline{c}_\theta \|z\|_{\theta,p,(\mathrm{Z},\mathrm{X});k} \leq \|z\|_{\theta,p,(\mathrm{Z},\mathrm{X});K} \leq \overline{c}_\theta \|z\|_{\theta,p,(\mathrm{Z},\mathrm{X});k}, \tag{5.6}$$

where the constants $\underline{c}_\theta, \overline{c}_\theta < \infty$ depend only on $\theta$; see for instance [6, Lem. 3.1.3].

Let us recall from (A.8) below that for $u, s > 0$, we have the continuous embedding

$$\left(\mathsf{H}_r^s(\mathbb{R}^d), \|\cdot\|_{\mathsf{H}_r^s(\mathbb{R}^d)}\right) \hookrightarrow \left(\mathsf{H}_r^{-u}(\mathbb{R}^d), \|\cdot\|_{\mathsf{H}_r^{-u}(\mathbb{R}^d)}\right) \tag{5.7}$$

so that real interpolation theory as presented above can be applied to the scale of Bessel potential spaces. However, the $K$-interpolation space $(\mathsf{H}_r^{-u}, \mathsf{H}_r^s)_{\theta,p;K}$ thus obtained is no longer a Bessel potential space, but rather a space in the Besov family; see e.g. [6, Thm. 6.4.5] or [52, Sec. 2.4.2]. To avoid having to introduce this family, we present the following lemma which allows us to return to the Bessel potential class by slightly lowering the parameter governing the regularity of functions.

**Lemma 5.2** | *Let $s, u > 0$ and $1 < r < 2$. Then for any $p \in [1, \infty)$, $\theta \in (0, 1)$ and $w < w_0 := \theta s - (1-\theta)u$ the intermediate space $(\mathsf{H}_r^{-u}, \mathsf{H}_r^s)_{\theta,p;K}$ continuously embeds into $\mathsf{H}_r^w$, in symbols*

$$(\mathsf{H}_r^{-u}, \mathsf{H}_r^s)_{\theta,p;K} \hookrightarrow \mathsf{H}_r^w.$$

*Equivalently, there exists a constant $c = c(p, r, w_0, w) < \infty$ such that*

$$\|f\|_{\mathsf{H}_r^w} \leq c \|f\|_{\theta,p,(\mathsf{H}_r^{-u},\mathsf{H}_r^s);K}$$

*for all $f \in (\mathsf{H}_r^{-u}, \mathsf{H}_r^s)_{\theta,p;K}$.*

**Proof** Details are found in Appendix B. □

## 5.2 Proof of Proposition 5.1

With the background from interpolation theory collected above, the proof of Proposition 5.1 is straightforward. We first observe that in order to deduce (5.3), we can use Lemma 5.2 to bound an appropriate $K$-interpolation norm. This in turn is more easily achieved by considering the equivalent discrete representation (5.5). We identify for each $n \in \mathbb{Z}$ in the sum in (5.5) a good decomposition $\lambda = a_\varepsilon + \varepsilon_\varepsilon$ in (5.1) by choosing $\varepsilon = \varepsilon_n$ appropriately. This in turn allows to obtain control over the $K$-functional $K(\varepsilon_n, \lambda)$ in (5.4) by using the hypothesized quantitative information (5.2) that we have on $a_\varepsilon$ and $e_\varepsilon$.

**Proof of Proposition 5.1** Fix $s > 0$ and $r \in (1, \infty)$. Consider real numbers $0 < \theta < 1$ and $1 \leq p < \infty$ to be chosen below. From the continuous embedding in Lemma 5.2 and the equivalence (5.6) of the discrete and continuous interpolation norms, we have for any

$$w < w_0(\theta) := \theta s - (1-\theta)u = \theta s - (1-\theta)(1 + d/r') \tag{5.8}$$

the bound

$$\|\lambda\|_{\mathsf{H}_r^w} \leq c(p, r, w_0, w) \|\lambda\|_{\theta,p,(\mathsf{H}_r^{-u},\mathsf{H}_r^s);K} \leq c(\theta) \|\lambda\|_{\theta,p,(\mathsf{H}_r^{-u},\mathsf{H}_r^s);k}. \tag{5.9}$$

We can thus aim to bound the right-hand side to gain control over the $\mathsf{H}_r^w$-norm of $\lambda$.

For part 1), we use the definition (5.5) of the discrete interpolation norm to see that

$$\|\lambda\|_{\theta,p,(\mathsf{H}_r^{-u},\mathsf{H}_r^s);k} = \left(\sum_{n=-\infty}^{\infty} 2^{-pn\theta} K(2^n, z)^p\right)^{1/p} \leq c(p)(K_- + K_+) \tag{5.10}$$

with

$$K_- := \left( \sum_{n \leq 0} 2^{-pn\theta} K(2^n, \lambda)^p \right)^{\frac{1}{p}}, \tag{5.11}$$

$$K_+ := \left( \sum_{n > 0} 2^{-pn\theta} K(2^n, \lambda)^p \right)^{\frac{1}{p}}. \tag{5.12}$$

Now let $(\varepsilon_n)_{n \in \mathbb{Z}}$ be a sequence of real numbers to be chosen below with $\varepsilon_n \leq \varepsilon_0$.

We start by controlling $K_-$. To this end, we first bound the $K$-functional (5.4) appearing in the sum in $K_-$ by using the assumed decomposition $\lambda = a_\varepsilon + e_\varepsilon$ for $\varepsilon \in (0, \varepsilon_0]$ and the inequality $(x+y)^p \leq c(p)(x^p + y^p)$ for $x, y \geq 0$ and $p \geq 1$ to find for all $t > 0$ that

$$K(t, \lambda)^p \leq (\|\lambda - a_\varepsilon\|_{\mathsf{H}_r^{-u}} + t\|a_\varepsilon\|_{\mathsf{H}_r^s})^p \tag{5.13}$$

$$\leq c(p)(\|e_\varepsilon\|^p_{\mathsf{H}_r^{-u}} + t^p \|a_\varepsilon\|^p_{\mathsf{H}_r^s}). \tag{5.14}$$

We can now bound $K_-$ in (5.11) by first changing the index of the summation from $-n$ to $n$ and then using the estimate from (5.14) to see that

$$K_-^p = \sum_{n=0}^{\infty} 2^{pn\theta} K(2^{-n}, \lambda)^p \leq c(p) \left( \sum_{n=0}^{\infty} 2^{pn\theta} \|e_{\varepsilon_{-n}}\|^p_{\mathsf{H}_r^{-u}} + \sum_{n=0}^{\infty} 2^{-pn(1-\theta)} \|a_{\varepsilon_{-n}}\|^p_{\mathsf{H}_r^s} \right). \tag{5.15}$$

Inserting the assumed norm bounds for $a_{\varepsilon_{-n}}$ and $e_{\varepsilon_{-n}}$ from (5.2) for the decomposition $\lambda = a_{\varepsilon_{-n}} + e_{\varepsilon_{-n}}$ gives

$$\sum_{n=0}^{\infty} 2^{pn\theta} \|e_{\varepsilon_{-n}}\|^p_{\mathsf{H}_r^{-u}} \leq c_e \sum_{n=0}^{\infty} 2^{pn\theta} \varepsilon_{-n}^{p(1+\xi)/2},$$

$$\sum_{n=0}^{\infty} 2^{-pn(1-\theta)} \|a_{\varepsilon_{-n}}\|^p_{\mathsf{H}_r^s} \leq c_a \sum_{n=0}^{\infty} 2^{-pn(1-\theta)} \varepsilon_{-n}^{-p(d/r'+s)/2}.$$

We next choose $(\varepsilon_{-n})_{n \in \mathbb{N}}$ and thus also the decomposition used in (5.13) and (5.15). For this, we first fix $\xi_0 < \xi$ and then set

$$\varepsilon_{-n} := 2^{-(2n\theta)/(1+\xi_0)} \varepsilon_0 \qquad \text{for all } n \in \mathbb{N}. \tag{5.16}$$

This choice ensures that $\varepsilon_{-n} \leq \varepsilon_0$ for all $n \geq 0$, so that the decomposition (5.1) is defined, (5.13) makes sense and the bound in (5.14) is valid. Moreover, we see that

$$\sum_{n=0}^{\infty} 2^{pn\theta} \varepsilon_{-n}^{p(1+\xi)/2} = \varepsilon_0^{p(1+\xi)/2} \sum_{n=0}^{\infty} 2^{-(\xi-\xi_0)pn\theta} < \infty, \tag{5.17}$$

and so the first infinite sum on the right-hand side of (5.15) is finite. Choosing $\theta$ with

$$0 < \theta < \bar{\theta}(\xi_0) := \frac{1+\xi_0}{1+\xi_0+s+d/r'} < 1 \tag{5.18}$$

gives $-\theta(d/r'+s)/(1+\xi_0) + (1-\theta) > 0$; so using the same sequence $(\varepsilon_{-n})_{n \in \mathbb{N}}$ leads to

$$\sum_{n=0}^{\infty} 2^{-pn(1-\theta)} \varepsilon_{-n}^{-p(d/r'+s)/2} = \varepsilon_0^{-p(d/r'+s)/2} \sum_{n=0}^{\infty} 2^{-np((1-\theta)-\theta(d/r'+s)/(1+\xi_0))} < \infty, \tag{5.19}$$

and also the second infinite sum in (5.15) is finite. With our choices of $(\varepsilon_{-n})_{n\in\mathbb{N}}$ and $\theta$ from (5.16) and (5.18), we find by combining (5.15) with (5.17) and (5.19) that

$$K_- \leq c_-(c_a\varepsilon_0^{-p(d/r'+s)/2} + c_e\varepsilon_0^{p(1+\xi)/2})^{1/p}, \tag{5.20}$$

where $c_- = c_-(p,\xi,\xi_0,d,\theta,r,s) < \infty$ is a constant that is independent of $\varepsilon_0$, $\lambda$ and $c_a$, $c_e$.

We next control the contribution of $K_+$ in (5.12). Note that for any $1 < r < \infty$ and $s, u > 0$, we have $\|a_{\varepsilon_0}\|_{\mathsf{H}_r^{-u}} \leq c_{r,s}\|a_{\varepsilon_0}\|_{\mathsf{H}_r^s}$ by (5.7). With this in mind, the triangle inequality and the assumed norm bounds from (5.2) imply

$$\|\lambda\|_{\mathsf{H}_r^{-u}} \leq \|a_{\varepsilon_0}\|_{\mathsf{H}_r^{-u}} + \|e_{\varepsilon_0}\|_{\mathsf{H}_r^{-u}} \leq c_a c_{r,s}\varepsilon_0^{-p(d/r'+s)/2} + c_e\varepsilon_0^{p(1+\xi)/2}.$$

In particular, $\lambda = a_{\varepsilon_0} + e_{\varepsilon_0} \in \mathsf{H}_r^{-u}$, so that the trivial decomposition $\lambda = \lambda + 0$ can be used to bound the $K$-functional by $K(t,\lambda)^p \leq (\|\lambda + 0\|_{\mathsf{H}_r^{-u}} + t\|0\|_{\mathsf{H}_r^s})^p = \|\lambda\|_{\mathsf{H}_u^{-r}}^p$ for all $t > 0$. Therefore, we get

$$K_+ \leq \bigg(\sum_{n=1}^\infty 2^{-pn\theta}\bigg)^{\frac{1}{p}}\|\lambda\|_{\mathsf{H}_r^{-u}} \leq c_+(c_a\varepsilon_0^{-(d/r'+s)/2} + c_e\varepsilon_0^{(1+\xi)/2})^{1/p} < \infty, \tag{5.21}$$

where $c_+ = c_+(p,\theta,r,s) := (1+c_{r,s})^{1/p}(\sum_{n=1}^\infty 2^{-pn\theta})^{1/p} < \infty$, which is clearly independent of $\varepsilon_0$, $\lambda$ and $c_a$, $c_e$.

We can now combine (5.10)–(5.12) with (5.20) and (5.21) to find

$$\|\lambda\|_{\theta,p,(\mathsf{H}_r^{-u},\mathsf{H}_r^s);k} \leq c_k(c_a\varepsilon_0^{-p(d/r'+s)/2} + c_e\varepsilon_0^{p(1+\xi)/2})^{1/p}$$

for some constant $c_k = c_k(p,\xi,\xi_0,d,\theta,r,s) < \infty$. The relation (5.9) between the various norms gives in turn that

$$\|\lambda\|_{\mathsf{H}_r^w} \leq c(c_a^{1/p}\varepsilon_0^{-(d/r'+s)/2} + c_e^{1/p}\varepsilon_0^{(1+\xi)/2}) \tag{5.22}$$

for all $w < w_0$, with a constant $c = c(p,\xi,\xi_0,d,\theta,r,s,w) < \infty$. We now see that $p$ is a nuisance parameter which can be chosen as $p = 1$. This establishes the bound (5.3).

For part 2), we must still establish the addendum that given $\xi > 0$ and $d \geq 1$, there exist choices of $1 < r < 2$ and $s > 0$ as well as subsequent choices $0 < \xi_0 < \xi$ and $0 < \theta < \bar{\theta}(\xi_0)$ such that $w_0(\theta) > 0$. To see this, we start with the definition of $w_0(\theta)$ from (5.8) and claim that it is enough to show that

$$\bar{\theta}(\xi_0)s - \big(1 - \bar{\theta}(\xi_0)\big)(d/r' + 1) > 0. \tag{5.23}$$

Indeed, if (5.23) holds, then choosing $\theta$ sufficiently close to but smaller than $\bar{\theta}(\xi_0)$ retains positivity. The inequality (5.23) on the other hand is easily seen to be equivalent to

$$(1+\xi_0)s > (d/r'+1)(s+d/r'), \tag{5.24}$$

using the definition of $\bar{\theta}(\xi_0)$ in (5.18). But for any $s > 0$, we can choose $1 < r < 2$ sufficiently small, or equivalently $2 < r' < \infty$ sufficiently large, to validate the inequality (5.24). In particular, we can choose $s > 0$ sufficiently small and $r' > 2$ sufficiently large to have $\gamma = (d/r' + s)/2 < 1$.

For part 3), we simply chase our chosen parameters and the decomposition we used through the proof. We first fix $\xi_0 \in (0, \xi)$ as, say, the midpoint $\xi_0^* := \xi/2$, and then $r$ and $s$ to satisfy (5.24). Note that as seen just above, the choices of $r$ and $s$ can be made to guarantee that $\gamma = (d/r' + s)/2 < 1$. The parameters $\xi_0^*$, $r$ and $s$ determine $\bar{\theta}(\xi_0^*)$ via (5.18), which depends on $r$ and $s$. This allows us to subsequently choose $\theta^* \in (0, \bar{\theta}(\xi_0^*))$ as, say, the midpoint $\theta^* := \bar{\theta}(\xi_0^*)/2$. Having chosen all these parameters, our bounds for $K_-$ and $K_+$ from (5.11) and (5.12) depend only on the decomposition and norm bounds for $\varepsilon = \varepsilon_0$ and $\varepsilon = \varepsilon_n$ given in (5.16); see (5.15) and (5.21), respectively. Thus $\alpha_0(\xi, r, s) := 2\theta^*/(1 + \xi_0^*)$ is the function whose existence we needed to show. This establishes part 3) and completes the proof. $\square$

**Remark 5.3** | The exponent $-(d/r' + s)$ on $\varepsilon_0$ we obtained in (5.3) is most likely not optimal. In the proof, we could instead have defined $\varepsilon_n = 2^{-(2n\theta)/(1+\xi_0)} \varepsilon_0^\rho$ for some additional parameter $\rho \in \mathbb{R}$. Following the same steps as in the proof, $\rho$ then appears in the exponent of the bounds on $K_-$ and $K_+$. We could then optimize over $\rho$ to improve the bounds on $K_-$ and $K_+$, obtaining $\rho^*$ and $\|\lambda\|_{\theta, p, (\mathsf{H}_r^{-u}, \mathsf{H}_r^s); K} \leq \tilde{c} \varepsilon_0^{-\rho^*(d/r'+s)/2}$ for $\varepsilon_0 \leq 1$.

In their seminal contribution, Fournier and Printems [25] showed that certain diffusive dynamics $S = (S_t)_{t \in [0,T]}$ have time-marginals $\lambda_t := \mathrm{Law}(S_t)$ for $t > 0$ which are absolutely continuous with respect to Lebesgue measure; see [25, Thm. 1.2]. We outline their strategy in Example 7.2 below and discuss the array of subsequent developments sparked by their idea in greater detail in Section 7.2. Proposition 5.1 is one of these subsequent developments, and in Remark 7.3 we argue that (5.24) leads to the assumptions in [25].

## 6 PROOF OF THE MAIN RESULT

We can now finally prove our main result in Theorem 1.5, which we repeat for ease of reference.

**Theorem 1.5** | *Impose Assumptions 1.1 and 1.2. Then there exist real numbers $w > 0$ and $r > 1$ such that for any cluster point $\mathbf{P}_{x_0}^\infty$ of the sequence $(\mathbf{P}_{x_0}^n)_{n \in \mathbb{N}}$, we have*

$$\mathbf{P}_{x_0}^\infty \left[ \left\{ \omega \in \Omega \ : \ \frac{\mathrm{d}\Lambda_t(\omega)}{\mathrm{d}x} = \bar{p}(t, \omega) \text{ is in } \mathsf{H}_r^w(\mathbb{R}^d) \text{ for almost all } t \in (0, T] \right\} \right] = 1, \quad (1.11)$$

*where $\bar{p} : (0, T] \times \Omega \to \mathsf{H}_r^w(\mathbb{R}^d)$ is a strongly measurable function.*

*More precisely, for each $t \in (0, T]$, we have $\mathbf{P}_{x_0}^\infty$-a.s. that $\Lambda_t \ll \mathrm{d}x$, and for any $q \geq 1$, $w$ and $r$ can be chosen in such a way that there exists $1 > \gamma > 0$ such that the function $t \mapsto (\omega \mapsto \mathrm{d}\Lambda_t(\omega)/\mathrm{d}x)$ defines a strongly measurable map*

$$p : (0, T] \to \mathbb{L}^q\big((\Omega, \mathcal{F}, \mathbf{P}_{x_0}^\infty); \mathsf{H}_r^w(\mathbb{R}^d)\big) \quad (1.12)$$

*which satisfies for some constant $c_{\mathrm{Thm.\ 1.5}} < \infty$ the bound*

$$\Big\| \|p(t)\|_{\mathsf{H}_r^w(\mathbb{R}^d)} \Big\|_{\mathbb{L}^q(\mathbf{P}_{x_0}^\infty)} \leq c_{\mathrm{Thm.\ 1.5}} (1 \wedge t)^{-\gamma} \quad (1.13)$$

*for all $t \in (0, T]$. In (1.13), $c_{\text{Thm. 1.5}}$ depends on $r$, $w$, $q$, $\gamma$, but not on $t$ and neither on $\mathbf{P}_{x_0}^\infty$.*

*Finally, $\bar{p}$ is unique up to $(\mathrm{d}t \otimes \mathrm{d}\mathbf{P}_{x_0}^\infty)$-a.e. equality and we have that $\bar{p}(t, \cdot) = p(t)$ in $\mathbb{L}^q((\Omega, \mathcal{F}, \mathbf{P}_{x_0}^\infty); \mathsf{H}_r^w(\mathbb{R}^d, \mathrm{d}x))$, for almost every $t \in (0, T]$. In particular, $\bar{p}(t, \omega)$ is the density of $\Lambda_t(\omega)$ $(\mathrm{d}t \otimes \mathrm{d}\mathbf{P}_{x_0}^\infty)$-a.e., and $\bar{p}(t, \cdot)$ satisfies the bound (1.13) for almost every $t \in (0, T]$.*

Before giving the proof, we outline the basic idea and clarify some notation. Theorem 1.5 is a statement applying to *all* cluster points of the sequence $(\mathbf{P}_{x_0}^n)_{n\in\mathbb{N}}$ of empirical measure flow laws. Let us therefore consider an arbitrary cluster point $\mathbf{P}_{x_0}^\infty$ of $(\mathbf{P}_{x_0}^n)_{n\in\mathbb{N}}$ and a subsequence $(n_k)_{k\in\mathbb{N}}$ with the property that $(\mathbf{P}_{x_0}^{n_k})_{k\in\mathbb{N}}$ converges to $\mathbf{P}_{x_0}^\infty$ narrowly in $\mathbf{M}_1^+(\Omega)$ as $k \to \infty$. Here $(\Omega, \mathcal{F})$ is the canonical setup from Section 1.3 with $\Omega = \mathbf{C}([0, T]; \mathcal{S}'(\mathbb{R}^d))$ and $\mathcal{F}$ the Borel-$\sigma$-algebra on $\Omega$. This is exactly the setting of (2.35), which is the first part of the two-step construction (2.35), (2.36) outlined in Section 2.3. Continuing with the second step (2.36), we now choose an arbitrary cluster point $\bar{\mathbf{P}}_{x_0}^\infty$ of the subsequence $(\bar{\mathbf{P}}_{x_0}^{n_k})_{k\in\mathbb{N}}$, where we insist that $(n_k)_{k\in\mathbb{N}}$ is as before. We also choose an arbitrary subsubsequence $(n_{k'})_{k'\in\mathbb{N}}$ of $(n_k)_{k\in\mathbb{N}}$ such that $(\bar{\mathbf{P}}_{x_0}^{n_{k'}})_{k'\in\mathbb{N}}$ converges to $\bar{\mathbf{P}}_{x_0}^\infty$ narrowly in $\mathbf{M}_1^+(\bar{\Omega})$ as $k' \to \infty$. This is now defined on the extended setup $(\bar{\Omega}, \bar{\mathcal{F}})$, introduced in Section 1.5 as

$$\bar{\Omega} = \mathbf{C}\Big([0, T]^\Delta; \mathcal{S}'(\mathbb{R}^d) \times \mathcal{S}'(\mathbb{R}^d) \times \mathbb{R}^m\Big)$$

and where $\bar{\mathcal{F}}$ is the Borel-$\sigma$-algebra on $\bar{\Omega}$. On $\Omega$ and $\bar{\Omega}$, we have the respective canonical processes $(\Lambda_t)_{t\in[0,T]}$ and $(\bar{\Lambda}_t, A_{\varepsilon,t}, Z_t)_{(\varepsilon,t)\in[0,T]^\Delta}$. The latter is used in (1.16) to define the decomposition

$$\bar{\Lambda}_t = A_{\varepsilon,t} + E_{\varepsilon,t} \quad \text{in } \mathcal{S}' \text{ for } (\varepsilon, t) \in [0, T]^\Delta. \tag{6.1}$$

Moreover, as in (2.37), we have by construction the identity

$$\text{Law}_{\bar{\mathbf{P}}_{x_0}^\infty}(\bar{\Lambda}) = \text{Law}_{\mathbf{P}_{x_0}^\infty}(\Lambda). \tag{6.2}$$

The reason we iteratively chose subsequences $(n_{k'})_{k'\in\mathbb{N}}$ of $(n_k)_{k\in\mathbb{N}}$, and $(n_k)_{k\in\mathbb{N}}$ of $(n)_{n\in\mathbb{N}}$ is to guarantee that (6.2) holds, which allows us to study the initially fixed cluster point $\mathbf{P}_{x_0}^\infty$ through the lens of $\bar{\mathbf{P}}_{x_0}^\infty$. This is achieved as follows.

We first appeal to Propositions 3.2 and 4.1 to obtain under $\bar{\mathbf{P}}_{x_0}^\infty$ estimates for $A_{\varepsilon,t}$ in terms of the $\mathsf{H}_r^s(\mathbb{R}^d)$-norm, and for $E_{\varepsilon,t}$ in terms of the $\mathsf{H}_r^{-u}(\mathbb{R}^d)$-norm. For appropriate values of $r$, $s$, $u$, the different rates at which the norms explode respectively vanish as $\varepsilon \to 0$, for fixed $t \in (0, T]$, can be exploited via the decomposition $\bar{\Lambda}_t = A_{\varepsilon,t} + E_{\varepsilon,t}$ in (6.1) and the abstract interpolation result in Proposition 5.1 to obtain, again under $\bar{\mathbf{P}}_{x_0}^\infty$, an intermediate-norm bound on $\bar{\Lambda}_t$; see Step 1 of the proof below.

The extension $(\bar{\Omega}, \bar{\mathcal{F}})$ of the original setup $(\Omega, \mathcal{F})$ is only a means to an end. We next use (6.2) to project back onto the original setup $(\Omega, \mathcal{F})$ and obtain from the norm bound for $\bar{\Lambda}_t$ under $\bar{\mathbf{P}}_{x_0}^\infty$ on $(\bar{\Omega}, \bar{\mathcal{F}})$ a norm bound for $\Lambda_t$ under $\mathbf{P}_{x_0}^\infty$, as in the statement of Theorem 1.5; see Step 2 of the proof below. For this step, we make the following useful observation.

Denote by $\bar{\mathbb{F}} := \mathbb{F}^{\bar{\Lambda}}$ the right-continuous version of the filtration generated by $\bar{\Lambda}$ on $\bar{\Omega}$ and write $\bar{\mathbb{F}} = (\bar{\mathcal{F}}_t)_{t\in[0,T]}$. From Section 1.3, we also have the right-continuous version of the filtration generated by $\Lambda$ on $\Omega$, which we denote by $\mathbb{F} = (\mathcal{F}_t)_{t\in[0,T]}$. To avoid confusion, we point out that $\mathcal{F}_T = \mathcal{F}$ by definition, but that $\bar{\mathcal{F}}_T \subsetneq \bar{\mathcal{F}}$. Note that $\bar{\mathbb{F}}$ can be identified with $\mathbb{F}$; indeed, every measurable set $B' \in \bar{\mathcal{F}}_t$ has the form $B' = \{\bar{\omega} \in \bar{\Omega} : \bar{\omega}_{0,\cdot}^1 \in B\}$ for some measurable set $B \in \mathcal{F}_t$ and vice versa. In consequence, we can identify any $\bar{\mathcal{F}}_t$-measurable random variable with an $\mathcal{F}_t$-measurable random variable and vice versa.

**Proof of Theorem 1.5** Let $\beta > 0$ be the Hölder-exponent of $\sigma$ and $\bar{\sigma}$ from Assumption 1.2(H). Fix $\xi := \beta/2$ and choose $q \geq 1$. Consider moreover $1 < r < 2$ and $s > 0$ as well as $\alpha > 0$, all to be chosen at the end of Step 1 below. We denote throughout by $r'$ the conjugate exponent of $r$ and set $u := 1 + d/r'$.

**Step 1** *Pointwise interpolation argument.* For all $\bar{\omega} \in \bar{\Omega}$, (1.16) gives the decomposition

$$\bar{\Lambda}_t(\bar{\omega}) = A_{\varepsilon,t}(\bar{\omega}) + E_{\varepsilon,t}(\bar{\omega}) \quad \text{in } \mathcal{S}' \text{ for } (\varepsilon, t) \in [0, T]^\Delta. \tag{6.3}$$

We first fix $(\varepsilon, t) \in [0, T]^\Delta$ with $0 < \varepsilon \leq 1 \wedge t$. In Proposition 3.2, we found a constant $c_{r,s} < \infty$, depending only on $s$ and $r$ but not on any other quantity chosen thus far, such that $\bar{\mathbf{P}}_{x_0}^\infty$-a.s.,

$$\|A_{\varepsilon,t}(\bar{\omega})\|_{\mathsf{H}_r^s(\mathbb{R}^d)} \leq c_{r,s}\varepsilon^{-(d/r'+s)/2}. \tag{6.4}$$

We can thus choose a $\bar{\mathbf{P}}_{x_0}^\infty$-nullset $N_{\varepsilon,t}^A \in \bar{\mathcal{F}}$ such that (6.4) holds for all $\bar{\omega} \in \bar{\Omega} \backslash N_{\varepsilon,t}^A$. Next, recall the auxiliary countable set

$$(0, 1 \wedge t]_\alpha = \{2^{-\alpha n}(1 \wedge t) : n \in \mathbb{N}_0\}$$

from the proof of Proposition 4.1. In that result, we found a real-valued random variable

$$\bar{C}_{t,r,\alpha,\xi} \in \mathbb{L}^q(\bar{\Omega}, \bar{\mathcal{F}}, \bar{\mathbf{P}}_{x_0}^\infty) \tag{6.5}$$

with the property that $\bar{\mathbf{P}}_{x_0}^\infty$-a.s.,

$$\|E_{\varepsilon,t}(\bar{\omega})\|_{\mathsf{H}_r^{-u}(\mathbb{R}^d)} \leq \bar{C}_{t,r,\alpha,\xi}\,\varepsilon^{(1+\xi)/2} \quad \text{simultaneously for all } \varepsilon \in (0, 1 \wedge t]_\alpha. \tag{6.6}$$

There thus exists a $\bar{\mathbf{P}}_{x_0}^\infty$-nullset $N_{\alpha,t}^E \in \bar{\mathcal{F}}$ such that (6.6) holds and $\bar{C}_{t,r,\alpha,\xi}(\bar{\omega}) < \infty$ for all $\bar{\omega} \in \bar{\Omega} \backslash N_{\alpha,t}^E$. It is important to note that the same choice of $r$, $s$ and $\alpha$ works simultaneously for all $\bar{\omega}$ and $t$.

For fixed $\bar{\omega}$, we now apply part 3) of Proposition 5.1, using $\varepsilon_0 = 1 \wedge t$ and $\xi = \beta/2$, the decomposition (6.3) for (5.1), and the norm bounds (6.4) and (6.6) for (5.2). We therefore find real numbers $r \in (1, 2)$ and $s > 0$ with $\gamma = (d/r' + s)/2 < 1$, and subsequently fix $\alpha := \alpha_0(\xi, r, s) > 0$ with $\alpha_0$ as in Proposition 5.1. As seen above, $r$, $s$ and $\alpha$ can be chosen to not depend on $(t, \bar{\omega}) \in (0, T] \times \bar{\Omega}$ and thus in particular to be deterministic. We fix such a choice. Next consider

$$\bar{\Omega}_t := \bar{\Omega} \backslash \left( N_{\alpha,t}^E \cup \left( \bigcup_{\varepsilon \in (0, 1 \wedge t]_\alpha} N_{\varepsilon,t}^A \right) \right). \tag{6.7}$$

Since $(0, 1 \wedge t]_\alpha$ is countable, we have $\bar{\mathbf{P}}_{x_0}^\infty[\bar{\Omega}_t] = 1$. As the bounds for $A_{\varepsilon,t}$ and $E_{\varepsilon,t}$ in (6.4) and (6.6) hold for all $\bar{\omega} \in \bar{\Omega}_t$, part 3) of Proposition 5.1 and the inequality $(1 \wedge t)^{(1+\xi)/2} \leq 1 \leq (1 \wedge t)^{-(d/r'+s)/2}$ imply that

$$\|\bar{\Lambda}_t(\bar{\omega})\|_{\mathsf{H}_r^w(\mathbb{R}^d)} \leq c\Big(c_{r,s}(1 \wedge t)^{-(d/r'+s)/2} + \bar{C}_{t,r,\alpha,\xi}(\bar{\omega})(1 \wedge t)^{(1+\xi)/2}\Big)$$
$$\leq c\Big(c_{r,s} + \bar{C}_{t,r,\alpha,\xi}(\bar{\omega})\Big)(1 \wedge t)^{-(d/r'+s)/2}$$

for all $w < w_0(d, \xi, s, r)$ with $w_0$ from Proposition 5.1. Since $\xi$, $s$, $r$ do not depend on $(t, \bar{\omega})$, neither do $\alpha$ and $w_0$; so we can and do choose $w$ to be independent of $(t, \bar{\omega})$. Thus also $c = c(w, d, \xi, s, r)$ from Proposition 5.1 is non-random. Upon writing

$$\bar{C}_t(\bar{\omega}) := c(c_{r,s} + \bar{C}_{t,r,\alpha,\xi}(\bar{\omega})) \tag{6.8}$$

for $\bar{\omega} \in \bar{\Omega}_t$, we have shown for our fixed $t \in (0, T]$ that $\bar{\mathbf{P}}_{x_0}^\infty$-a.s.,

$$\|\bar{\Lambda}_t\|_{\mathsf{H}_r^w(\mathbb{R}^d)} \leq \bar{C}_t(1 \wedge t)^{-\gamma}, \tag{6.9}$$

where $\bar{C}_t$ is in $\mathbb{L}^q(\bar{\mathbf{P}}_{x_0}^\infty)$ in view of its definition (6.8) and the bound (6.5). In fact, by Proposition 4.1, we have $\sup\{\|\bar{C}_{t,r,\alpha,\xi}\|_{\mathbb{L}^q(\bar{\Omega}, \bar{\mathcal{F}}, \bar{\mathbf{P}}_{x_0}^\infty)} : t \in [0, T]\} < \infty$. Hence

$$\sup_{t \in (0, T]} \|\bar{C}_t\|_{\mathbb{L}^q(\bar{\Omega}, \bar{\mathcal{F}}, \bar{\mathbf{P}}_{x_0}^\infty)} < \infty. \tag{6.10}$$

**Step 2** *Projection to original setup.* To transfer our findings from $(\bar{\Omega}, \bar{\mathcal{F}})$ to $(\Omega, \mathcal{F})$, define

$$C_t' := \mathbf{E}^{\bar{\mathbf{P}}_{x_0}^\infty}[\bar{C}_t \mid \bar{\mathcal{F}}_t]$$

and identify this with an $\mathcal{F}_t$-measurable random variable $C_t$. Because (6.9) implies via conditioning on $\bar{\mathcal{F}}_t$ that $\|\bar{\Lambda}_t\|_{\mathsf{H}_r^w(\mathbb{R}^d)} \leq C_t'(1 \wedge t)^{-\gamma}$ $\bar{\mathbf{P}}_{x_0}^\infty$-a.s. and because the definition of $\bar{\Lambda}_t$ in (1.14) gives $\mathrm{Law}_{\bar{\mathbf{P}}_{x_0}^\infty}(\bar{\Lambda}, C_t') = \mathrm{Law}_{\mathbf{P}_{x_0}^\infty}(\Lambda, C_t)$, we first get

$$\|\Lambda_t\|_{\mathsf{H}_r^w(\mathbb{R}^d)} \leq C_t(1 \wedge t)^{-\gamma} \qquad \mathbf{P}_{x_0}^\infty\text{-a.s.}, \tag{6.11}$$

then via Jensen's inequality also

$$\mathbf{E}^{\mathbf{P}_{x_0}^\infty}[C_t^q] = \mathbf{E}^{\bar{\mathbf{P}}_{x_0}^\infty}[(C_t')^q] \leq \mathbf{E}^{\bar{\mathbf{P}}_{x_0}^\infty}[\bar{C}_t^q] < \infty \tag{6.12}$$

and finally in view of (6.10) even the uniform bound

$$c_{\text{Thm. 1.5}} := \sup_{t \in (0, T]} \|C_t\|_{\mathbb{L}^q(\Omega, \mathcal{F}_t, \mathbf{P}_{x_0}^\infty)} < \infty. \tag{6.13}$$

Together (6.11) and (6.13) show for our fixed $t$ that

$$\left\|\|\Lambda_t\|_{\mathsf{H}_r^w(\mathbb{R}^d)}\right\|_{\mathbb{L}^q(\Omega, \mathcal{F}_t, \mathbf{P}_{x_0}^\infty)} \leq c_{\text{Thm. 1.5}}(1 \wedge t)^{-\gamma}. \tag{6.14}$$

**Step 3** *Absolute continuity and density.* In (6.14), we have obtained a norm bound for $\Lambda_t$ viewed as an $\mathcal{S}'(\mathbb{R}^d)$-valued random variable. We now want to connect it to a random function obtained as the Radon–Nikodým derivative with respect to Lebesgue measure when we view $\Lambda_t$ as an $\mathbf{M}_1^+(\mathbb{R}^d)$-valued random variable. Due to (1.5) and (6.14), we can find for each $t \in (0, T]$ a $\mathbf{P}_{x_0}^\infty$-nullset $N_t \in \mathcal{F}$ such that $\Lambda_t(\omega) \in \mathbf{M}_1^+(\mathbb{R}^d) \subseteq \mathcal{S}'(\mathbb{R}^d)$ and $\|\Lambda_t(\omega)\|_{\mathsf{H}_r^w(\mathbb{R}^d, \mathrm{d}x)} < \infty$ for all $\omega \in \Omega \setminus N_t$. The duality bound

$$|\Lambda_t(\omega)[\phi]| = |\langle \Lambda_t(\omega); \phi\rangle_{\mathcal{S}' \times \mathcal{S}}| \leq \|\Lambda_t(\omega)\|_{\mathsf{H}_r^w(\mathbb{R}^d)} \|\phi\|_{\mathsf{H}_{r'}^{-w}(\mathbb{R}^d)} \qquad \text{for all } \phi \in \mathcal{S}$$

implies that for each $\omega \in \Omega \backslash N_t$ there there exists a d$x$-a.e. unique measurable function $q(t, \omega) : \mathbb{R}^d \to [0, \infty]$ such that $q(t, \omega) \in \mathsf{H}_r^w(\mathbb{R}^d)$ and

$$\Lambda_t(\omega)[\phi] = \int_{\mathbb{R}^d} \phi(x) \big(q(t, \omega)(x)\big) \, \mathrm{d}x \qquad \text{for all } \phi \in \mathcal{S}\,. \tag{6.15}$$

Since $\Lambda_t(\omega) \in \mathbf{M}_1^+(\mathbb{R}^d)$, we have that $q(t, \omega) \in \mathsf{H}_r^w(\mathbb{R}^d) \cap \mathbb{L}^1(\mathbb{R}^d, \mathrm{d}x)$, and so $q(t, \omega)(x) < \infty$ d$x$-a.e. Via approximation, (6.15) extends to

$$\Lambda_t(\omega)[A] = \int_A \big(q(t, \omega)(x)\big) \, \mathrm{d}x \qquad \text{for all } A \in \mathcal{B}(\mathbb{R}^d)\,.$$

Hence $\Lambda_t(\omega) \ll \mathrm{d}x$ and $\mathrm{d}\Lambda_t(\omega)/\mathrm{d}x = q(t, \omega)$ in $\mathbb{L}^1(\mathrm{d}x)$ by the Radon–Nikodým theorem.

***Step 4*** *Measurability properties.* We next establish some required measurability properties by using results from analysis in Banach spaces. First, from Lemma 2.5, we have $\mathsf{H}_r^w(\mathbb{R}^d)' = \mathsf{H}_{r'}^{-w}(\mathbb{R}^d)$, where $r'$ is the conjugate exponent of $r$. Since the separable space $\mathcal{S}(\mathbb{R}^d)$ is norm-dense in $\mathsf{H}_r^w(\mathbb{R}^d)'$, it is also weak*-dense so that $\mathsf{H}_r^w(\mathbb{R}^d)'$ is weakly separable; see Bergh and Löfström [6, Thm. 6.2.3] and Rudin [47, Cor. to Thm. 3.12]. Let $N_t \in \mathcal{F}_t$ be as in Step 3. By (6.15), we have for each $t \in (0, T]$ and $\omega \in \Omega \backslash N_t$ that $\Lambda_t(\omega) \in \mathsf{H}_r^w(\mathbb{R}^d)$. In addition, for each $\phi \in \mathcal{S}$, the duality pairing $\Lambda_t(\omega)[\phi] = \langle q(t, \omega); \phi \rangle_{\mathcal{S}' \times \mathcal{S}} = \langle q(t, \omega); \phi \rangle_{\mathsf{H}_r^w \times (\mathsf{H}_r^w)'}$ takes the form (6.15). Hence for each $t \in (0, T]$, the map

$$\Omega \backslash N_t \ni \omega \mapsto \langle q(t, \omega); \phi \rangle_{\mathsf{H}_r^w \times (\mathsf{H}_r^w)'} \in \mathbb{R}$$

is $\mathcal{F}_t/\mathcal{B}(\mathbb{R}^d)$-measurable and thus also $\mathcal{F}/\mathcal{B}(\mathbb{R}^d)$-measurable, and Pettis' theorem now implies that $\Omega \backslash N_t \ni \omega \mapsto p(t)(\omega) := q(t, \omega) \in \mathsf{H}_r^w(\mathbb{R}^d)$ is strongly measurable; see Hytönen et al. [29, Thm. 1.1.20] or Theorem A.1. So together with Step 2, we have for each $t \in (0, T]$ that $p(t) \in \mathbb{L}^q((\Omega, \mathcal{F}, \mathbf{P}_{x_0}^\infty); \mathsf{H}_r^w(\mathbb{R}^d))$.

Similarly, by [29, Prop. 1.3.3], the dual of the Bochner space $\mathbb{L}^q((\Omega, \mathcal{F}, \mathbf{P}_{x_0}^\infty); \mathsf{H}_r^w(\mathbb{R}^d))$ is $\mathbb{L}^{q'}((\Omega, \mathcal{F}, \mathbf{P}_{x_0}^\infty); \mathsf{H}_r^w(\mathbb{R}^d)')$, where now $q'$ is the conjugate exponent of $q$. Since simple functions (finite linear combinations of functions of the form $\mathbb{1}_A \phi$ with $A \in \mathcal{F}$ and $\phi \in \mathcal{S}$) are norm-dense in $\mathbb{L}^{q'}((\Omega, \mathcal{F}, \mathbf{P}_{x_0}^\infty); \mathsf{H}_r^w(\mathbb{R}^d)')$, see [29, Lem. 1.2.19, Prop. 1.2.29], and since the map

$$(0, T] \ni t \mapsto \mathbf{E}^{\mathbf{P}_{x_0}^\infty}\left[\Lambda_t(\omega)[\phi] \mathbb{1}_A\right] = \mathbf{E}^{\mathbf{P}_{x_0}^\infty}[\langle q(t, \omega); \phi \rangle_{\mathsf{H}_r^w \times \mathsf{H}_{r'}^{-w}} \mathbb{1}_A] \in \mathbb{R}$$

is measurable, the same reasoning using Pettis' theorem shows that $t \mapsto p(t) = q(t, \cdot)$ is strongly measurable as a map $(0, T] \to \mathbb{L}^q((\Omega, \mathcal{F}, \mathbf{P}_{x_0}^\infty); \mathsf{H}_r^w(\mathbb{R}^d))$. Via (6.14) and (6.15), this gives (1.12) and (1.13).

***Step 5*** *Identification with a product-measurable function.* From Step 4, we have the map

$$p : (0, T] \to \mathbb{L}^q\big((\Omega, \mathcal{F}, \mathbf{P}_{x_0}^\infty); \mathsf{H}_r^w(\mathbb{R}^d)\big)\,,$$

which is strongly measurable. We now use [29, Prop. 1.2.25] to conclude that there exists a strongly measurable and $(\mathrm{d}t \otimes \mathrm{d}\mathbf{P}_{x_0}^\infty)$-a.e. unique function

$$\bar{p} : (0, T] \times \Omega \to \mathsf{H}_r^w(\mathbb{R}^d)$$

with the property that $p(t) = \bar{p}(t, \cdot)$ in $\mathbb{L}^q((\Omega, \mathcal{F}, \mathbf{P}_{x_0}^\infty); \mathsf{H}_r^w(\mathbb{R}^d))$, almost everywhere on $(0, T]$. As a particular consequence of this we get that $\bar{p}(t, \cdot)$ satisfies the bound in (1.13) $dt$-a.e. In addition, by Steps 3 and 4 we have $\mathbf{P}_{x_0}^\infty$-a.s. that $p(t)$ is the density $d\Lambda_t/dx$ for all $t \in (0, T]$. Thus we can use Tonelli's theorem to find for each $\phi \in \mathcal{S}(\mathbb{R}^d)$ that

$$\begin{aligned}
0 &= \int_0^T \mathbf{E}^{\mathbf{P}_{x_0}^\infty} \left[ \left| \Lambda_t(\omega)[\phi] - \int_{\mathbb{R}^d} \phi(x) \big( p(t)(\omega)(x) \big) \, dx \right| \right] dt \\
&= \int_0^T \mathbf{E}^{\mathbf{P}_{x_0}^\infty} \left[ \left| \Lambda_t(\omega)[\phi] - \int_{\mathbb{R}^d} \phi(x) \big( \bar{p}(t, \omega)(x) \big) \, dx \right| \right] dt \\
&= \mathbf{E}^{\mathbf{P}_{x_0}^\infty} \left[ \int_0^T \left| \Lambda_t(\omega)[\phi] - \int_{\mathbb{R}^d} \phi(x) \big( \bar{p}(t, \omega)(x) \big) \, dx \right| dt \right], \quad (6.16)
\end{aligned}$$

from which we deduce that $\bar{p}(t, \omega) = d\Lambda_t(\omega)/dx$ up to $dx$-a.e. equivalence, $dt \otimes d\mathbf{P}_{x_0}^\infty$-a.e. All this establishes the last part of Theorem 1.5, and it remains to prove (1.11), which is no longer difficult. We first observe that since $\gamma < 1$, we have

$$c_\gamma := \int_0^T (1 \wedge t)^{-\gamma} \, dt < \infty.$$

Now let $K > 0$. Using Markov's inequality, Tonelli's theorem, $q > 1$ and (6.14) gives

$$\begin{aligned}
\mathbf{P}_{x_0}^\infty \left[ \int_0^T \|\bar{p}(t, \cdot)\|_{\mathsf{H}_r^w} \, dt > K \right] &\leq \frac{1}{K} \mathbf{E}^{\mathbf{P}_{x_0}^\infty} \left[ \int_0^T \|\bar{p}(t, \cdot)\|_{\mathsf{H}_r^w} \, dt \right] \\
&= \frac{1}{K} \int_0^T \mathbf{E}^{\mathbf{P}_{x_0}^\infty} [\|\bar{p}(t, \cdot)\|_{\mathsf{H}_r^w}] \, dt \\
&= \frac{1}{K} \int_0^T \mathbf{E}^{\mathbf{P}_{x_0}^\infty} [\|p(t)\|_{\mathsf{H}_r^w}] \, dt \\
&\leq \frac{1}{K} \int_0^T \big\| \|p(t)\|_{\mathsf{H}_r^w} \big\|_{\mathbb{L}^q(\mathbf{P}_{x_0}^\infty)} \, dt \\
&\leq \frac{c_\gamma \, c_{\text{Thm. 1.5}}}{K}.
\end{aligned}$$

By letting $K \to \infty$, we deduce that

$$\mathbf{P}_{x_0}^\infty \left[ \int_0^T \|\bar{p}(t, \cdot)\|_{\mathsf{H}_r^w} \, dt < \infty \right] = 1,$$

which in turn implies that $\|\bar{p}(t, \cdot)\|_{\mathsf{H}_r^w} < \infty$ and thus $\bar{p}(t, \cdot) \in \mathsf{H}_r^w$, for almost every $t \in (0, T]$, $\mathbf{P}_{x_0}^\infty$-a.s. Since (6.16) also implies that $\bar{p}(t, \omega) = d\Lambda_t(\omega)/dx$ up to $dx$-a.e. equivalence, for almost every $t \in (0, T]$, $\mathbf{P}_{x_0}^\infty$-a.s., we have (1.11) and the proof is complete. □

**Remark 6.2** | Let us now clarify the dependence of the constant $c_{\text{Thm. 1.5}}$ appearing in Theorem 1.5 on the various quantities mentioned in Remark 1.6. First recall (6.13), namely

$$c_{\text{Thm. 1.5}} = \sup_{t \in (0, T]} \|C_t\|_{\mathbb{L}^q(\Omega, \mathcal{F}_t, \mathbf{P}_{x_0}^\infty)} < \infty, \quad (6.13)$$

and from (6.12) and (6.8) that we have

$$\|C_t\|_{\mathbb{L}^q(\Omega, \mathcal{F}_t, \mathbf{P}_{x_0}^\infty)} \leq \|\bar{C}_t\|_{\mathbb{L}^q(\bar{\Omega}, \bar{\mathcal{F}}, \bar{\mathbf{P}}_{x_0}^\infty)} \leq c \big( c_{r,s} + \|\bar{C}_{t,r,\alpha,\xi}\|_{\mathbb{L}^q(\bar{\Omega}, \bar{\mathcal{F}}, \bar{\mathbf{P}}_{x_0}^\infty)} \big).$$

This reveals the full dependence of $c_{\text{Thm. 1.5}}$ on all parameters, as follows:

1) The estimates for the regular and distributional parts combined in Step 1 of the proof via the interpolation result in Proposition 5.1 introduce a dependence on the dimension $d$ as well as on the quantities $r$ and $s$ fixed in Step 1 to satisfy $\gamma = (d/r' + s)/2 < 1$. In addition, we chose $\xi = \beta/2$ and $w < w_0(d, \xi, s, r)$, all of which enter the constant $c = c(w, d, \xi, s, r)$ from Proposition 5.1.

2) The constant $c_{r,s}$ from the bound on the regular part in (6.4) is obtained via Proposition 3.2. This uses Lemma 3.4, so that $c_{r,s}$ depends in addition on the ellipticity bound $\kappa$ for $\sigma$ (but not on $\bar{\sigma}$) from Assumption 1.2(E) and on the dimension $d$.

3) The random constant $\bar{C}_{t,r,\alpha,\xi}$ appearing in the bound on the distributional part in (6.6) is obtained from Proposition 4.1, for which we used Lemma 4.2. This gives a dependence on the time horizon $T$, the bounds $\|b\|_\infty$, $\|\sigma\|_\infty$, $\|\bar{\sigma}\|_\infty$, as well as on the $\beta$-Hölder-seminorms of $(x, \mu) \mapsto \sigma_t(x, \mu)$ and $(t, x, \mu) \mapsto \sigma_t(x, \mu)$. Both are controlled by $\sup_{t \in [0,T]} [\sigma_t]_\beta$ and $[\bar{\sigma}]_\beta$. All these quantities are controlled by Assumption 1.2(H). Finally, the random constant $C_{t,r,\alpha,\xi}$ from Proposition 4.1 appears in (6.14) through its $\mathbb{L}^q$-norm, and this introduces the dependence on $q$.

## 7 DISCUSSION AND CONCLUDING REMARKS

In (1.11), we established that any cluster point $\mathbf{P}^\infty_{x_0}$ of the sequence of empirical measure flow laws $(\mathbf{P}^n_{x_0})_{n \in \mathbb{N}}$ concentrates its mass on flows consisting of absolutely continuous probability measures, more precisely, that

$$\mathbf{P}^\infty_{x_0}\left[\left\{\omega \in \Omega \,:\, \frac{\mathrm{d}\Lambda_t(\omega)}{\mathrm{d}x} \text{ is in } \mathsf{H}^w_r(\mathbb{R}^d) \text{ for almost all } t \in (0, T]\right\}\right] = 1\,. \tag{1.11}$$

We achieved this by a direct study of the $n$-particle systems (1.1). To contrast our approach to the existing literature, and also to discuss its merits and its applications, it is worthwhile to present a different strategy from the literature that can be fruitfully pursued and also leads to a result as in (1.11).

### 7.1 Two approaches using the propagation of chaos

A classical starting point is the propagation of chaos property for McKean–Vlasov systems. We discuss this concept at an informal level and explain how it can be used. The *propagation of chaos* links under appropriate assumptions on the coefficients $b$, $\sigma$, $\bar{\sigma}$ the particle systems (1.1),

$$\begin{cases} X^{i,n}_t = x_0 + \int_0^t b_s(X^{i,n}_s, \mu^n_s)\,\mathrm{d}s + \int_0^t \sigma_s(X^{i,n}_s, \mu^n_s)\,\mathrm{d}B^{i,n}_s + \int_0^t \bar{\sigma}_s(X^{i,n}_s, \mu^n_s)\,\mathrm{d}Z^n_s\,, \\ \mu^n_t = \dfrac{1}{n}\sum_{i=1}^n \delta_{X^{i,n}_t}\,, \end{cases} \tag{7.1}$$

with the McKean–Vlasov SDE (4),

$$\begin{cases} X_t = x_0 + \int_0^t b_s(X_s, U_s)\,\mathrm{d}s + \int_0^t \sigma_s(X_s, U_s)\,\mathrm{d}B_s + \int_0^t \bar{\sigma}_s(X_s, U_s)\,\mathrm{d}Z_s\,, \\ U_t = \mathrm{Law}(X_t \,|\, \mathcal{F}_t^Z)\,, \end{cases} \quad (7.2)$$

by showing that the pair $(X^{1,n}, \mu^n)$ converges to $(X, U)$ in an appropriate sense as $n \to \infty$; see Sznitman [50] for a classical account in the case without common noise, i.e. when $\bar{\sigma} \equiv 0$, Carmona and Delarue [9, Ch. 2.1] for the case with common noise, or the extensive survey in Chaintron and Diez [10, 11] for more on the topic.

A prototypical assumption on the coefficients in (7.1) under which the propagation of chaos holds is Lipschitz-continuity of $(x, \mu) \mapsto b_t(x, \mu), \sigma_t(x, \mu), \bar{\sigma}_t(x, \mu)$, uniformly in $t \in [0, T]$. In this case, (7.2) has a unique strong solution $(X, U)$ and $(\mu^n)_{n \in \mathbb{N}}$ converges to the unique limit $U$. The key point is that from this convergence, one obtains with (7.2) an *explicit limit dynamics* which one can analyze. For that, different techniques are available, and we mention two central examples.

A first strategy to obtain (1.11) from (7.2), already briefly mentioned in the introduction, uses techniques from Malliavin calculus to analyze the flow of conditional marginals $(\mathrm{Law}(X_t | \mathcal{F}_t^Z))_{t \in [0,T]}$. In a recent contribution, Crisan and McMurray [15] assume that $\sigma$ is uniformly elliptic and the coefficients $b, \sigma, \bar{\sigma}$ are smooth in all variables, with smoothness in the measure argument defined in the sense of P.-L. Lions; see Carmona and Delarue [8, Ch. 5] or Cardaliaguet et al. [7, Sec. 2.2]. They then develop a generalization of the Malliavin calculus based on Kusuoka and Stroock [36] to show that $U_t$ possesses a smooth density for all $t > 0$, almost surely. They also establish additional regularity properties for flows of McKean–Vlasov SDEs which we do not consider here. We refer to [15] for a discussion of these results and also for the connection to further works approaching the regularity question for McKean–Vlasov SDEs using Malliavin calculus.

A second strategy to obtain (1.11) from (7.2) uses techniques for SPDEs. This exploits the link of (7.2) to the *stochastic Fokker–Planck equation* which governs the evolution of $t \mapsto U_t$; see Lacker et al. [37, Sec. 1.2]. Indeed, if we write $U_t[\phi] := \mathbf{E}[\phi(X_t) \,|\, \mathcal{F}_t^Z]$ for $\phi \in \mathcal{S}$, then Itô's formula and a version of the stochastic Fubini theorem yield

$$U_t[\phi] = \phi(x_0) + \int_0^t U_s[(L_s \phi)(U_s)]\,\mathrm{d}s + \int_0^t U_s[(R_s \phi)(U_s)]\,\mathrm{d}Z_s\,, \quad (7.3)$$

where $L_s$ and $R_s$ are differential operators given in terms of the coefficients $b_s, \sigma_s, \bar{\sigma}_s$ appearing in (7.2); see e.g. [9, 37] or also Hammersley et al. [26]. The collection of real-valued flows $\mathcal{U} := \{U[\phi] := (U_t[\phi])_{t \in [0,T]} : \phi \in \mathcal{S}\}$ from (7.3) characterizes the measure-valued flow $U$ in (7.2), and this opens the door for a study using for example a symbolic calculus for pseudo-differential operators as in Chaleyat-Maurel et al. [13] and Chaleyat-Maurel and Michel [12]. While those works are developed with applications to filtering equations in mind, the techniques developed there apply with modifications also to (7.3). Indeed, if $x \mapsto b_t(x, \mu), \sigma_t(x, \mu), \bar{\sigma}_t(x, \mu)$ are smooth and $\mu \mapsto b_t(x, \mu), \sigma_t(x, \mu), \bar{\sigma}_t(x, \mu)$ are uniformly Lipschitz, then [13, Thm. 1.2 and estimates in Sec. 2.5] can be used to deduce that $U_t$ again possesses a smooth density. For this particular result, $\sigma$ need not be elliptic, but the initial condition must have a smooth density. The case of elliptic $\sigma$ and an irregular initial condition (like a point mass $\delta_{x_0}$) requires an adaptation of [12]. We refer to those papers for further discussion and the connection to related works in this direction, such as e.g. Pardoux [42].

## 7.2 The approach here, using directly the particle systems

The common feature of both strategies above is that they operate on the explicit dynamics provided by (7.2), which allows us to arrive at (1.11). Moreover the assumptions imposed in either of these strategies imply that the propagation of chaos holds. Under the weaker assumptions we impose in the present work, it is in general not true that cluster points $\mathbf{P}_{x_0}^\infty$ of $(\mathbf{P}_{x_0}^n)_{n\in\mathbb{N}}$ are related to solutions of the associated McKean–Vlasov equation, as the following trivial but nonetheless valid example shows.

**Example 7.1** | Let $d=1$ and $\eta_t := \mathcal{N}(0,t)$ be the centered normal law with variance $t$. Consider $x_0 = 0$, $b_t(\mu) := b_t(x,\mu) := \mathbb{1}_{\{\eta_t\}}(\mu)$, $\sigma \equiv 1$ and $\bar\sigma \equiv 0$. These coefficients clearly satisfy Assumption 1.1 and 1.2. Each empirical measure $\mu^n = \frac{1}{n}\sum_{i=1}^n \delta_{X_t^{i,n}}$ is purely atomic and thus cannot be equal to $\eta_t$ for $t > 0$. The associated $n$-particle system (7.2) thus takes for any $n \in \mathbb{N}$ the form

$$X_t^{i,n} = x_0 + \int_0^t b_s(\mu_s^n)\,ds + B_t^{i,n} = B_t^{i,n} \qquad \text{for all } i \in [n]. \tag{7.4}$$

With the notation in Section 1.3, the Glivenko–Cantelli theorem shows that we have $\Lambda_t = \eta_t$ for all $t \in [0,T]$, $\mathbf{P}_{x_0}^\infty$-a.s.; see Parthasarathy [43, Thm. II.7.1]. Consider now the (formal) McKean–Vlasov SDE with the same coefficients $b$, $\sigma$, $\bar\sigma$. It takes the form

$$X_t = \int_0^t \mathbb{1}_{\{\eta_s\}}(U_s)\,ds + B_t \qquad \text{with } U_t = \mathrm{Law}(X_t). \tag{7.5}$$

One can show that (7.5) does not have a weak solution $(X,U)$, and thus we also cannot have $(X^{i,n}, \mu) \to (X, U)$ in any sense as $n \to \infty$.

To see informally that we cannot have convergence to dynamics satisfying (7.5), note that (7.4) gives $(X^{1,n}, \mu^n) \to (B, \eta)$ in law on $\mathbf{C}([0,T]; \mathbb{R} \times \mathcal{S}'(\mathbb{R}))$ as $n \to \infty$, where $B$ is a Brownian motion and $\eta = (\eta_t)_{t\in[0,T]}$ is the flow of Brownian time-marginals. If the equality $(B,\eta) = (X,U)$ were true and if $(X,U)$ were to satisfy the dynamics in (7.5), then $X = B$ and (7.5) would give $\mathbb{1}_{\{\eta_s\}}(U_s) = 0$ $ds$-a.e., and hence

$$\mathrm{Law}(X_s) = U_s \neq \eta_s = \mathcal{N}(0,s) \qquad ds\text{-a.e.}$$

But since $B$ is a Brownian motion and $X = B$, (7.5) also implies the opposite, namely

$$U_s = \mathrm{Law}(X_s) = \mathrm{Law}(B_s) = \mathcal{N}(0,s) = \eta_s \qquad \text{for all } s \in (0,T].$$

Of course the example is oversimplified in that we know the sequence $(\mu^n)_{n\in\mathbb{N}}$ to converge in law to the unique limit $\eta = (\eta_t)_{t\in[0,T]}$, where $\eta_t = \mathcal{N}(0,t)$ is the time-marginal law of a Brownian motion. In more complicated cases, however, we have in general neither uniqueness of limits nor access to an explicit dynamics.

*The point we want to emphasize is that a regularity result of type* (1.11) *can hold for cluster points even if one does not know the limiting dynamics and thus, in particular, also* *if the propagation of chaos does not hold.* Developing a strategy to obtain (1.11) which circumvents the propagation of chaos and works directly with cluster points $\mathbf{P}_{x_0}^\infty$ of $(\mathbf{P}_{x_0}^n)_{n\in\mathbb{N}}$ is therefore desirable. We are aware only of the works of Varadhan [54, Thm. 3.2(i)] and Olla and Varadhan [41] that have the same philosophy and obtain insights of comparable nature as (1.11). However, these works consider rather specific systems of particles interacting on a circle only through their drifts and in a setting without common noise. Their results are obtained via an entropy technique whose applicability appears limited to certain dynamics and which requires more regularity than we impose in Assumptions 1.1 and 1.2.

Here, we chose instead to use an *approximation and interpolation technique* that has its origins in the seminal work of Fournier and Printems [25]. Let us discuss this technique, and our adaptation of it, by presenting it in the basic case of simple diffusive dynamics.

**Example 7.2** | For simplicity, we take $d = 1$. Let $(\Omega, \mathcal{F}, \mathbf{P})$ be a generic probability space with filtration $\mathbb{F} = (\mathcal{F}_t)_{t\in[0,T]}$, on which we are given a real-valued process satisfying the dynamics
$$S_t = s_0 + \int_0^t a_r(S_r)\,\mathrm{d}r + \int_0^t m_r(S_r)\,\mathrm{d}W_r \qquad \text{for } t \in [0,T],$$
where $a, m : [0,T] \times \mathbb{R} \to \mathbb{R}$ are bounded measurable functions and $W$ is a real-valued $(\mathbb{F}, \mathbf{P})$-Brownian motion. In addition, we assume that

a) for some $\beta > 0$, the function $x \mapsto m_t(x)$ is $\beta$-Hölder-continuous, uniformly in $t \in [0,T]$;

b) for some $\kappa > 0$, we have $m_t(x) \geq \kappa$ for all $(t,x) \in [0,T] \times \mathbb{R}$.

In this setting, the approximation and interpolation approach can be used to show that $\mathrm{Law}_{\mathbf{P}}(S_t)$ possesses a density for all $t > 0$. This is achieved as follows.

i) Choose an appropriate family $\mathcal{T}$ of test functions $\phi : \mathbb{R} \to \mathbb{R}$. Fix $t \in (0,T]$ and consider the terms
$$\lambda_t^\phi := \mathbf{E}[\phi(S_t)] \qquad \text{for } \phi \in \mathcal{T}.$$

ii) For $0 \leq \varepsilon \leq t$, define an approximation $V_{\varepsilon,t} = S_{t-\varepsilon} + \int_{t-\varepsilon}^t m_r(S_{t-\varepsilon})\,\mathrm{d}W_r$ and consider
$$a_{\varepsilon,t}^\phi := \mathbf{E}[\phi(V_{\varepsilon,t})] \quad \text{and} \quad e_{\varepsilon,t}^\phi := \mathbf{E}[\phi(S_t)] - \mathbf{E}[\phi(V_{\varepsilon,t})] \qquad \text{for } \phi \in \mathcal{T}.$$
This clearly gives
$$\lambda_t^\phi = a_{\varepsilon,t}^\phi + e_{\varepsilon,t}^\phi \qquad \text{for all } \phi \in \mathcal{T} \text{ and } 0 \leq \varepsilon \leq t.$$

iii) The specific form of the approximation $V_{\varepsilon,t}$ shows that
$$\mathrm{Law}_{\mathbf{P}}(V_{\varepsilon,t} \mid \mathcal{F}_{t-\varepsilon}) = \mathcal{N}\left(S_{t-\varepsilon}, \int_{t-\varepsilon}^t m_r^2(S_{t-\varepsilon})\,\mathrm{d}r\right) \qquad \mathbf{P}\text{-a.s.}$$

This is used with the tower property, the ellipticity condition b), and specific properties of the test function $\phi$ to obtain a bound of the form
$$|a_{\varepsilon,t}^\phi| \leq \mathbf{E}[|\phi(V_{\varepsilon,t})|] \leq \mathbf{E}\!\left[\left|\mathbf{E}[\phi(V_{\varepsilon,t}) \mid \mathcal{F}_{t-\varepsilon}]\right|\right] \leq c(\phi)p_a(\varepsilon), \tag{7.6}$$
where $c(\phi)$ is some constant depending on $\phi$ and $p_a(\varepsilon)$ is a function depending on $\varepsilon$ with $p_a(\varepsilon) \to \infty$ as $\varepsilon \to 0$.

iv) The specific form of the error $S_t - V_{\varepsilon,t}$ allows us to use first the triangle inequality, the Itô isometry and the $\beta$-Hölder-regularity of $\sigma$ and then standard estimates for SDEs, e.g. in the form of Lemma 2.1, to deduce that

$$\|S_t - V_{\varepsilon,t}\|_{\mathbb{L}^2(\mathbf{P})} \leq \varepsilon c(m) + [\sigma]_\beta \left(\mathbf{E}\left[\int_{t-\varepsilon}^t |S_r - S_{t-\varepsilon}|^{2\beta} \mathrm{d}r\right]\right)^{1/2}$$
$$\leq c(m)\varepsilon + c(T,a,m,\beta)\varepsilon^{(1+\beta)/2}$$
$$\leq c(T,a,m,\beta)\varepsilon^{(1+\beta)/2}.$$

With specific properties of the test function $\phi$, we next obtain a bound of the form

$$|e^\phi_{\varepsilon,t}| \leq c(\phi)\|S_t - V_{\varepsilon,t}\|_{\mathbb{L}^2(\mathbf{P})} \leq c(\phi)c(T,a,m,\beta)p_e(\varepsilon,\beta), \tag{7.7}$$

where $p_e(\varepsilon,\beta)$ is a function with $p_e(\varepsilon,\beta) \to 0$ as $\varepsilon \to 0$, for any $\beta > 0$.

v) Combine the estimates (7.6) and (7.7) for $a^\phi_{\varepsilon,t}$ and $e^\phi_{\varepsilon,t}$ to obtain an estimate for their sum $\lambda^\phi_t$ which is valid for all test functions $\phi \in \mathcal{T}$. If the family $\mathcal{T}$ is separating, use the estimates for $\lambda^\phi_t$ with $\phi \in \mathcal{T}$ to deduce regularity of $\mathrm{Law}_\mathbf{P}(X_t)$ and from this that $\mathrm{Law}_\mathbf{P}(X_t) \ll \mathrm{d}x$ for all $t \in (0,T]$.

In the literature (see below), this simple strategy has been extended to more complicated diffusive dynamics such as e.g. path-dependent cases, and also to certain SPDEs. The strategy thus turns out to be quite versatile and flexible. For step v), the key is to devise a good way of combining the estimates obtained in iii) and iv), and this goes hand-in-hand with choosing in i) a good family of test functions.

In the case of Fournier and Printems [25], test functions are complex exponentials and an $\mathbb{L}^2(\mathrm{d}x)$-regularity result is deduced by developing a hands-on interpolation scheme to control the Fourier transform; see Remark 7.3 below. Subsequent works such as Debussche and Fournier [20], Debussche and Romito [21], Romito [46] use as test functions a family of Besov functions and in v) functional-analytic means to combine the estimates, and Bally [4], Bally and Caramellino [5] use certain families of smooth functions and real interpolation theory to combine the estimates. For a fuller discussion, applications and additional references, we refer to [5, 25] and in particular to the insightful and thorough overview in [46], which also discusses the merits and limitations of the approximation and interpolation technique.

The approximation and interpolation technique requires ellipticity and some regularity for the diffusion coefficient to obtain existence of a density for *every* $t > 0$. These assumptions cannot be significantly weakened as the counterexample in Fabes and Kenig [23] shows. Classical PDEs techniques such as those in Aronson [3], Porper and Èidel'man [45], Fabes and Rivière [24], Stroock and Varadhan [49, Ch. 9] demand less regularity on the diffusion coefficient while giving stronger estimates for the density, however, but only for *almost every* $t > 0$. One of the merits of the approximation and interpolation strategy is next to its simplicity also that it is capable of handling cases that PDE techniques cannot, e.g. path-dependent examples.

Our aim here is not to study a *given dynamics* such as the one in Example 7.2, but rather analyze a *random flow of measures arising as a cluster point* of a sequence of empirical measure flows. Here again the approximation and interpolation strategy turns out to be sufficiently robust to be applied to this case. With Example 7.2 in mind, we now discuss how the basic case considered there naturally paves the way to approach the problem we have in mind, and how the $\mathcal{S}'$-valued setting provides a convenient and systematic framework for this analysis.

i) First we must obtain an object to analyze via the interpolation scheme. Considering as in Example 7.2 the law of a given random variable for a fixed $t \in (0,T]$ is not feasible in our case. Instead, our starting point is for each $n \in \mathbb{N}$ and $\phi \in \mathcal{S}$ the process

$$\mu_t^n[\phi] = \frac{1}{n}\sum_{i=1}^n \delta_{X_t^{i,n}}[\phi] = \frac{1}{n}\sum_{i=1}^n \phi(X_t^{i,n}) \qquad \text{for } t \in [0,T]. \tag{7.8}$$

Since we look in (1.11) for a regularity result not for a fixed $t \in (0,T]$ under some cluster point, but for a result that applies under any fixed cluster point to almost every $t$, we must construct our object of study globally. To this end we first aggregate (7.8) to a sequence of $\mathcal{S}'$-valued processes $\mu^n = (\mu_t^n)_{t \in [0,T]}$ for which we show continuity of paths and then tightness of the sequence $(\mathbf{P}_{x_0}^n)_{n \in \mathbb{N}}$ of laws associated with $(\mu^n)_{n \in \mathbb{N}}$. As we showed in Step 1 of the proof of Theorem 1.5, this identifies through a limiting procedure the $\mathcal{S}'$-valued process $\Lambda$ on $(\Omega, \mathcal{F}, \mathbf{P}_{x_0}^\infty)$ that we can study using $\phi \in \mathcal{S}$ as test functions. In this way, we can fix a cluster point $\mathbf{P}_{x_0}^\infty$, under which all our subsequent considerations are valid for any $t \in (0,T]$ and $\phi \in \mathcal{S}$. So in contrast to Example 7.2, where we considered a fixed measure for each $t \in (0,T]$, we consider here a measure-valued random process $\Lambda$.

ii) As in Example 7.2, we next construct an appropriate decomposition. We do not have access to a decomposition directly (i.e. under $\mathbf{P}_{x_0}^\infty$), but only for each $n \in \mathbb{N}$; so we again must use a limiting procedure. We thus consider the sequence of two-parameter $\mathcal{S}'$-valued processes $((\nu_{\varepsilon,t}^n)_{(\varepsilon,t) \in [0,T]^\Delta})_{n \in \mathbb{N}}$, for which we show continuity to subsequently deduce tightness. Here the two-stage procedure in (2.35) and (2.36) is important: We first choose a cluster point $\mathbf{P}_{x_0}^\infty$ and only afterwards construct the associated $\bar{\mathbf{P}}_{x_0}^\infty$ and the decomposition. With the same reasoning as in i), we also must construct the approximation jointly for all $(\varepsilon,t) \in [0,T]^\Delta$ under $\bar{\mathbf{P}}_{x_0}^\infty$. These considerations lead us in Step 2 of the proof of Theorem 1.5 and bring us to the $\mathcal{S}'$-valued two-parameter processes $A$ and $E$ on $(\bar{\Omega}, \bar{\mathcal{F}}, \bar{\mathbf{P}}_{x_0}^\infty)$.

iii)–iv) The objects obtained in ii) satisfy a priori no dynamics whose structure can be used to deduce estimates as in parts iii) and iv) of Example 7.2 in a parallel way. Having obtained $A$ and $E$ on $(\bar{\Omega}, \bar{\mathcal{F}}, \bar{\mathbf{P}}_{x_0}^\infty)$ through a limiting procedure, we instead use the $n$-particle systems $(X^{i,n})_{i \in [n]}$ in (1.17) and the explicit approximation $(Y^{i,n;\varepsilon,t})_{i \in [n]}$ in (1.18) to adapt the way in which we obtain estimates.

For $A_{\varepsilon,t}$, Proposition 3.1 shows that $A_{\varepsilon,t}$ is a Gaussian mixture that can be represented in terms of $Z_{\cdot \wedge t}$ and $\Lambda_{t-\varepsilon}$. Hence at a fundamental level, $A_{\varepsilon,t}$ shares a similar conditional Gaussian structure as the conditionally Gaussian term $V_{\varepsilon,t}$ from Example 7.2. The origin of this is in both cases the independent Brownian noise, but the way we obtain this structure is different. Proposition 3.1 uses a conditional law of large numbers argument; see the application of Hoeffding's inequality in Step 2 of the proof of Proposition 3.1. The way we use the Gaussian structure is then analogous to Example 7.2, namely via the ellipticity of $\sigma$ which translates with the scaling property of the $\mathsf{H}_r^s$-norm into a quantitative estimate for $A_{\varepsilon,t}[\phi]$ in Proposition 3.2.

The estimate for the error $E_{\varepsilon,t}$ is more nuanced than that in part iv) of Example 7.2. Recall that $\text{Law}_{\bar{\mathbb{P}}^n_{x_0}}(E_{\varepsilon,t}) = \text{Law}_{\mathbb{P}^n_{x_0}}(\mu^n_t - \nu^n_{\varepsilon,t})$ for all $n \in \mathbb{N}$. The specific property of the test functions $\phi \in \mathcal{S}$ we use is $\rho$-Hölder-continuity for arbitrary $\rho \in (0,1)$. As in (4.6), this allows us to obtain in a first step for each $n \in \mathbb{N}$ the bound

$$|\mu^n_t[\phi] - \nu^n_{\varepsilon,t}[\phi]| \leq \frac{1}{n}\sum_{i=1}^{n} |\phi(X^{i,n}_t) - \phi(Y^{i,n;\varepsilon,t}_{\varepsilon,t})| \leq [\phi]_\rho \frac{1}{n}\sum_{i=1}^{n} |X^{i,n}_t - Y^{i,n;\varepsilon,t}_t|^\rho, \qquad (7.9)$$

which we then expand in terms of the dynamics of $X^{i,n}$ and $Y^{i,n;\varepsilon,t}$ from (1.17) and (1.18) to find a bound by random terms of the form

$$|\mu^n_t[\phi] - \nu^n_{\varepsilon,t}[\phi]| \leq [\phi]_\rho \left( \varepsilon^\rho + \left| \frac{1}{n}\sum_{i=1}^{n} \int_{t-\varepsilon}^{t} \left( \sigma_r(X^{i,n}_r, \mu^n_r) - \sigma_{t-\varepsilon}(X^{i,n}_{t-\varepsilon}, \mu^n_{t-\varepsilon}) \right) dZ^n_r \right|^\rho + o_n \right),$$

where $o_n$ denotes a random contribution that stems from the idiosyncratic Brownian motions $(B^{1,n}, \ldots, B^{n,n})$ with an $\mathbb{L}^2(\mathbb{P}^n_{x_0})$-norm no larger than $(c\varepsilon/n)^\rho$ for some global constant $c < \infty$. In contrast to Example 7.2, we thus do not need to control a *norm*, but must instead obtain an *exact pathwise estimate*, uniformly for a family of stochastic processes. Here a Kolmogorov-type argument appears natural. However, as mentioned before proving Proposition 4.1, a standard argument does not give a bound of order $\rho(1+\beta)/2$ as we need, but only a bound of order $\rho/2$ which is too poor. We thus used for Proposition 4.1 a direct estimate inspired by the Kolmogorov argument to obtain the more powerful local bound for $|\mu^n_t[\phi] - \nu^n_{\varepsilon,t}[\phi]|$; cf. (4.11). This bound is unfortunately not valid for all $\varepsilon \in (0,t]$, but only for $\varepsilon$ in the countable set $(0,t]_\alpha \subsetneq (0,t]$ for $\alpha \in (0,1)$.

v) The final step is the combination of the estimates obtained from iii) and iv). In contrast to Example 7.2, we have estimates only for certain pairs $(\varepsilon, t) \in [0,T]^\Delta$. In addition, the estimates we obtain are random and may also fail on nullsets which can depend on $(\varepsilon, t)$; see Propositions 3.2 and 4.1. These added technical intricacies restrict us in the way we can combine the estimates, and also demand additional care to ensure measurability and integrability of the estimate obtained from the combination.

We resolved this by using part 3) of Proposition 5.1. This result is specifically tailored to our needs. It makes use only of the estimates we have access to and combines them in a way that respects measurability; cf. (6.7). At the same time, Proposition 5.1 is transparent about how constants evolve in the interpolation. This clarifies how the random constant $\bar{C}_{t,r,\alpha,\xi}$ in (4.1) from the bound on the distributional part transforms, allowing us to obtain the $\mathbb{L}^q(\mathbf{P}^\infty_{x_0})$-moment-bound on the random $\mathsf{H}^w_r(\mathbb{R}^d)$-norm of the density in Theorem 1.5; see (5.3) and (6.8).

For our study, we chose the $\mathcal{S}'$-valued framework. Stochastic calculus and interpolation theory set in $\mathcal{S}'(\mathbb{R}^d)$ have a long-established history; see Itô [30] and Bergh and Löfström [6, Ch. 5]. We feel that our development is most natural in this framework. It allows us to a priori carry out the analysis in $\mathcal{S}'(\mathbb{R}^d)$ and restrict to a specific subspace $\mathsf{H}^s_r(\mathbb{R}^d)$ a posteriori whenever this is convenient or technically needed. This framework also makes for an efficient presentation of Proposition 5.1, which is formulated at an abstract level and yet in an elementary way. The steps we followed in the proof of that result suggest themselves naturally within the standard interpolation framework that is a well known and standard tool in harmonic analysis. Our hope is that this argument sheds some light on the fundamental workings of the approximation and interpolation approach, which compels in its original form in Fournier and Printems [25] through its striking simplicity.

We should mention that Bally [4], Bally and Caramellino [5] also use an interpolation argument. While their argument requires less regularity for the diffusion coefficient ($m$ in Example 7.2), this comes at the expense of a substantial amount of additional technicalities that we wish to avoid here. The approach in Debussche and Fournier [20], Debussche and Romito [21], Romito [46] set in the three-parameter Besov family gives another result that is comparable to part 2) of Proposition 5.1. However, part 2) of Proposition 5.1 cannot be applied to prove Theorem 1.5. In addition, we saw following (5.22) that to obtain Proposition 5.1, a third parameter can be avoided, i.e. that the two-parameter family of Bessel potential spaces is sufficient for our purposes.

**Remark 7.3** | Let us also briefly compare Proposition 5.1 to the original approach in Fournier and Printems [25]. In their seminal contribution, [25] showed for the time-marginals $\lambda_t := \text{Law}(S_t)$ of the diffusive dynamics $S = (S_t)_{t \in [0,T]}$ from Example 7.2 absolute continuity and that $d\lambda_t/dx$ is in $\mathbb{L}^2(dx)$ for all $t > 0$; see [25, Thm. 1.2]. The approach taken in [25] is in essence a bare-handed interpolation argument. In our Proposition 5.1, square-integrability corresponds to the case $r = 2$ and $0 = w < w_0$. The argument in [25] requires a one-dimensional setting; so we take $d = 1$ in Proposition 5.1. With these choices, we get in (5.24) the condition

$$(1 + \xi)s > (1/2 + 1)(s + 1/2),$$

which is equivalent to $(\xi - 1/2)s > 1/2$. For this to be valid for some $s > 0$, we must have $\xi > 1/2$. With $\beta = \xi$ in Example 7.2, this corresponds to the $\beta$-Hölder-regularity of the diffusion coefficient $m$ of $S$ in Example 7.2 being better than $1/2$, which is precisely what is assumed in [25].

## 7.3   Applications and future directions

To conclude, let us mention an application of Theorem 1.5, and possible future directions.

We already mentioned in the introduction that Theorem 1.5 is a key ingredient in forthcoming work, where we use it as follows. Starting with the sequence $(\mu^n)_{n \in \mathbb{N}}$ of empirical measure flows in (7.1), and using that the bound (1.13) in Theorem 1.5 is valid for any cluster point of $(\mathbf{P}^n_{x_0})_{n \in \mathbb{N}}$, we can show for of coefficients having only low regularity the existence of weak solutions and the validity of an appropriate version of the propagation of chaos for the McKean–Vlasov SDE (7.2). Our strategy thus obtains the *propagation of chaos from the emergence of regularity*. This reverses the direction of reasoning discussed in Section 7.1 for the regime of smooth coefficients, where the strategy is to start with a solution of (7.2), obtain from it a regularity result for $d\text{Law}(X_t \mid \mathcal{F}^Z_t)/dx$ in the flavor of (1.13) and identify this as a property satisfied by the limit of $(\mu^n)_{n \in \mathbb{N}}$ in (7.1) by appealing to the propagation of chaos.

It may be possible to extend our approach to more general McKean–Vlasov dynamics, e.g. with jumps, or even to other particle systems. It also appears within reach to show that any solution of the McKean–Vlasov SDE (7.2) must satisfy the conclusions of Theorem 1.5. Let us also note that in principle, nothing in our arguments prevents us from allowing an explicit dependence of the coefficients $b$, $\sigma$, $\bar{\sigma}$ on $n$ as long as Assumptions 1.1 and 1.2 are satisfied uniformly in $n$.

A good discussion of the limitations of the approximation and interpolation scheme can be found in [46]. That discussion shows that the boundedness assumption on the drift $b$ can in principle be weakened, but this extension is technical. It may also be possible to use a different proof technique for the representation of the approximation $A$ in Proposition 3.2 that does not necessitate the $\beta$-Hölder-regularity in the time variable for $\bar{\sigma}$ from Assumption 1.2(H). In fact, we imposed this assumption to get continuity of the representation in the common noise $Z$; see (3.16). This avoids technicalities, but leads to an asymmetry in the auxiliary process $Y^{i,n;\varepsilon,t}$ from (1.18) between the common and the idiosyncratic noise.

While our result also applies to the case without common noise, i.e. when $\bar{\sigma} \equiv 0$, we expect that classical results on the existence of time-marginal densities for SDEs such as [3, 24, 45, 49] in the discussion following Example 7.2 can be used to get better bounds than the one we found in (1.11). Finally, it may be possible to adapt some of these classical PDE techniques also to the case with common noise.

# A BACKGROUND

## A.1 Functional analysis

**The spaces of Schwartz functions an tempered distributions.** Let $\alpha \in \mathbb{N}_0^d$ be a multi-index, define $|\alpha| = \alpha_1 + \cdots + \alpha_d$, write $\partial^\alpha := \partial_1^{\alpha_1} \cdots \partial_d^{\alpha_d}$ and set

$$\langle x \rangle := (1 + |x|^2)^{1/2} \quad \text{for } x \in \mathbb{R}^d. \tag{A.1}$$

For each $m \in \mathbb{N}_0$ and $\phi \in \mathbf{C}^\infty(\mathbb{R}^d)$, we define the seminorm

$$\|\phi\|_m^* := \max_{|\alpha| \leq m} \sup_{x \in \mathbb{R}^d} |\langle x \rangle^m \partial^\alpha \phi(x)|. \tag{A.2}$$

The *space of rapidly decreasing functions* or Schwartz functions is the vector space

$$\mathcal{S} := \mathcal{S}(\mathbb{R}^d) := \{\phi \in \mathbf{C}^\infty(\mathbb{R}^d) : \|\phi\|_m^* < \infty \text{ for all } m \in \mathbb{N}_0\}.$$

We endow $\mathcal{S}$ with the Fréchet topology generated by the family $\|\cdot\|_m^*$ for $m \in \mathbb{N}_0$, i.e. $\phi_n \to \phi$ in $\mathcal{S}$ if and only if $\|\phi_n - \phi\|_m^* \to 0$ when $n \to \infty$, for each $m \in \mathbb{N}$. As usual, a set $B \subseteq \mathcal{S}$ is said to be *bounded* if $\sup_{\phi \in B} \|\phi\|_m^* < \infty$ for all $m \in \mathbb{N}_0$.

It is well known that $\mathcal{S}$ is metrizable, complete and separable. Indeed, let $\mathbf{C}_c^\infty = \mathbf{C}_c^\infty(\mathbb{R}^d)$ be the Fréchet space of smooth and compactly supported functions $\mathbb{R}^d \to \mathbb{R}$; see Rudin [47, Sec. 6.2] for details. Then $\mathbf{C}_c^\infty \hookrightarrow \mathcal{S}$, which is to say that $\mathbf{C}_c^\infty$ embeds continuously into $\mathcal{S}$. Moreover, this embedding is dense. As a consequence we also have the continuous embedding $\mathcal{S} \hookrightarrow \mathbb{L}^r(\mathrm{d}x) = \mathbb{L}^r(\mathbb{R}^d, \mathrm{d}x)$ for all $r \in [1, \infty]$, which is dense if $r \in [1, \infty)$.

Let $\lambda : \mathcal{S} \to \mathbb{R}$ be a linear functional. We denote the duality pairing by $\langle \lambda; \phi \rangle_{\mathcal{S}' \times \mathcal{S}}$ for $\phi \in \mathcal{S}$, and usually write $\lambda[\phi] := \langle \lambda; \phi \rangle_{\mathcal{S}' \times \mathcal{S}}$. The functional $\lambda$ is *continuous* if $\phi_n \to \phi$ in $\mathcal{S}$ implies $\lambda[\phi_n] \to \lambda[\phi]$. Equivalently, there exist $C > 0$ and $m \in \mathbb{N}$ such that

$$|\lambda[\phi]| \leq C \|\phi\|_m^* \quad \text{for all } \phi \in \mathcal{S}(\mathbb{R}^d). \tag{A.3}$$

The minimal index $m$ validating this inequality is called the *order* of $\lambda$. Continuous linear functionals are called *(tempered) distributions*; they form the strong topological dual space

$$\mathcal{S}' := \mathcal{S}'(\mathbb{R}^d) := \{\lambda : \mathcal{S}(\mathbb{R}^d) \to \mathbb{R} \ : \ \lambda \text{ is continuous}\}.$$

Its strong topology is generated by uniform convergence on bounded sets, i.e. the topology is generated the family of seminorms $\|\lambda\|_{\mathcal{S}',B} := \sup_{\phi \in B} |\lambda[\phi]|$ as $B$ ranges over the bounded sets of $\mathcal{S}$.

The strong topology generates the Borel-$\sigma$-field on $\mathcal{S}'$, denoted $\mathcal{B}(\mathcal{S}'(\mathbb{R}^d))$ or simply $\mathcal{B}(\mathcal{S}')$. It can be shown that $\mathcal{B}(\mathcal{S}')$ equals the $\sigma$-algebra generated by the collection of sets $\{\lambda \in \mathcal{S}' : \lambda[\phi] < a\}$ for $\phi \in \mathcal{S}$ and $a \in \mathbb{R}$; see e.g. Kallianpur and Xiong [32, Thm. 3.1.1].

Let $f$ be a measurable function. We define

$$\lambda_f[\phi] := \int_{\mathbb{R}^d} \phi(x) f(x) \, \mathrm{d}x \,, \tag{A.4}$$

for $\phi \in \mathcal{S}$, provided the integral exists. Consider $\lambda \in \mathcal{S}'$ and suppose that there exists a $f \in \mathbb{L}^1_{\mathrm{loc}}(\mathrm{d}x)$ such that $\lambda = \lambda_f$ in $\mathcal{S}'$, that is to say $\lambda[\phi] = \lambda_f[\phi]$ for all $\phi \in \mathcal{S}$. Then $f$ is determined uniquely up to $\mathrm{d}x$-a.e. equivalence; see Hörmander [27, Thm. 1.2.5]. In this case it is customary to say that $\lambda$ *is a function*, or that $\lambda$ is representable by $f$. The *canonical embedding* $\mathcal{S} \hookrightarrow \mathcal{S}'$ mapping $f \in \mathcal{S}$ to $\lambda_f \in \mathcal{S}'$ is thus well-defined and, in fact, continuous. This inclusion allows us to identify $\mathcal{S}$ with a subset of $\mathcal{S}'$: If $\lambda \in \mathcal{S}'$ is represented by a function in $\mathcal{S}$, then we freely write $\lambda \in \mathcal{S}$. Similarly, we have the canonical embedding $\mathbb{L}^r(\mathrm{d}x) \hookrightarrow \mathcal{S}'$ for all $r \in [1, \infty]$, and if $\lambda \in \mathcal{S}'$ is represented by a function in $\mathbb{L}^r(\mathrm{d}x)$, we simply write $\lambda \in \mathbb{L}^r(\mathrm{d}x)$. In the same fashion, if $\mu$ is a measure on the Borel sets of $\mathbb{R}^d$, we may define

$$\lambda_\mu[\phi] := \int_{\mathbb{R}^d} \phi(x) \, \mu(\mathrm{d}x) \tag{A.5}$$

for $\phi \in \mathcal{S}$, provided the integral exists. In analogy to the above, we say that $\lambda \in \mathcal{S}'$ is a measure if there is a measure $\mu$ such that $\lambda = \lambda_\mu$ in $\mathcal{S}'$. Let $\mathbf{M}_1^+ = \mathbf{M}_1^+(\mathbb{R}^d)$ denote the set of probability measures on the Borel sets of $\mathbb{R}^d$. It is standard to verify that $\mathbf{M}_1^+ \subseteq \mathcal{S}'$ by virtue of (A.5). In fact, the inclusion is sequentially continuous if $\mathbf{M}_1^+$ carries the narrow topology; see Lemma A.3 below.

**Bessel potential spaces** We base the following discussion on Bergh and Löfström [6, Ch. 6]. The Fourier transform $\mathscr{F}$ is defined as the map that sends $f \in \mathcal{S}$ to $\mathscr{F}f : \mathbb{R}^d \to \mathbb{R}$ given by

$$\mathscr{F}f(\xi) := \int_{\mathbb{R}^d} f(x) \exp(-ix \cdot \xi) \, \mathrm{d}x \quad \text{for } \xi \in \mathbb{R}^d \,. \tag{A.6}$$

It can be shown that $\mathscr{F}$ is an $\mathcal{S}$-valued continuous linear map that admits an inverse $\mathscr{F}^{-1}$. In other words, $\mathscr{F} : \mathcal{S} \to \mathcal{S}$ is a linear automorphism in $\mathcal{S}$.

The Fourier transform, initially defined on $\mathcal{S}$, extends to a linear automorphism in $\mathcal{S}'$ by duality. It is again denoted by $\mathscr{F}$. More explicitly, $\mathscr{F} : \mathcal{S}' \to \mathcal{S}'$ is uniquely defined by the relation

$$\langle \mathscr{F}g; f \rangle_{\mathcal{S}' \times \mathcal{S}} := \langle g, \mathscr{F}f \rangle_{\mathcal{S}' \times \mathcal{S}} \quad \text{with } g \in \mathcal{S}' \text{ and } f \in \mathcal{S},$$

in which the left-hand side is defined by the right-hand side which has meaning since $\mathscr{F}$ in (A.6) is a map $\mathcal{S} \to \mathcal{S}$. Next, define the *Bessel potential*

$$J^s f = \mathscr{F}^{-1}\big(h^s(\mathscr{F} f)\big)$$

with $s \in \mathbb{R}$ and $f \in \mathcal{S}'$, and where $h(\xi) := \langle \xi \rangle = (1 + |\xi|^2)^{1/2}$. Let us then set

$$\|f\|_{\mathsf{H}^s_r(\mathbb{R}^d)} := \|J^s f\|_{\mathbb{L}^r(\mathrm{d}x)} \tag{A.7}$$

and define the *Bessel potential spaces* for $s \in \mathbb{R}$ and $1 < r < \infty$ by

$$\mathsf{H}^s_r(\mathbb{R}^d) := \{f \in \mathcal{S}'(\mathbb{R}^d) : \|f\|_{\mathsf{H}^s_r(\mathbb{R}^d)} < \infty\}.$$

From the fact that $\mathscr{F}$ is an automorphism and $\mathcal{S} \subseteq \mathbb{L}^r$ is dense, we can see that $\mathsf{H}^0_r(\mathbb{R}^d) = \mathbb{L}^r(\mathbb{R}^d, \mathrm{d}x)$. In addition, we have the continuous embeddings

$$\mathsf{H}^{s_1}_r(\mathbb{R}^d) \hookrightarrow \mathsf{H}^{s_0}_r(\mathbb{R}^d) \tag{A.8}$$

for any $s_1 > s_0$ and

$$\mathcal{S}(\mathbb{R}^d) \hookrightarrow \mathsf{H}^s_r(\mathbb{R}^d) \hookrightarrow \mathbb{L}^r(\mathbb{R}^d) \hookrightarrow \mathsf{H}^{-s}_r(\mathbb{R}^d) \hookrightarrow \mathcal{S}'(\mathbb{R}^d) \tag{A.9}$$

for any $s > 0$. When there is no risk of confusion, we simply write $\mathsf{H}^s_r$ instead of $\mathsf{H}^s_r(\mathbb{R}^d)$.

If $s \in \mathbb{N}$, it is a classical fact that there are universal constants $\underline{c}$ and $\overline{c}$, both depending on $s$, $r$ and $d$, such that

$$\underline{c}\bigg(\|f\|_{\mathbb{L}^r} + \sum_{i=1}^d \|\partial_i^s f\|_{\mathbb{L}^r}\bigg) \leq \|f\|_{\mathsf{H}^s_r} \leq \overline{c}\bigg(\|f\|_{\mathbb{L}^r} + \sum_{i=1}^d \|\partial_i^s f\|_{\mathbb{L}^r}\bigg), \tag{A.10}$$

where derivatives are understood in the distributional sense; see [6, Thm. 6.2.3].

**Analysis in Banach spaces** The following material is taken from Hytönen et al. [29, Ch. 1.1.a, 1.1.b]. Let $(S, \mathscr{S}, \mu)$ be a measure space, $(Y, \|\cdot\|_Y)$ be a Banach space and $Y'$ its topological dual space. A function $f : S \to Y$ is called *weakly measurable* if $s \mapsto \langle f, y^* \rangle(s) := \langle f(s), y^* \rangle$ is measurable for all $y^* \in Y^*$, *measurable* if $f^{-1}(B) \in \mathscr{S}$ for all Borel-sets $B$ in $Y$, and *strongly measurable* if there exists a sequence of simple functions $f_n : S \to Y$ such that $\lim_{n \to \infty} f_n = f$ pointwise in $S$.

In non-separable Banach spaces these notions of measurability are increasingly strict. Pettis' measurability theorem, which we give next, shows that in the following case the various notions coincide.

A function $f : S \to Y$ is said to be *$\mu$-essentially separably valued* if there exists a closed separable subspace $Y_0$ of $Y$ such that $f(s) \in Y$ for $\mu$-almost all $s \in S$. Note that if $Y$ is a separable Banach space, then $f : S \to Y$ is trivially $\mu$-essentially separably valued.

**Theorem A.1** | *For a function $f : S \to Y$ the following are equivalent:*
  1) *$f$ is strongly measurable;*
  2) *$f$ is $\mu$-essentially separably valued and there exists a weak* dense subspace $X$ of $Y^*$ such that $\langle f, x^* \rangle$ is $\mu$-measurable for all $x^* \in X$.*

This result can be found in [29, Thm. 1.1.20] or Dunford and Schwartz [22, Thm. III.6.11].

## A.2 Probability theory

**Spaces of probability measures.** A special role is our study is taken by positive measures, and in particular probability measures. We thus take a moment to clarify how they relate to our exposition so far. For a topological space X, we write $\mathbf{M}_1^+(X)$ for the set of probability measures on $\mathcal{B}(X)$. If $X = \mathbb{R}^d$, then we oftentimes write $\mathbf{M}_1^+ := \mathbf{M}_1^+(\mathbb{R}^d)$.

The *narrow topology* $\tau_{\mathrm{wk}^*}$ is induced by duality with continuous bounded functions $\mathbf{C}_b(X)$. Specifically, if $(\lambda_m)_{m \in \mathbb{N}}$ is a sequence in $\mathbf{M}_1^+(X)$, then $\lambda_m \to \lambda$ as $m \to \infty$ relative to $\tau_{\mathrm{wk}^*}$ if $\lambda_m[f] \to \lambda[f]$ for all $f \in \mathbf{C}_b(X)$. We write $\mathcal{P}_{\mathrm{wk}^*}(X) := (\mathbf{M}_1^+(X), \tau_{\mathrm{wk}^*})$ for the space $\mathbf{M}_1^+(X)$ endowed with the narrow topology $\tau_{\mathrm{wk}^*}(X)$. If X is Polish, then so is $\mathcal{P}_{\mathrm{wk}^*}(X)$; see, e.g. Aliprantis and Border [2, Thm. 15.15].

Since $\mathbf{M}_1^+ \subseteq \mathcal{S}'$, we can consider the trace topology of $\mathcal{S}'$ on $\mathbf{M}_1^+$. This topology is, however, different from narrow convergence. To discuss this difference, denote by $\mathbf{M}_{\leq 1}^+ = \mathbf{M}_{\leq 1}^+(\mathbb{R}^d)$ the space of subprobability measures on $\mathcal{B}(\mathbb{R}^d)$, i.e. positive measures of no more than unit mass. Evidently $\mathbf{M}_1^+ \subseteq \mathbf{M}_{\leq 1}^+ \subseteq \mathcal{S}'$.

**Remark A.2** | The set $\mathbf{M}_1^+$ is *not closed* in $\mathcal{P}_{\leq 1, \mathcal{S}'}$. Indeed, if $\delta_n$ denotes the Dirac point mass at $(n, \ldots, n) \in \mathbb{R}^d$, and $\mathbf{0}$ the zero measure on $\mathbb{R}^d$, then $\delta_n \to \mathbf{0}$ in $\mathcal{P}_{\leq 1, \mathcal{S}'}$ since $\delta_n[\phi] \to 0$ as $n \to \infty$ for any $\phi \in \mathcal{S}$. On the other hand, $\mathbf{M}_{\leq 1}^+$ is a closed subset of $\mathcal{S}'$. Therefore non-negativity of measures is preserved under convergence in $\mathcal{P}_{\leq 1, \mathcal{S}'}$, but mass may be lost.

The set $\mathbf{M}_{\leq 1}^+$ equipped with the trace topology of $\mathcal{S}'$ is denoted by $\mathcal{P}_{\leq 1, \mathcal{S}'} := (\mathbf{M}_{\leq 1}^+, \tau_{\mathcal{S}'})$. Specifically, if $(\lambda_m)_{m \in \mathbb{N}}$ is a sequence in $\mathbf{M}_{\leq 1}^+$, then $\lambda_m \to \lambda$ as $m \to \infty$ in $\mathcal{P}_{\leq 1, \mathcal{S}'}$ if for all $\phi \in \mathcal{S}$, we have that $\lambda_m[\phi] \to \lambda[\phi]$.

**Lemma A.3** | *If $(\nu_m)_{m \in \mathbb{N}}$ is a sequence in $\mathbf{M}_1^+$, then:*
 1) *If $\nu_m \to \nu$ in $\mathcal{P}_{\mathrm{wk}^*}$, then $\nu \in \mathbf{M}_1^+$ and $\nu_m \to \nu$ in $\mathcal{P}_{\mathcal{S}'}$.*
 2) *If $\nu_m \to \nu$ in $\mathcal{P}_{\leq 1, \mathcal{S}'}$ and $\nu[\mathbb{R}^d] = 1$, then $\nu_m \to \nu$ in $\mathcal{P}_{\mathrm{wk}^*}$.*

*Proof* Part 1) follows from the fact that $\mathcal{S} \subseteq \mathbf{C}_b$ and that the notions of weak and strong convergence in $\mathcal{S}'$ coincide for sequences; see Huang and Yan [28, Thm. 3.12]. Part 2) follows from classical properties of narrow convergence; see Klenke [34, Thm. 13.16]. □

For $\mu, \nu \in \mathbf{M}_1^+$ having first moments, i.e. $\int_{\mathbb{R}^d} |x| \mathrm{d}\mu, \int_{\mathbb{R}^d} |x| \, \mathrm{d}\nu < \infty$, denote by

$$\mathrm{d}_{\mathbb{W}_1}(\mu, \nu) := \sup \left\{ \int f \, \mathrm{d}(\mu - \nu) : \|f\|_{\mathsf{Lip}} \leq 1 \right\}, \tag{A.11}$$

the *Kantorovich–Rubinstein* or 1-Wasserstein metric relative to the Euclidean distance on $\mathbb{R}^d$. We write $\mathbb{W}_1 := (\mathbf{M}_1^+, \mathrm{d}_{\mathbb{W}_1})$ for the associated metric space of probability measures. Details on this distance are found in, e.g., Villani [55, Thm. 7.3, Eq. 7.1].

# B OMITTED PROOFS

*Proof of Lemma 2.4* Combining Triebel [52, Sec. 2.5.6, Eq. (2)] with [52, Sec. 2.3.2, Eq. (9)], we have $\mathsf{H}_{r'}^u \hookrightarrow \mathsf{B}_{r', \max\{r', 2\}}^u$, where $\mathsf{B}_{r', \max\{r', 2\}}^u$ is a space in the so-called Besov family.

Now from [52, Sec. 2.3.2, Eq. (7)], we get for any $\varepsilon > 0$ that $\mathsf{B}^u_{r',\max\{r',2\}} \hookrightarrow \mathsf{B}^{u-\varepsilon}_{r',q}$ for arbitrary $q \geq 1$. But by [52, Sec. 2.7.1, Eq. (12)], we find $\mathsf{B}^{u-\varepsilon}_{r',q} \hookrightarrow \text{Höl}_{1-\varepsilon}$ if $u = 1 + d/r'$. This chain of continuous inclusions gives $\mathsf{H}^u_{r'} \hookrightarrow \text{Höl}_{1-\varepsilon}$, and choosing $\varepsilon$ sufficiently small concludes the proof. □

***Proof of Lemma 3.3*** Part 1) follows from the correspondence of translation and modulation under the Fourier transform, $(\mathscr{F}(\tau_m f))(\xi) = e^{-i\xi \cdot m}(\mathscr{F}f)(\xi)$; see e.g. Strichartz [48, Sec. 3.1(1)]. From the definition of the Bessel potential and the formula for the inverse Fourier transform $(\mathscr{F}^{-1}g)(x) = (2\pi)^{-d}\int_{\mathbb{R}^d} g(\xi)\exp(ix \cdot \xi)\,\mathrm{d}\xi$ for $g \in \mathcal{S}$, it can be shown that $J^s(\tau_m f) = \tau_{-m}(J^s f)$. Then the claim follows from the translation-invariance of the $\mathbb{L}^r(\mathrm{d}x)$-norm. For part 2) refer to the monograph of Triebel [52, Thm. 2.5.6(i) and Prop.1(ii) in Sec. 3.4.1]. □

***Proof of Lemma 3.4*** We start by observing three basic facts from linear algebra.

**Step 1** First, since $V \in \mathcal{V}$ is symmetric, we can write $V = PDP^{-1}$, where $D$ is a diagonal matrix consisting of the eigenvalues $a_1, \ldots, a_d$ of $V$ and $P$ is a unitary matrix of orthogonal eigenvectors so that in particular $P^\top = P^{-1}$. Let $e_i$ be the $i$-th standard basis vector of $\mathbb{R}^d$ and set $y = Pe_i$. Using the ellipticity of $A$ gives

$$c = c\,y^\top y \leq y^\top V y = (P^{-1}y)^\top D P^{-1} y = e_i^\top D e_i = a_i\,.$$

Since $i$ was arbitrary, we find

$$\underline{a} := \min\{a_1, \ldots, a_d\} \geq c\,, \tag{B.1}$$

and so the smallest eigenvalues are all bounded from below by $c$.

Secondly, since $P$ is unitary and invertible, we have $z = Px$ if and only if $x = P^{-1}z$, and then $|x| = |z|$. Denote the matrix operator norm by $\|V^{-1}\| := \sup_{|x|=1}|V^{-1}x|$. Using (B.1), we have $\|V^{-1}\| = \sup_{|z|=1}|D^{-1}z| \leq a_1^{-1} + \cdots + a_d^{-1} \leq dc^{-1}$. Since all norms on $\mathbb{R}^{d\times d}$ are equivalent, we get that all entries of the inverse of $V$ are bounded above, i.e.

$$\max_{i,j\in[d]}(V^{-1})_{i,j} \leq c_d c^{-1}\,, \tag{B.2}$$

for some constant $c_d$ depending only on the dimension $d$.

Thirdly and lastly, let $v_i$ be an eigenvector of unit length corresponding to $a_i$. Since $Vv_i = a_i v_i$ and $v_i$ has unit length, it follows that $\|V\| = \sup_{|x|=1}|Vx| \geq a_i$. Since all norms on $\mathbb{R}^{d\times d}$ are equivalent and $\mathcal{V}$ is assumed to be bounded, there exists $c_0 < \infty$ such that $\|V\| \leq c_0$ for all $V \in \mathcal{V}$. Thus,

$$\bar{a} := \max\{a_1, \ldots, a_d\} \leq c_0\,,$$

and so the largest eigenvalues are all bounded from above. Since

$$x^\top V^{-1} x = (P^{-1}x)^\top D^{-1}(P^{-1}x) \geq c_0^{-1}(P^{-1}x)^\top(P^{-1}x)$$

for all $x \in \mathbb{R}^d$ and $(P^{-1})^\top = P$, we have for all $V \in \mathcal{V}$ that

$$\exp\left(-\frac{1}{2} x^\top V^{-1} x\right) \leq \exp\left(-\frac{1}{2} c_0^{-1} x^\top x\right). \tag{B.3}$$

**Step 2**  We are now ready to prove (3.25). Let $s_0$ be the minimal natural number such that $s < s_0$. Since (A.8) gives $\|g(0,V)\|_{\mathsf{H}_r^s} \leq \|g(0,V)\|_{\mathsf{H}_r^{s_0}}$ for all $V \in \mathcal{V}$, we might as well control the norm on the right-hand side. This can be done transparently using the bound (A.10) on the $\mathsf{H}_r^{s_0}$-norm for integer values of $s_0$. Recall that $\partial_i^{s_0}$ denotes the $s_0$-fold partial derivative with respect to the $i$-th coordinate. Then by (A.10),

$$\|g(\,\cdot\,;0,V)\|_{\mathsf{H}_r^{s_0}} \leq \bar{c}\|g(\,\cdot\,;0,V)\|_{\mathbb{L}^r} + \bar{c}\sum_{i=1}^d \|\partial_i^{s_0} g(\,\cdot\,;0,V)\|_{\mathbb{L}^r}. \tag{B.4}$$

To bound the right-hand side, we use the explicit form of the normal density as $g(y;0,V) = (2\pi)^{-d/2} \det(V)^{-1/2} \exp(-y^\top V^{-1} y/2)$ and the fact that $\det(V) = a_1 \cdots a_d$. The lower bound (B.1) on the eigenvalues and the upper bound (B.3) on the exponential term then show for all $y \in \mathbb{R}^d$ and $V \in \mathcal{V}$ that

$$|g(y;0,V)| \leq (2\pi c)^{-d/2} \exp\left(-(y^\top y)/(2c_0)\right),$$

and $\sup_{V \in \mathcal{V}} \|g(\,\cdot\,;0,V)\|_{\mathbb{L}^r} < \infty$ follows immediately. Moreover, an explicit computation shows that for each $i \in [d]$, we have $\partial_i^{s_0} g(y;0,V) = p_i(y,V^{-1}) g(y;0,V)$, for some polynomial $p_i(y,V^{-1})$ of order $s_0$ in the variables $y$ and $V^{-1}$. Since all entries of $V^{-1}$ are bounded uniformly for $V \in \mathcal{V}$ in view of (B.2), we see that also $\sup_{V \in \mathcal{V}} \|\partial_i^{s_0} g(\,\cdot\,;0,V)\|_{\mathbb{L}^r} < \infty$. In summary, (B.4) is bounded uniformly for $V \in \mathcal{V}$, establishing (3.25) and completing the proof. □

**Proof of Lemma 5.2**  On the one hand, by combining Triebel [52, Thms. 2.5.6(i) and 2.4.2(ii)], we get that $(\mathsf{H}_r^{-u}, \mathsf{H}_r^s)_{\theta,p;K} = \mathsf{B}_{r,p}^{w_0}$, where $\mathsf{B}_{r,p}^{w_0}$ is a space in the Besov family and $w_0 = \theta s - (1-\theta) u$. On the other hand, [52, Sec. 2.3.2, Prop. 2(ii) and 2(iii)] shows that $\mathsf{B}_{r,p}^{w_0} \hookrightarrow \mathsf{B}_{r,r}^w \hookrightarrow \mathsf{H}_r^w$. So in sum, $(\mathsf{H}_r^{-u}, \mathsf{H}_r^s)_{\theta,p;K} = \mathsf{B}_{r,p}^{w_0} \hookrightarrow \mathsf{B}_{r,r}^w \hookrightarrow \mathsf{H}_r^w$, giving the claimed continuous embedding. □